\documentclass[11pt]{amsart}

\usepackage{amsaddr}

\usepackage{amsmath,amsfonts,euscript,amssymb,amsthm,verbatim}
\usepackage{mathrsfs}
\usepackage{bm,bbm}
\usepackage[active]{srcltx}
\usepackage[colorlinks,linkcolor={blue},citecolor={blue},urlcolor={black}]{hyperref}
\usepackage{parskip}
\usepackage[shortlabels]{enumitem}
\usepackage[normalem]{ulem}

\allowdisplaybreaks

\newtheorem{definition}{Definition}[section] 
\newtheorem{theorem}{Theorem}[section]
\newtheorem{corollary}[theorem]{Corollary}

\newtheorem{lemma}[theorem]{Lemma}
\newtheorem{remark}[theorem]{Remark}
\newtheorem{proposition}[theorem]{Proposition}
\newtheorem{example}[theorem]{Example}

\newtheorem{hypothesis}[theorem]{Hypothesis}

\def\St{\mathbb{S}^2}
\def\R{{\mathbb R}}
\def\N{{\mathbb N}}
\def\E{{\mathbb E}}
\def\P{\mathbb P}
\def\T{\mathbb T}

\def\mc{\mathscr}

\def\b{\mathbb}
\def\c{\mathcal}
\def\eps{\varepsilon}
\def\Y{\mc Y}
\def\m{\boldsymbol{m}}
\def\u{\boldsymbol{u}}
\def\v{\boldsymbol{v}}
\def\w{\boldsymbol{w}}

\def\la{\left(}
\def\ra{\right)}
\newcommand{\lb}{\left\langle}
\newcommand{\rb}{\right\rangle}

\newcommand{\ie}{\setlength{\parskip}{0cm} \setlength{\itemsep}{0cm}}
\DeclareMathAlphabet{\mathcal}{OMS}{cmsy}{m}{n}

\allowdisplaybreaks

\title[Pathwise solvability and bubbling in 2D SLLG]{Pathwise Solvability and Bubbling in 2D Stochastic Landau-Lifshitz-Gilbert Equations}

\begin{document}
\maketitle

\begin{center}
	BEN GOLDYS \\ [6pt]
	School of Mathematics and Statistics, \\The University of Sydney, Sydney 2006, Australia. \\ Email: \texttt{beniamin.goldys@sydney.edu.au} \\ [15pt]
	CHUNXI JIAO \\ [6pt]
	School of Mathematics and Statistics, \\ UNSW, Sydney 2052, Australia. \\ Email: \texttt{chunxi.jiao@unsw.edu.au} \\ [15pt]
	CHRISTOF MELCHER \\ [6pt]
	Lehrstuhl für Angewandte Analysis, \\ RWTH Aachen University, Aachen 52062, Germany. \\ Email: \texttt{melcher@math1.rwth-aachen.de}
\end{center}

\let\thefootnote\relax\footnotetext{\textit{Corresponding author:} Chunxi Jiao}

\let\thefootnote\relax\footnotetext{\textit{Key words}. Stochastic Landau-Lifshitz equation, Multiplicative noise, Struwe solutions, Bubbling, Doss-Sussman type transformation.}

\let\thefootnote\relax\footnotetext{\textit{2020 Mathematics Subject Classification}. 35A02, 35A21, 35K59, 35Q60, 60H15.}

\let\thefootnote\relax\footnotetext{\textbf{Statements and Declarations.} The authors have no competing interests to declare that are relevant to the content of this article. This  work was supported by the Australian Research Council Project DP200101866 and the Deutsche Forschungsgemeinschaft (DFG, German Research Foundation) Project 442047500 - SFB 1481. CM is grateful for the support and hospitality of the Sydney Mathematical Research Institute (SMRI).} 

\let\thefootnote\relax\footnotetext{\textbf{Data availability.} No data were used in this article.}

\maketitle

\begin{abstract}
    We investigate the stochastic Landau-Lifshitz-Gilbert (LLG) equation on a periodic 2D domain, driven by infinite-dimensional Gaussian noise in a Sobolev class. We establish strong local well-posedness in the energy space and characterize blow-up at random times in terms of energy concentration at small scales (bubbling). By iteration, we construct pathwise global weak solutions, with energy evolving as a c{\`a}dl{\`a}g process, and prove uniqueness within this class. These results offer a stochastic counterpart to the deterministic concept of Struwe solutions. The approach relies on a transformation that leads to a magnetic Landau-Lifshitz-Gilbert equation with random gauge coefficients.
\end{abstract}

\tableofcontents
\textbf{\today}

\section{Introduction}

	The formation and structure of singularities are among the key questions in the theory of (dissipative) harmonic flows. 
	An interesting prototype is the Landau-Lifshitz-Gilbert equation (LLG)
	\begin{equation*}
		\partial_t \m = -\m \times \la \Delta \m  + \alpha \m \times \Delta \m \ra
	\end{equation*}
	which has its roots in the continuum theory of ferromagnetism. 
	Mathematically it combines the heat and Schrödinger flow of harmonic maps $\m$ from a domain $\mathbb{D}$ into the unit sphere $\mathbb{S}^2 \subset \R^3$. 
	The governing energy is the Dirichlet energy $E(\m)= \frac{1}{2} \int_{\mathbb{D}} |\nabla \m|^2 \, dx$. 
	The damping constant $\alpha>0$ highlights the dissipative term which is nothing but the tension field $\Delta \m + |\nabla \m|^2 \m$, known from the theory of harmonic maps. 
	As a hallmark of precessional spin dynamics, the Hamiltonian terms of LLG gives rise to a large number of interesting oscillatory phenomena. 
	These terms, however, make the equation quasilinear and mathematically more challenging than the semi-linear harmonic map heat flow equation (HMHF).
	
	Stochastic Landau-Lifshitz-Gilbert equation (SLLG)
	\begin{equation*}
		\partial_t \m = - \m \times \left( \Delta \m + \xi + \alpha \m \times \Delta \m \right). 
	\end{equation*}
	with a random field $\xi$ and a suitable interpretation in the framework of stochastic calculus are examined in the context of random fluctuations and thermal activated processes in ferromagnetism. 
	A particular focus is on the switching between stable equilibrium states over energy barriers and related stochastic optimal control problems \cite{Dunst}.
	More recently, singularity mediated processes such as the nucleation, annihilation and collapse of topological field configurations such as magnetic vortices and skyrmions came into focus. 
	The singular nature of stochastic forcing terms in the spatially multi-dimensional case and their effect on the blow-up behaviour of dissipative harmonic flows is a widely unmapped area of research.

	The theory of partial regularity in the energy-critical two-dimensional case is  largely based on the seminal work of Struwe on the harmonic map heat flows from surfaces \cite{Struwe_1985}.
    There is a canonical class of energy-decreasing weak solutions for which the initial value problem is uniquely and globally solvable \cite{Freire,Struwe_1985}, nowadays known as \textit{Struwe solutions}. 
	In this class, the blow-up scenario is described as the concentration and drop of a quantized amount of energy accompanied by changes of the local topology of the solution. This so-called bubbling 
    occurs at a finite number of space-time points depending on the initial energy \cite{ChangDingYe}. In particular, there is a precise energy threshold below which solutions remain globally regular. 
    In a refined approach, the so-called bubble-tree construction restores the energy identity at singular times \cite{Qing}. Finally, the so-called reverse bubbling construction \cite{Bertsch,Topping} 
    serves to construct a weak solution to the initial value problem that are not in the class of Struwe solutions. The theory has been extended to LLG \cite{ChenGuo, Harpes}. Recently, a very general construction 
    of blow-up solutions has been developed in \cite{Wei}. Regularity and blow-up criteria are crucial to understand the stability and dynamics of magnetic vortices and skyrmions \cite{DoeringMelcher, Gustafson, KurzkeMelcherMoser, KurzkeMelcherMoserSpirn,  Melcher_Sakellaris}. However, there is limited theory about their behaviour in the presence of noise \cite{BMS_stochastic, Chugreeva_Melcher}, and we refer the reader to \cite[Section 2.3.3]{Banas_book} for numerical experiments on pathwise (singular) dynamics under different noise intensities.
	
	Passing to a stochastic framework for bubbling and Struwe type solutions requires first of all a strong well-posedness concept of pathwise solutions in energy space for small (stopping) times. This challenges existing mathematical strategies for the SLLG based on a Faedo-Galerkin scheme and stochastic compactness arguments leading to weak martingale solutions as developed in \cite{BGJ1}. Improvements regarding pathwise uniqueness based on Doss-Sussmann type transformations \cite{BGJ2} and the concept of rough paths \cite{GussettiHocquet} are typically restricted to the (subcritical) one-dimensional case.    
	In comparison, the stochastic HMHF is a semilinear parabolic system that allows a mild formulation using the underlying semigroup and the Duhamel principle. 
	Existence and uniqueness of local pathwise solutions in higher order Sobolev spaces can then be derived from stochastic convolution estimates and a fixed point argument. This is the starting point for a program to construct Struwe-type solutions for the stochastic HMHF developed by Hocquet in \cite{Hocquet1}. A following up work \cite{Hocquet2} shows that the comparison argument used to construct blow-up solutions for the axially symmetric HMHF extends to the case of random perturbations, proving that finite time blow-up happens with positive probability. 
	
	In this work we develop a framework to extend the concept of Struwe-type solutions to SLLG and allows for a precise bubbling analysis. 
	More precisely, we examine processes $\m$ of maps from the two-dimensional flat torus $\T^2$ to the unit sphere $\St$ that satisfy 
	\begin{equation}\label{eq: intro sLLG}
		d\m = -\m \times \left( \Delta \m + \alpha \m \times \Delta \m \right) + \m \times \circ \ dW.
	\end{equation}
	Here, $ W $ is a Wiener process on a probability space $(\Omega, \mc F, \P)$ taking values in $ H^\sigma(\T^2; \R^3) $ for some $ \sigma > 0 $, which is formally defined in Section \ref{Section: Wiener process}. 
	The equation is understood in the framework of Stratonovich calculus
	that guarantees that the geometric constraint 
	$\m(t,x,\omega) \in \mathbb{S}^2$ is preserved along the flow. 
	
	Our main task is the construction of strong (in the PDE sense) solutions of \eqref{eq: sLLG} for given initial conditions $\m_0 \in H^1(\T^2; \St)$. 
	Our approach is based on a Doss-Sussmann type transformation whose conceptual application in the context of SLLG was developed by one of us in a series of work \cite{GGL,GLT}. 
    The key idea is that the reduced problem $\partial_t \m = \xi \times \m$ gives rise to a pure precession of $\m$
	about the field $\xi$ that can be captured by a parametrized rotation $ \Y $ which is a unique solution of a stochastic differential equation taking values in $ C^{\frac{1}{2}-}([0,T]; H^\sigma(\T^2; SO(3))) $, $ \P $-a.s. 
	The random transformation $\u= \Y^{-1} \m$ gives rise to a new process $\u$ and 
	a random PDE in form of a gauged LLG equation
	\begin{equation}\label{eq: intro u}
		\partial_t \u + \u \times \left( \Delta_A \u + \alpha \u \times \Delta_A \u \right) =0,
	\end{equation}
	with the covariant gradient and Laplacian 
	$ \nabla_A \u = \nabla \u + A \u $ and $ \Delta_A = \nabla_A \cdot( \nabla_A \u) $, respectively, 
	where $ A= \Y^{-1} \nabla \Y $. 
	This is the LLG equation for the magnetic Dirichlet energy  
	\begin{equation*}
		E(\u, A):= \frac{1}{2} \int_{\mathbb{T}^2} |\nabla_A \u|^2 \, dx.
	\end{equation*}
	For $\sigma \geq 1$ and $\nabla \u \in L^\infty(0,T; L^2(\T^2))$, the energies
	$E(\u, A)$ and $E(\u)$ differ by an $ L^2 $-estimate of $ A $ which is $ \P $-a.s. finite, so that the conventional weak solution 
	concept in the space $H^1((0,T) \times \T^2; \St) \cap L^\infty(0,T; H^1(\T^2; \St))$ applies in a pathwise 
	fashion. 
	It is not difficult to obtain an equivalence result between pathwise weak solutions of $ \eqref{eq: intro u} $ and the SLLG (see Lemma \ref{Lemma: u vs m weak}), which allows us to focus on constructing solutions of \eqref{eq: intro u} instead of \eqref{eq: sLLG}.  
	
	The idea of using a gauge transformation to convert an SPDE into a random PDE has also been applied in the context of stochastic nonlinear Schr\"odinger equations \cite{BarbuRoecknerZhang, deBouardFukuzumi} where $\Y$ is unitary and $A$ takes the form $A=ia$ for a real valued process $a$.
	
	For local strong solutions $\u$ of \eqref{eq: intro u} in the critical Sobolev space, which is $H^1(\mathbb{T}^2; \R^3)$ in our case, it is customary 
	to first construct solutions in higher order Sobolev spaces where the equation behaves subcritical. 
	We observe that \eqref{eq: intro u} can be considered as a perturbation of the ungauged LLG equation, which is parabolic with an additional tangent field
	$F(\u)$ that is bounded by
	\begin{equation*}
		|F(\u)| \lesssim \la |A|^2+|\nabla \cdot A| + |A||\nabla \u| \ra,
	\end{equation*}
	motivating the regularity assumption $\sigma \geq 2$ that also guarantees a pointwise bound of $A$.
	The method of choice for quasilinear parabolic systems is based on a linearization procedure and fixed point argument that uses maximal regularity properties.
	Following Struwe's program, extension to the critical case by approximation requires an $\eps$-regularity result arising from the special analytic structure of the geometric nonlinearities in combination with a localized version of Ladyzhenskaya's interpolation inequality, and a covering argument. 
	The combination of small energy compactness, the removability of singularities and energy gap featured by harmonic maps between $\St$ allows one to describe the bubbling scenario. 
	Finally, the merging of local solutions leads to the Struwe solution. 
	A key issue in the case of a random PDE is to keep track of the dependence on the infinite-dimensional parameter arising from stochasticity. 
	
	We remark that in \cite{Hocquet1}, the stochastic HMHF features an $H^1$ noise which seems to be critical for the basic energy estimated. 
	This, however, seems to require a special symmetry assumption leading to a cancellation in order to avoid $L^\infty$ properties. 
	Our approach for SLLG requires higher regularity assumptions on the noise but no structural conditions. 
	Moreover, we provide the precise energy threshold for bubbling (see Theorem \ref{Theorem: Struwe u}). 
	
	Next, we provide a precise formulation of \eqref{eq: intro sLLG} and state our main result.

\subsection{Wiener process}\label{Section: Wiener process}
	Let $(\Omega,\mc F,\la\mc F_t\ra_{t \geq 0},\P)$ be a probability space with the filtration \smash{$\b F=(\mc F_t)_{t \geq 0}$}.  
	We assume that this probability space supports a sequence $(W_j)$ of independent, real-valued Brownian motions. 
	Note that $W=(W_j)$ defines a cylindrical Wiener process on the Hilbert space $\ell^2$ of square-summable sequences. 
	We also assume that for each $ t \geq 0 $, $\mc F_t$ is a completion of the $\sigma$-algebra $\sigma\la W_j(s);\,j \ge 1,\,s\le t\ra$. This implies 
	\begin{equation*}\label{eq_ft}
		\mc F_{t-}=\mc F_t=\mc F_{t+}\,,
	\end{equation*}
	a fact that we will need later. 
	
	For simplicity, let $ \b H^\sigma := H^\sigma(\T^2; \R^3) $ and $ {\b L}^p := L^p(\T^2; \R^3) $ for $ \sigma \geq 0 $ and $ p \in [1,\infty] $. 
	Let $\la g_j\ra_{j \in \b N}\subset\b L^2$ be a sequence that satisfies $ \sum_{j=1}^\infty|g_j|_{\b L^2}^2<\infty $. 
	Let $P:\ell^2\to\b L^2$ be defined as 
	\begin{align*}
		P l=\sum_{j=1}^\infty\lb l,h_j\rb g_j,\quad l\in \ell^2\,,
	\end{align*}
	where $(h_j)$ is a complete orthonormal basis of $\ell^2$. 
	Then the process 
	\[(P W)(t)=\sum_{j=1}^\infty W_j(t)g_j, \quad t \geq 0,\]
	defines an $\b L^2$-valued Wiener process with the covariance operator 
	$ Q= P^\star P= \sum_{j=1}^\infty g_j\otimes g_j $, 
	and we have 
	$ \E[ |P W(t)|^2_{\b L^2} ] = t \mathrm{Tr\, Q} =t\sum_{j=1}^\infty|g_j|_{\b L^2}^2<\infty $ for any $ t \geq 0 $.
	
	The main result (Theorem \ref{Theorem: Struwe m}) will require a stronger condition that for a certain $ \sigma \geq 2 $, 
	\begin{equation}\label{eq_q3}
		q^2(\sigma):= \sum_{j=1}^\infty |g_j|^2_{\b H^{\sigma}}<\infty \,.
	\end{equation}
	Assume that \eqref{eq_q3} hold for some $ \sigma \geq 2 $. 
	For every $ j \geq 1 $, we define a bounded linear operator $G_j:\b L^2\to\b L^2$ by
	\begin{align*}
		G_j\varphi=\varphi\times g_j,
		\quad 
		\varphi \in \b L^2. 
	\end{align*}
	Then, \smash{$\|G_j\|_{\mc L(\b L^2, \b L^2)}=|g_j|_{\b L^\infty} \lesssim |g_j|_{\b H^2}$} and $ G_j $ has adjoint \smash{$G_j^\star=-G_j$} for every $ j \geq 1 $, and 
	$ \sum_{j=1}^\infty\left\|G_j\right\|^2_{\mc L(\b H^\sigma,\b H^\sigma) }<\infty $, where $ \mc L(X,Y) $ denotes the space of bounded linear operators mapping from $ X $ to $ Y $. 
	Similarly, we define $ G: \b L^2 \to {\mc L}(\ell^2,\b L^2) $ (and $ G(\varphi): \ell^2 \to \b L^2 $) by
	\begin{align*}
		G(\varphi)l=\varphi \times Pl=\sum_{j=1}^\infty \lb l,h_j \rb G_j\varphi, \quad \varphi \in \b L^2, \ l \in \ell^2\,.
	\end{align*}
	In particular, 
	\begin{align*}
		G(\varphi)W(t)=\sum_{j=1}^\infty G_j \varphi \ W_j(t)\,.
	\end{align*}
	Thus, if $\v$ is an $ \b F $-progressively measurable process such that $\v\in L^\infty(0,T,\b L^2)$, $ \P $-a.s. then we can define the Stratonovich integral 
	\begin{align*}
		\int_0^t G(\v(s))\circ dW(s)=\int_0^t S(\v(s))\,ds+\int_0^t G(\v(s))\,dW(s)\,,
	\end{align*}
	where $S$ is the Stratonovich correction term $ S(\v)=\smash{\frac{1}{2}\sum_{j=1}^\infty G_j^2\v} $.
	We refer the reader to \cite{dz1} for details of the It\^o integral $ \int_0^t G(\v(s))\,dW(s) $ on Hilbert space.

\subsection{Main result}
	Formally, the SLLG equation reads 
	\begin{equation}\label{eq: sLLG}
		\begin{aligned}
			d\m(t) 
			&= -\m(t) \times \left[ \Delta \m(t) + \alpha \la \m(t) \times \Delta \m(t) \ra \right] dt \\
			&\quad + G(\m(t))\circ dW(t), \quad t \geq 0,  
		\end{aligned}	
	\end{equation}
	with initial condition $ \m(0)= \m_0 \in H^1(\T^2; {\b S}^2) $.
	
\begin{definition}[Weak martingale solution]\label{Defn: martingale} 
    For every $ T \in (0,\infty) $, a system $ (\Omega, \mc F, (\mc F_t)_{t \in [0,T]}, \P, W, \m) $ consisting of 
    a filtered probability space supporting a Wiener process $ W $ as in Section \ref{Section: Wiener process} 
    and a progressively measurable process $ \m: [0,T] \times \Omega \to \b L^2 $, is said to be a weak martingale solution of \eqref{eq: sLLG} on $ [0,T] $ if $ \P $-a.s.
		\begin{enumerate}[leftmargin=*,parsep=1pt]
			\item[(a)]
			$ |\m(t,x)|=1 $ for a.e.-$ (t,x) \in [0,T] \times \T^2 $,
			
			\item[(b)]
			$ \m \in C([0,T]; \b L^2) \cap L^\infty(0,T; \b H^1) $ and $ \m \times \Delta \m \in L^2(0,T; \b L^2) $, 
			
			\item[(c)] 
			for every $ t \in [0,T] $ and $ \psi \in C_0^\infty(\T^2; \R^3) $, the equality
			\begin{align}
				\lb \m(t)-\m_0, \psi \rb_{\b L^2} 
				&= \int_0^t \lb \m \times \nabla \m, \nabla \psi - \alpha \nabla(\m \times \psi)\rb_{\b L^2} ds \label{eq: sLLG weak} \\
				&\quad + \int_0^t \lb \psi, G(\m) \circ dW(s) \rb_{\b L^2}. \nonumber
			\end{align}
		\end{enumerate}
	\end{definition}
	
	The nonlinear term $ \m \times \Delta \m $ in Definition \ref{Defn: martingale}(b) is understood in the sense of distribution
	\begin{align*}
		\m \times \Delta \m := \nabla \cdot ( \m \times \nabla \m), 
	\end{align*}
	consistent with the weak formulation \eqref{eq: sLLG weak},  
	so $\b L^2$-regularity of this term can be seen as a form of finite dissipation property and will play a crucial role in the question of uniqueness of suitable weak solutions.
	
	\begin{definition}[Local strong solution]\label{Defn: local strong} 
		Given $ T \in (0,\infty) $ and a filtered probability space supporting a Wiener process $ (\Omega, {\mc F}, (\mc F_t)_{t \in [0,T]}, \P, W) $ as in Section \ref{Section: Wiener process}. A pair $ (\m,[\tau_0,\tau)) $ defined on $ (\Omega, {\mc F}, (\mc F_t)_{t \in [0,T]}, \P) $ is said to be a local strong solution of \eqref{eq: sLLG} if $ \m $ is a progressively measurable process and $ \tau_0, \tau\in (0, T]$ are stopping times such that $ \P $-a.s. for every $ t \in (\tau_0,\tau) $, 
		\begin{enumerate}[leftmargin=*,parsep=1pt]
			\item[(a)] 
			$ \m \in C([\tau_0, t];H^1(\T^2; \St)) $, 
			
			\item[(b)]
			$ \Delta \m \in L^2(\tau_0,t; \b L^2) $,
			
			\item[(c)]
			the equality
			\begin{align*}
				\m(t) 
				&= \m(\tau_0) - \int_{\tau_0}^t \left[ \m \times \la \Delta \m + \alpha \m \times \Delta \m \ra \right](s) \ ds \\
				&\quad + \int_{\tau_0}^t G(\m(s)) \circ dW(s), 
			\end{align*}
			holds in $ \b L^2 $.
		\end{enumerate}
	\end{definition}

	Now we state the main result of this paper, where the energy concentration is described by a lower bound for energies over balls $ B_r(x_n) \subset \R^2 $ of arbitrary radii $ r $ centred at random locations $ x_n \in \T^2 $. 
	The result follows directly from inverse transformation of the Struwe solution $ \u $ of \eqref{eq: intro u}, for which we postpone the details to Theorem \ref{Theorem: Struwe u} and \ref{Theorem: u unique} in Section \ref{Section: well-posedness}. 
	
	\begin{theorem}\label{Theorem: Struwe m}
		Assume that \eqref{eq_q3} holds for $ \sigma = 4 $. 
		For every $\m_0 \in H^1(\T^2; \St)$, there exists a weak martingale solution 
		$ (\Omega, \mc F, (\mc F_t)_{t \in [0,T]}, \P, W, \m) $ of \eqref{eq: sLLG} with $ \m(0) = \m_0$ in the sense of Definition \ref{Defn: martingale}
		and an increasing sequence of $ (\mc F_t)_{t \in [0,T]} $-stopping times 
		$(\tau_k)_{k \in \N}$, 
		such that 
		$ T \in (\tau_M, \tau_{M+1}) $ for some $ M \geq 0 $ with $ \tau_0 = 0 $, and that for all $ k \leq M $,
		\begin{enumerate}[leftmargin=*,parsep=1pt]
			\item[(a)]  
			$ \m|_{[\tau_k, \tau_{k+1} \wedge T)} $ is a (local) strong solution in the sense of Definition \ref{Defn: local strong}, 
			
			\item[(b)] 
			$ \lim_{t_i \nearrow \tau_k} \m(t) = \m(\tau_k) $ weakly in $ \b H^1 $, $ \P $-a.s. for any sequence of stopping times $ (t_i)_{i \in \b N} $ with $ t_i < \tau_k $ for all $ i \in \b N $ and $ t_i \nearrow \tau_k $, $ \P $-a.s.
			
			\item[(c)] 
			there exist a random variable $ N_k \in [0,\infty) $, $ \P $-a.s. and random points $ x_1, \ldots, x_{N_k} \in \T^2 $ where 
			for every $r>0$, 
			\begin{align*}
			\limsup_{t \nearrow \tau_k} \frac{1}{2} \int_{B_r(x_n)} |\nabla \m(t)|^2 \, dx \geq 4 \pi 
			\quad n=1, \dots, N_k, \quad \P\text{-a.s.} 
			\end{align*}
		\end{enumerate} 
	\end{theorem}

	\begin{theorem}\label{Theorem: unique m}
		For $ i=1,2 $, let $ (\Omega, \mc F, (\mc F_t)_{t \in [0,T]}, \P, W, \m_i) $ be a weak martingale solution to \eqref{eq: sLLG} with c\`adl\`ag energy process $ \{|\nabla \m_i(t)|^2_{\b L^2} : t \in [0,T]\} $. Then, for every $ t \in [0,T] $, 
		\begin{align*}
			|\m_1(t) - \m_2(t)|_{\b L^2} = 0, \quad \P\text{-a.s.} 
		\end{align*}
	\end{theorem}
	
	Applying an infinite-dimensional Yamada-Watanabe theorem (for instance, \cite[Theorem 2]{Ondrejat_YW}), we obtain a solution of \eqref{eq: sLLG} strong in the probabilistic sense.

	The rest of the paper is devoted to the proof of Theorem \ref{Theorem: Struwe m} and \ref{Theorem: unique m}. 
	In Section \ref{Section: Transformation}, we introduce the transformation operator $ \Y $ (and $ A $) and an equivalence result between the transformed equation and \eqref{eq: sLLG} (Lemma \ref{Lemma: u vs m weak}). 
	In Section \ref{Section: solution of transformed equation}, we show the existence of a maximal solution $ \u $ of the gauged LLG equation \eqref{eq: intro u} in $ \b H^2 $ and verify measurability of the local solutions and the maximal time. 
	To prepare for the approximations of the critical case, we collect energy and higher order estimates of $ \u $ in Sections \ref{Section: energy estimates} and \ref{Section: regularity small energy}. 
	We formally describe the Struwe solution of \eqref{eq: intro u}, or equivalently \eqref{eq: u}, in Section \ref{Section: well-posedness} (Theorem \ref{Theorem: Struwe u} and \ref{Theorem: u unique}), and leave the proofs to Sections \ref{Section: Proof existence} and \ref{Section: Proof uniqueness}. 
	We conclude with some remarks on the regularity of noise in Section \ref{Section: Concluding remarks}.

\section{Transformation}\label{Section: Transformation}

\subsection{Transformation operator $ \Y $ and $ A $}\label{Section: properties of Y}
We start with recalling some facts about the equation 
\begin{equation}\label{eq: dy}
	dY(t)= S(Y(t)) \,dt+ G(Y(t))\,dW(t),\quad t \geq 0\,,
\end{equation} 
with initial condition $ Y(0)=y_0\in\b L^2 $. 

\begin{lemma}\label{Lemma: Y solution}
	Assume that $y_0\in\b H^\sigma$ and \eqref{eq_q3} holds for some $\sigma \geq 2$. 
	There exists a unique solution $ Y(\cdot;y_0) $ of \eqref{eq: dy} taking values in $ \b H^\sigma $, $ \P $-a.s. 
	Moreover, for every $ T \in (0,\infty) $ and $ p \in [1,\infty) $,  
	\begin{align*}
		\E\left[\sup_{t\in [0,T]}\left|Y(t;y_0) \right|_{\b H^\sigma}^p \right]<\infty\,,
	\end{align*}
	and $Y(\cdot;y_0) \in C^\gamma ([0,T];\b H^\sigma) $, $ \P $-a.s. for $ \gamma \in (0,\frac{1}{2}) $. 
\end{lemma}
\begin{proof}
	By the assumption \eqref{eq_q3}, the equation \eqref{eq: dy} has Lipschitz coefficients in $\b H^\sigma$ for $\sigma \geq 2$. The result is then standard following from \cite[Theorem 7.2]{dz1} and Kolmogorov's continuity criterion. 
\end{proof}

With an abuse of notation, we still write $ \b H^\sigma $ and $ \b L^p $ for matrix or tensor-valued functions when the codomain is clear.
Similar to $ G_j $ in Section \ref{Section: Wiener process}, for $ \varphi \in L^2(\T^2; \R^3 \otimes \R^3) $ and $ h \in \R^3 $, define
\[
(\c G_j\varphi) h=\la\varphi h\ra\times g_j = G_j (\varphi h),\quad j \geq 1\,.
\] 
The operator \smash{$\c G_j: L^2(\T^2; \R^3 \otimes \R^3) \to L^2(\T^2; \R^3 \otimes \R^3) $} is bounded with 
\begin{equation*}
	\|\c G_j \|_{\b L^2 \to \b L^2}=|g_j|_{\b L^\infty}, 
\end{equation*}
and adjoint \smash{$\c G_j^\star=-\c G_j$}. 
Moreover, 
$ \sum_{k=1}^\infty\|\c G_k\|^2_{\b H^\sigma \to\b H^\sigma }<\infty $ under the assumption \eqref{eq_q3}. 
Then we can define operators $ \c G $ and $ \c S $ 
in a similar manner to $ G $ and $ S $, respectively, with $ \c G_k $ in place of $ G_k $. 
Consider the equation
\begin{equation}\label{eq: mc Y}
	d\, \Y(t)= \c S(\Y(t)) \,dt+ \c G(\Y(t)) \,dW(t),\quad t \geq 0 \,,
\end{equation} 
with initial condition $ \Y(0)=I_{\R^3} $. 

\begin{lemma}\label{Lemma: Y flow}
	Assume that \eqref{eq_q3} holds for some $ \sigma \geq 2 $. 
	There exists a unique solution $\Y$ of \eqref{eq: mc Y} taking values in $ H^\sigma(\T^2; \R^3 \otimes \R^3) $, $ \P $-a.s. with the following properties. 
	\begin{enumerate}[leftmargin=*,parsep=0pt]
		\item[(a)]
		For every $ T\in(0,\infty) $ and $ p \in [1,\infty) $,  
		\begin{align*}
			\E \left[ \sup_{t \in [0,T]} |\Y(t)|_{H^\sigma(\T^2; \R^3 \otimes \R^3)}^p \right] < \infty \,, 
		\end{align*}
		and $\Y\in C^\gamma\la [0,T]; H^\sigma(\T^2; \R^3 \otimes \R^3) \ra$, $ \P $-a.s. for $\gamma\in\la 0,\frac{1}{2}\ra$. 
		
		
		\item[(b)]
		For every $ t \geq 0 $, $ x \in \b T^2 $, $ h \in \R^3 $ and a.e.-$ \omega \in \Omega $, 
		\[
		|\Y(t,x,\omega)h|_{\R^3}=|h|_{\R^3} \,.
		\]
		
		\item[(c)]
		For every $ t \geq 0 $ and a.e.-$ \omega \in \Omega $, the mapping $\Y(t,\omega) \in H^\sigma(\T^2; \R^3 \otimes \R^3)$ can be identified with a bounded linear operator $\Y(t,\omega):\b L^2\to\b L^2$ by 
		\[
		(\Y(t,\omega)\varphi)(x)= \Y(t,x,\omega)\varphi(x),\quad \varphi\in\b L^2\,,
		\]
		where $ \Y(t,\omega)$ defines an isometric isomorphism on $\b L^2$ and an isomorphism on $\b H^\sigma$ for $ \sigma \geq 2 $, with $ \Y^{-1}(t,\omega)= \Y^\star(t,\omega)$ governed by the equation
		\begin{equation}\label{eq: Y^* eq}
			d \Y^\star(t) 
			= -\Y^\star(t) \c G \circ dW(t), \quad t \geq 0, \ \Y^\star(0) = I_{\R^3}.
		\end{equation}
		
		\item[(d)]
		For every $ t \geq 0 $ and a.e.-$ \omega \in \Omega $,
		\begin{align*}
			\Y(t,\omega) (f\v) &= f \Y(t,\omega) \v, \\
			\Y(t,\omega) (\v \times \varphi) &= \Y(t,\omega) \v \times \Y(t,\omega) \varphi, 
		\end{align*}
		for any $ f: \T^2 \to \R $ and $ \v: \T^2 \to \R^3 $ such that $ f \v \in \b L^2 $, and $ \varphi \in \b L^2 $. The same equalities hold for $ \Y^\star(t,\omega) $. 
	\end{enumerate} 
\end{lemma}
\begin{proof}
	Part (a) follows from Lemma \ref{Lemma: Y solution}. 
	Part (b) and (c) are analogous to the proof of \cite[Lemma 3.3]{GGL}, noting that the matrix $ \Y^\star(t,x,\omega) $ is the transpose of $ \Y(t,x,\omega) $ for every $ (t,x) \in [0,T] \times \T^2 $ and a.e.-$ \omega \in \Omega $. 
	Part (d) is similar to the proof of \cite[Lemma 3.5]{GGL}.
\end{proof}

We define the process $ A $ by
\begin{equation}\label{eq: defn A}
	A:= \Y^\star (\nabla \Y): [0,T] \times \T^2 \times \Omega \to \R^2 \otimes \R^{3 \times 3},
\end{equation}
where its components 
$ A_i = \Y^\star (\nabla_{x_i} \Y) $, $ i=1,2 $,  
take values in $ \mathfrak{so}(3) $ the space of skew-symmetric $ 3 \times 3 $ matrices, for every $ (t,x) \in [0,T] \times \T^2 $, $ \P $-a.s. 
By Lemma \ref{Lemma: Y flow}(a), we have 
\begin{align*}
	A &\in L^p(\Omega; L^\infty(0,T;H^{\sigma-1}(\T^2; \R^2 \otimes \mathfrak{so}(3)))), \\
	A &\in C^\gamma([0,T]; H^{\sigma-1}(\T^2; \R^2 \otimes \mathfrak{so}(3))), \quad \P \text{-a.s.}
\end{align*}
for every $ T \in (0,\infty) $, $ p \in [1,\infty) $ and $ \gamma \in (0,\frac{1}{2}) $. 
Moreover, applying It{\^o}'s lemma, 
\begin{align*}
	A(t) 
	&= \frac{1}{2} \int_0^t \sum_{k=1}^\infty \Y^\star(s) \la (\nabla \c G_k) \c G_k - \c G_k (\nabla \c G_k) \ra \Y(s) \ ds \\
	&\quad + \int_0^t \sum_{k=1}^\infty \Y^\star(s) (\nabla \c G_k) \Y(s) \ dW_k(s), \quad t \geq 0. 
\end{align*}
The relevance of the process $ A $ manifests itself in the behaviour of the spatial gradient under transformation $ Y $. 
Here, we have the covariant derivative
\begin{equation}\label{eq: defn gradA}
	\nabla_A := \Y^\star \nabla \Y, 
\end{equation}
and the covariant Laplacian
\begin{equation}\label{eq: defn D2A}
	\Delta_A 
	:= \nabla_A \cdot \nabla_A 
	= \Y^\star \Delta \Y.
\end{equation}
Taking into account $ (Y^\star Y)(t,x)= I_{\R^3} $, $ \P $-a.s. we have 
\begin{align*}
	\nabla_{A} = \nabla + A, \quad
	\Delta_{A} = \Delta + A \cdot \nabla + \nabla \cdot A + A^2,
	\quad \P\text{-a.s.}
\end{align*}  
where the dot product means contracting over the $ \R^2 $ components, and 
$ A^2= A \cdot A = A_1^2 + A_2^2: [0,T] \times \T^2 \times \Omega \to \R^{3\times3} $ mapping to a symmetric matrix $ \P $-a.s.   
Hence, transformation by means of $\Y$ is only reflected in the equation by means of the process $A$. 

\begin{remark}\label{Remark: A curvature} 
	From a coordinate transformation on $\St$, $A = A(\omega)$ is a pure gauge for a.e.-$ \omega \in \Omega $, so that the curvature  
	\[
	\partial_1 A_2 - \partial_2 A_1 + [A_1, A_2] = 0.
	\]
	With $ \nabla^\perp = (\partial_2, -\partial_1) $, this provides a pointwise control $|\nabla^\perp A| = |[A_1, A_2]|$ so that Gagliardo-Nirenberg type interpolation estimates can be reduced to divergence bounds. 
	In fact, for any component $ b=\smash{(A^{(j,k)}_1, A^{(j,k)}_2)}: \T^2 \to \R^2 $, $1 \le j,k \le 3$, 
	it holds that 
	\[
	|\nabla b|^2 =|\nabla \cdot b|^2+|\nabla^\perp \cdot b|^2 +
	\nabla^\perp \cdot \left( b_1 \nabla b_2 - b_2 \nabla b_1 \right),
	\]
	and thus 
	\begin{equation} \label{eq:div_curl}
		|\nabla A|^2 =|\nabla \cdot A|^2+|[A_1, A_2]|^2 -
		\nabla^\perp \cdot \mathrm{Tr} \left( A_1  \nabla A_2 - 
		A_2  \nabla A_1 \right).
	\end{equation}
	A gauge transformation is a smooth map $X:\R \times \T^2 \to \mathrm{SO}(3)$ acting on $ \v $ and $ A $. 
	Then for $A'= X^\star A X + X^\star \nabla X$, 
	\begin{align*}
		X^\star \nabla_A X = \nabla + A', \quad  
		X^\star \partial_t X = \partial_t + X^\star (\partial_t X).
	\end{align*}
	The rough dependence of $A$ on $t$ limits the spectrum of possible gauge transformations that bring $A$ into a particularly convenient form.  
\end{remark}

\begin{example}\label{Example: 1D}
	\rm 
	Let $ W $ be a real-valued Wiener process. 
	The transformation $\Y=e^{WG}$ with $G\u= \u \times g$ for some $g \in H^1(\T^2; \St)$ 
	reduces to
	\[
	\Y(t) = I_{\R^3} + (\sin W(t)) G + (1- (\cos W(t))) G^2.
	\]
	A straightforward calculation shows that the corresponding gauge potential reads
	\[
	A= \sin W(t) \nabla G + (1- \cos W(t)) (\nabla G \, G - G \, \nabla G), 
	\]
	and $ \nabla \cdot A = (\sin W(t)) \Delta G + (1- \cos W(t)) ((\Delta G) G - G \Delta G) $.
\end{example}

\subsection{Transformed equation}
Fix $ T \in (0,\infty) $ and $ \gamma \in (0,\frac{1}{2}) $. 
Let $ \Y $ be the solution of \eqref{eq: dy} described in Lemma \ref{Lemma: Y flow}, and let $ A $ be given by \eqref{eq: defn A}.
Consider the following equation with random coefficients
\begin{equation}\label{eq: u}
	\begin{aligned}
		\partial_t \u(t)
		&=\alpha \la \Delta_A \u+|\nabla_A \u|^2 \u \ra(t) - \la \u \times\Delta_A \u \ra (t), \quad 
		t \in [0,T].
	\end{aligned}
\end{equation}
The right-hand side can be reduced to the original Landau-Lifshitz terms (precession and damping) with a lower-order perturbation 
\begin{equation}\label{eq: du = LL + F}
	\partial_t \u
	=\alpha \la \Delta \u+|\nabla \u|^2\u \ra - \u \times\Delta \u + F(t,A,\u),
\end{equation}
where 
\begin{align}
	F(t,a,\v)
	&= \alpha \la a(t) \cdot \nabla \v + (\nabla \cdot a(t))\v + a^2(t)\v \ra \label{eq: F} \\
	&\quad + \alpha \la |a(t) \v |^2v +2 \lb \nabla \v , a(t)\v \rb \v \ra \nonumber \\
	&\quad - \v \times \la a(t) \cdot \nabla \v + (\nabla \cdot a(t)) \v + a^2(t) \v \ra, \nonumber
\end{align}
for $ t \in [0,T] $, and sufficiently regular functions $ a: [0,T] \times \T^2 \to \R^2 \otimes \mathfrak{so}(3) $ and $ \v : \T^2 \to \R^3 $. 
If $ \v : \T^2 \to {\b S}^2 $, then there exists a universal constant $ c > 0 $ such that  
\begin{equation}\label{eq: |F|}
	|F(t,a,\v)| \leq c\la  |a(t)|^2 +|\nabla \cdot a(t)| + |a(t)||\nabla \v| \ra,
\end{equation}
and similarly,
\begin{equation} \label{eq: |DF|}
	\begin{aligned}
		|\nabla F(t,a,\v)| 
		\leq c &\Big[ |a(t)||\nabla a(t)| +|\nabla (\nabla \cdot a(t))|	\\
		& \ + (|\nabla a(t)| + |a(t)|^2)|\nabla \v | \\
		& \ +|a(t)| \la |\nabla^2 \v| + |\nabla \v |^2 \ra \Big],
	\end{aligned}
\end{equation}
for $ t \in [0,T] $, provided that $ a $ and $ \v $ have well-defined derivatives. 

\begin{remark}\label{Remark: Coulomb gauge}
	Under the Coulomb gauge condition $\nabla \cdot A=0$, the highest order $A$ derivatives in $ F$ and $\nabla F$ are eliminated. 
\end{remark}

Given a suitable progressively measurable solution $ \u $ of the transformed equation \eqref{eq: u}, we can show that there exists a weak martingale solution to the SLLG equation \eqref{eq: sLLG} by transformation. 
This equivalence has been shown in \cite{GGL,GLT}. 
We adapt this result to our problem setting below.

	%
	%

\begin{lemma}\label{Lemma: u vs m weak}
	Assume that \eqref{eq_q3} holds for $ \sigma=2 $. 
	If a process $ \u: [0,T] \times \T^2 \times \Omega \to \b R^3 $ defined on a filtered probability space $ (\Omega, \mc F, (\mc F_t)_{t \in [0,T]}, \P) $ satisfies $ \u(0) = \m_0 \in H^1(\T^2; \St) $ and the following properties:
	\begin{enumerate}[(a), leftmargin=*]
		\item 
		$ \u $ is progressively measurable with respect to $ (\mc F_t)_{t \in [0,T]} $,
		
		\item 
		$ |\u(t,x)|=1 $ for a.e.-$ (t,x) \in [0,T] \times \T^2 $, $ \P $-a.s.
		
		\item 
		$ \u \in C([0,T]; \b L^2) \cap L^\infty(0,T; \b H^1) \cap H^1(0,T; \b L^2) $, $ \P $-a.s. and it satisfies \eqref{eq: u}, $ \P $-a.s.
	\end{enumerate}
	then $ \m = \Y\u $ is a weak martingale solution of \eqref{eq: sLLG} in the sense of Definition \ref{Defn: martingale}. 
\end{lemma}
\begin{proof}
	By Lemma \ref{Lemma: Y flow}(b), 
	$ |\m(t,x)| = |(\Y\u)(t,x)| = 1 $ and 
	$ |\nabla \m(t,x)| = |\nabla_A \u(t,x)|$ for a.e.-$ (t,x)\in [0,T] \times \T^2 $, $ \P $-a.s. 
	Recall that $ \nabla_A = \nabla +A $, $ \P $-a.s. and $ A \in L^p(\Omega; L^\infty(0,T; \b H^1)) $ for $ p \in [1, \infty) $. 
	Then it is immediate that $ \m $ satisfies part (a) and (b) of Definition \ref{Defn: martingale}. 
	For part (c), applying the It{\^o} formula to $ \m = \Y \u $, 
	\begin{align*}
		d\m
		&= 
		\c S(\Y) \u\,dt
		+ \c G(\Y) \u \,dW
		+ \Y d\u.
	\end{align*}	
	This implies that for any $ \psi \in C_0^\infty(\T^2; \R^3) $,
	\begin{align*}
		\lb \m(t)-\m_0, \psi \rb_{\b L^2}
		&= 
		\int_0^t \lb \Y(s) (\partial_s \u(s)), \psi \rb_{\b L^2} \ ds
		+ \int_0^t \lb S(\m(s)),\psi \rb_{\b L^2} ds \\
		&\quad + \int_0^t \lb \psi, G(\m(s)) \, dW(s) \rb_{\b L^2}, 
	\end{align*}
	where by \eqref{eq: defn gradA}, Lemma \ref{Lemma: Y flow}(c) and (d), 
	\begin{align*}
		\int_0^t \lb \Y (\partial_s \u), \psi \rb_{\b L^2} ds 
		&= \int_0^t \lb \partial_s \u, \Y^\star \psi \rb_{\b L^2} ds \\
		&= \alpha \int_0^t \la \lb |\nabla \m|^2 \m, \psi \rb_{\b L^2} - \lb \nabla \m, \nabla \psi \rb_{\b L^2} \ra ds \\
		&\quad + \int_0^t \lb \m \times \nabla \m, \nabla \psi \rb_{\b L^2} ds.
	\end{align*}
	Thus, $ \m $ satisfies \eqref{eq: sLLG weak}, $ \P $-a.s. 
\end{proof}

Therefore, we study the equation \eqref{eq: u} instead of the original SLLG, where PDE tools can be applied directly to the former when $ \omega \in \Omega $ is fixed.

\section{Solution of transformed equation}\label{Section: solution of transformed equation}

\subsection{Maximal solution in $ \b H^2 $}
Fix $ \sigma = 2 $. 
Let $ (\Omega, \mc F, (\mc F_t)_{t \in [0,T]}, \P, W) $ be given as in Section \ref{Section: Wiener process}. 
We follow the approach in \cite[Chapter 8]{lunardi} to show the existence and uniqueness of a maximally defined solution in $ \b H^2 $ of the transformed equation \eqref{eq: u} at a.e.-$ \omega \in \Omega $. 

Let $ \mc E := C^\gamma([0,T];H^1(\T^2; \R^2 \otimes \mathfrak{so}(3))) $ be the parameter space. 
Consider a family of nonlinear equations: for a.e.-$ \omega \in \Omega $,
\begin{equation}\label{eq: v(t,y)}
	\begin{aligned}
		\partial_t \varphi(t,\omega)
		&=H(t,z(\omega),\varphi(t,\omega)), \quad t \in [t_0(\omega),T], \\
		\varphi(t_0,\omega) 
		&= u_0(\omega),
	\end{aligned}
\end{equation}
where 
\begin{equation}\label{eq: H}
	H(t,a,\v)=\alpha \la \Delta \v + |\nabla \v|^2 \v \ra - \v \times\Delta \v + F(t,a,\v),
\end{equation}
for $ t \in [0,T] $, $ a \in \mc E $ and $ \v \in H^1(\T^2; {\b S}^2) $ or $ \v \in \b H^2 $. 

Taking $ t_0 \equiv 0 $ and $ z = A $, then \eqref{eq: v(t,y)} reduces to \eqref{eq: u}. 
Fix $ \omega $, the equations become deterministic, for which we can introduce the notion of a maximal solution without the concept of stopping times. 

\begin{definition}\label{Definition: maximal solution}
	\rm 
	For every $ \omega \in \Omega $ with $ z(\omega) \in \mc E $ and $ \u_0(\omega) \in \b H^2 $, the pair $ (\u,\tau)(\omega) $ depending on $ z $ and $ \u_0 $, is said to be an $ \b H^2 $-maximal solution of equation \eqref{eq: v(t,y)} at $ \omega $ if either 
	\begin{enumerate}[leftmargin=*]
		\item[(a)]
		$ \tau = T $ and $ \u \in C([t_0,T]; \b H^2) $ satisfies \eqref{eq: v(t,y)} in $ [t_0,T] $, or
		
		\item[(b)]
		$ \tau \in (t_0,T] $ and $ \u \in C([t_0,\tau); \b H^2) $ satisfies \eqref{eq: v(t,y)} in $ [t_0,\tau) $ with no continuable extension in $ \b H^2 $.
	\end{enumerate}
	For clarity, we write $ (\u,[t_0,T])(\omega) $ in case (a), and $ (\u,[t_0,\tau)) $ in case (b) or when the case (a) is not known. 
\end{definition}

A maximal solution is well-defined if there exists a unique local solution. 
Local existence and uniqueness is guaranteed if the function $ H $ has sufficiently smooth derivatives and one of them can generate an analytic semigroup. 
We collect properties of $ H $ below. 

From the definitions \eqref{eq: H} and \eqref{eq: F}, $ H : [0,T] \times \mc E \times \b H^2 \to \b L^2 $ is continuous. 
For every $ t \in [0,T] $ and $ \v \in {\b H}^2 $, $ H(t,\cdot,\v): \mc E \to \b L^2 $ is Fr{\'e}chet differentiable  with the derivative 
\begin{equation}\label{eq: Ha}
	H_a(t,a,\v)b = F_a(t,a,\v) b, \quad b \in \mc E. 
\end{equation}
For every $t\in[0,T]$ and $a\in \mc E$, $H(t,a,\cdot):\b H^2 \to\b L^2$ is Fr\'echet differentiable with the derivative
\begin{align}
	H_{\v}(t,a,\v)\varphi
	&=\la \alpha \Delta \varphi - \v \times\Delta \varphi - \varphi \times\Delta \v \ra \label{eq: Hv} \\
	&\quad + \alpha \la |\nabla \v|^2\varphi +2\lb\nabla \v,\nabla \varphi\rb_{\R^3} \v \ra 
	+F_{\v}(t,a,\v)\varphi \nonumber
\end{align}
for any $\varphi\in\b H^2$.  
Let 
\begin{align*}
	\mc E_{\lambda_1}(a)&:= \{b\in \mc E;\,|a-b|_{\mc E}< \lambda_1\}, \\
	\b B_{\lambda_2}(v) &:= \{\varphi\in\b H^2;\,|\varphi-\v|_{\b H^2}<\lambda_2\}, 
\end{align*}
denote the open balls in $ \mc E $ and $ \b H^2 $, centred at $ a \in \mc E $ and $ \v \in \b H^2 $, with radii $ \lambda_1, \lambda_2>0 $, respectively.   
For simplicity, let $ \mc E_{\lambda_1} := \mc E_{\lambda_1}(0) $ and $ \b B_{\lambda_2} := \b B_{\lambda_2}(0) $. 
\begin{lemma}\label{lem-coerc}
	For any $ a \in \mc E_{\lambda_1}$ and $ \v\in \b B_{\lambda_2} $, there exist constants $c_0=c_0(\lambda_1,\lambda_2)>0$ and $c_1=c_1(T,\lambda_1,\lambda_2)>0$, such that for every $t \in [0,T]$ and $ \varphi \in \b H^2 $, 
	\[\lb H_{\v}(t,a,\v)\varphi,\varphi\rb_{\b L^2}\le-c_0|\varphi|_{\b H^1}^2+c_1|\varphi|^2_{\b L^2}\,.\]
\end{lemma}
\begin{proof}
	Recall the components of $ H_{\v}(t,y,\v) $ in \eqref{eq: Hv}. 
	For $ \v \in \b B_{\lambda_2} $, we define the operators $ f_1(\v): \b H^2 \to \b L^2 $ and $ f_2(\v): \b H^2 \to \b L^2 $ by
	\begin{align}
		f_1(\v) \varphi &= \alpha \Delta \varphi - \v \times\Delta \varphi - \varphi \times\Delta \v, \label{eq: I(v)} \\
		f_2(\v) \varphi &= \alpha \la |\nabla \v|^2\varphi +2\lb\nabla \v,\nabla \varphi\rb_{\R^3} \v \ra. \label{eq: J(v)}
	\end{align} 
	Then for $ \varphi \in \b H^2 $, we have
	\begin{align*}
		\lb f_1(\v)\varphi,\varphi\rb_{\b L^2}&=-\alpha|\nabla \varphi|^2_{\b L^2} + \lb \varphi \times\nabla \v,\nabla \varphi\rb_{\b L^2},
	\end{align*}
	where 
	\begin{align*}
		|\lb \varphi \times \nabla \v,\nabla \varphi \rb_{\b L^2}|
		&\le |\nabla \v|_{\b L^4} |\nabla \varphi|_{\b L^2} |\varphi|_{\b L^4} \\
		&\lesssim \lambda_2 \la |\nabla \varphi|_{\b L^2}^\frac{3}{2} |\varphi|_{\b L^2}^\frac{1}{2} + |\nabla \varphi|_{\b L^2} |\varphi|_{\b L^2} \ra \\
		&\le \delta |\nabla \varphi|^2_{\b L^2} + c(\delta) (\lambda_2^4+\lambda_2^2) |\varphi|_{\b L^2}^2\,.
	\end{align*}
	Therefore, for $ \v \in \b H^2 $ and a certain small $ \delta \in (0,1) $,
	\begin{align*}
		\lb f_1(\v)\varphi,\varphi\rb_{\b L^2} &\le (-\alpha+\delta)|\nabla \varphi|_{\b L^2}^2+c(\lambda_2,\delta)|\varphi|_{\b L^2}^2\,.
	\end{align*}
	Similar arguments 
	show that 
	\begin{align*}
		|\lb f_2(\v)\varphi, \varphi\rb_{\b L^2}| &\le \delta|\nabla \varphi|^2_{\b L^2}+c(\lambda_2,\delta)|\varphi|^2_{\b L^2}\,.
	\end{align*}
	For the last term in \eqref{eq: Hv}, it is easy to see that 
	\begin{equation}\label{eq: |Fv|}
	\begin{aligned}
		|F_{\v}(t,a,\v)\varphi|&\lesssim \la |a(t)| + |a(t)|^2 + |\nabla \cdot a(t)|^2 \ra \\
		&\quad \ \cdot f(|v|,|\nabla v|) \la |\varphi|+ |\nabla\varphi| \ra \,,
	\end{aligned}
	\end{equation}
	where 
	$f(\xi,\eta)$ is a linear combination of monomials of the form $\xi^p\eta^q$ with $p\le 2$ and $q\le 1$. 
	Since $ |a|_{\mc E} \leq \lambda_1 $ and $|v|_{\b H^2} \leq \lambda_2$, there exists a constant $ c(\lambda_1,\lambda_2)>0$ such that
	\begin{align*}
		|\lb F_{\v}(t,y,\v)\varphi,\varphi\rb_{\b L^2}|
		&\le c(\lambda_1,\lambda_2) \la \delta |\nabla\varphi|^2_{\b L^2}+ c(\delta)|\varphi|^2_{\b L^2} \ra .
	\end{align*} 
	Then we choose a sufficiently small $\delta \in(0,1)$ such that $c(\lambda_1, \lambda_2)\delta + 2 \delta < \alpha$. 
\end{proof}

\begin{lemma}\label{lem-Lip}
	Fix $ \lambda_1, \lambda_2>0 $. 
	There exists $c=c(T,\lambda_1,\lambda_2)>0$ such that for any $ s,t \in [0,T] $, $ a \in \mc E_{\lambda_1} $ and $v, w \in \b B_{\lambda_2}$, 
	\begin{align*}
		&\Vert H_{a}(t,a,\v) - H_{a}(s,a,\w) \Vert_{\mc E \to \b L^2} + \Vert H_{v}(t,a,\v) - H_{v}(s,a,\w) \Vert_{\b H^2 \to \b L^2} \\
		&\leq c \la |t-s|^\gamma +  |\v -\w|_{\b H^2} \ra. 
	\end{align*}
\end{lemma}
\begin{proof}
	Recall the derivative $ H_a(t,a,\v) $ in \eqref{eq: Ha} and $ F $ in \eqref{eq: F}. 
	We have
	\begin{align*}
		H_a(t,a,\v)b
		&= \alpha h_1(t,a,b,\v) - \v \times h_1(t,a,b,\v) + 2\alpha h_2(t,a,b,\v),
	\end{align*}
	where
	\begin{align*}
		h_1(t,a,b,\v) &= b(t) \cdot \nabla \v + (\nabla \cdot b(t))\v + 2(a(t) \cdot b(t))\v, \\
		h_2(t,a,b,\v) &= \lb a(t) \v, b(t) \v \rb \v + \lb \nabla \v, b(t) \v \rb \v,
	\end{align*}
	for $ t \in [0,T] $, $ a \in \mc E_{\lambda_1} $, $ \v \in \b B_{\lambda_2} $ and any $ b \in \mc E $. 
	Then for $ i=1,2 $, $ h_i $ is $ \gamma $-H{\"o}lder continuous on $ [0,T] $ for every $ \v \in \b B_{\lambda_2} $ due to the time regularity of $ a, b \in \mc E $. 
	Also, $ h_i $ is linear in $ \nabla \v $ and (at most) cubic in $ \v $, such that for $ \v, \w \in \b B_{\lambda_2} $, 
	\begin{align*}
		&|h_i(t,a,b,\v)-h_i(t,a,b,\w)|_{\b L^2} \\
		&\lesssim (1+|a(t)|_{\b H^2})|b(t)|_{\b H^2} \la 1+ |\v|_{\b H^2}^2 + |\w|_{\b H^2}^2 \ra |\v-\w|_{\b H^2} \\
		&\leq c(T,\lambda_1,\lambda_2) |b|_{\mc E} |\v-\w|_{\b H^2}.
	\end{align*} 
	These imply that $ H_a(t,a,\v) $ is $ \gamma $-H{\"o}lder continuous in $ t $ and locally Lipschitz in $ \v $.   
	For the derivative $ H_{\v}(t,y,\v) $, 	
	its first two components $ f_1 $ in \eqref{eq: I(v)} and $ f_2 $ in \eqref{eq: J(v)} do not depend on $ t $ explicitly and they consist of terms linear in $ \v $ and $ \Delta \v $ and quadratic in $ \nabla \v $,
	which imply the local Lipschitz property: for $ \v,\w \in \b B_{\lambda_2} $,
	\begin{align*}
		\Vert f_1(\v) - f_1(\w) \Vert_{\b H^2 \to \b L^2} + \Vert f_2(\v) - f_2(\w) \Vert_{\b H^2 \to \b L^2} \le c(\lambda_2) |\v-\w|_{\b H^2}.  
	\end{align*}
	For the third component $ F_{\v}(t,a,\v) $, it is linear in $ \nabla \cdot a(t) $, quadratic in $ a(t) $, and by \eqref{eq: |Fv|}, linear in $ \nabla \v $ and quadratic in $ \v $. 
	Then similarly, for any $ a \in \mc E_{\lambda_1} $, $ F_{\v}(t,a,\v) $ is $ \gamma $-H{\"o}lder continuous on $ [0,T] $ for every $ \v \in \b B_{\lambda_2} $, and Lipschitz on $ \b B_{\lambda_2} $ for every $ t \in [0,T] $. 
\end{proof}

With the above properties of $ H $, it follows immediately from \cite[Theorem 8.1.1]{lunardi} that for a.e.-$ \omega \in \Omega $, there exists a unique $ \b H^2 $-local solution $ \u_\zeta $ of \eqref{eq: v(t,y)}, where 
\begin{align*}
	\u_\zeta = \u_\zeta(\omega) \in C([t_0(\omega),\zeta];\b H^2) \cap C_\gamma^\gamma((t_0(\omega),\zeta]; \b H^2)
\end{align*}
for some $ \zeta = \zeta(\omega) \in (t_0(\omega),T] $ that depends on $ a(\omega) $ and initial data $ u_0(\omega) $. 
For any $ t_1,t_2 \in [0,\infty) $, the space $ C^\gamma_\gamma((t_1,t_2]; \b H^2) $ is consisted of bounded and $ \gamma $-H{\"o}lder continuous functions $ f: (t_1,t_2] \to \b H^2 $ with the map $ s \mapsto s^\gamma f(s) $ in $ C^\gamma((t_1,t_2]; \b H^2) $. 
Then the maximal solution $ (\u,[t_0,\tau))(\omega) $ is constructed as follows:
\begin{itemize}
	\item 
	$ \u(t,\omega) = \u_\zeta(t)(\omega) $ for $ t\in [t_0(\omega),\zeta] $,
	
	\item
	$ \tau(\omega) $ is the supremum of all $ \zeta > t_0 $ such that there exists a unique $ \b H^2 $-local solution in $ [0,\zeta] $. 
\end{itemize} 
We summarise the result in the following theorem, where part (b) holds by \cite[Corollary 8.3.3]{lunardi}.

\begin{theorem}\label{Theorem: Lunardi maximal soln}
	Fix $ \omega \in \Omega $. 
	Assume that \smash{$ z(\omega) \in \mc E_{\lambda_1} $} and \smash{$ \u_0(\omega)\in \b B_{\lambda_2} $} with \smash{$ |\u_0(x,\omega)| =1 $} for all $ x \in \T^2 $ for some $ \lambda_1 = \lambda_1(\omega) > 0 $ and $ \lambda_2 = \lambda_2(\omega) > 0 $.   
	\begin{enumerate}[leftmargin=*,parsep=5pt]
		\item[(a)] 
		There exists a unique $ \b H^2 $-maximal solution $ \la \u, [t_0,\tau)\ra (\omega) $ of \eqref{eq: v(t,y)} in the sense of Definition \ref{Definition: maximal solution}, where $ \u $ and $ \tau $ depend on $ (z,u_0)(\omega) $. 
		
		\item[(b)]
		Let $ D_{\lambda_1,\lambda_2}(\omega) := \{(t,z,\u_0)(\omega) : (z,\u_0) \in \mc E_{\lambda_1} \times \b B_{\lambda_2}, t \in [t_0,\tau(z,\u_0)) \} $. 
		Then the solution map 
		$ D_{\lambda,\rho}(\omega) \ni (t,z,\u_0)(\omega) \mapsto \u(t; z,\u_0)(\omega) \in \b H^2 $ is continuous.
	\end{enumerate}
\end{theorem}
Taking a non-random initial data $ \u_0 \in \b H^2 $ at $ t_0=0 $ with $ \lambda_2 = |\u_0|_{\b H^2} $ and an $ \mc E $-valued stochastic process $ z = A $ with $ \lambda_1 = |A|_{\mc E} $, by Theorem \ref{Theorem: Lunardi maximal soln} the equation \eqref{eq: u} admits a unique $ \b H^2 $-maximal solution.  
In the case $ \tau < T $, we require an iteration scheme to construct a global solution (see Section \ref{Section: well-posedness}). 
To this end, we consider a more general setting: random initial data $ \u_0 $ at a stopping time $ \tau_0 $. 

\begin{corollary}\label{Coro: u S2}
	Assume that $ \tau_0 \in [0,T) $, $ \P $-a.s. is an $ \b F $-stopping time.  
	Let $ \u_0: \Omega \to H^2(\T^2; \St) $ be $ \mc F_{\tau_0} $-measurable, and $ z = A \in \mc E $, $ \P $-a.s. 
	Then for a.e-$ \omega \in \Omega $, there exists a unique $ \b H^2 $-maximal solution $ (\u, [\tau_0,\tau))(\omega) $ of \eqref{eq: v(t,y)}, satisfying $ |\u(t,x,\omega)|=1 $ for a.e.-$ (t,x) \in [\tau_0,\tau)(\omega) \times \T^2 $. 
\end{corollary}
\begin{proof}	
	Fix $ \omega \in \Omega $ such that $ A(\omega) \in \mc E $. 
	For simplicity, in this proof we write $ f $ in place of $ f(\omega) $ for any random variable $ f $. 
	Let $ \phi(\u(t)) = \frac{1}{4}|1-|\u(t)|^2|_{\b L^2}^2 $ for $ t \in [\tau_0,\tau) $. 
	We have
	\begin{align*}
		\phi(\u) &\lesssim 1 + |\u|_{\b L^2}^2 + |\u|_{\b L^4}^4 < \infty, \\
		\phi'(\u)f &= -\lb (1-|\u|^2) \u,f \rb_{\b L^2}, \quad f \in \b L^2. 
	\end{align*} 
	Then by \eqref{eq: defn gradA}, \eqref{eq: defn D2A} and integration-by-parts, 
	\begin{align*}
		\frac{d}{dt}\phi(\u)
		&= -2\alpha \int_{\T^2} \lb \u, \nabla \u \rb^2 dx
		+ \alpha \int_{\T^2} (1-|\u|^2)^2 |\nabla_{A} \u|^2 dx.
	\end{align*}
	For the second integral, applying the interpolation inequality to $ |1-|\u|^2|_{\b L^4} $, we have
	\begin{align*}
		&\int_{\T^2} (1-|\u|^2)^2 |\nabla_{A} \u|^2 dx \\
		&\lesssim |\nabla_{A} \u|_{\b L^4}^2 \la |1-|\u|^2|_{\b L^2} \int_{\T^2} \lb \u, \nabla \u\rb^2 dx + |1-|\u|^2|_{\b L^2}^2 \ra \\
		&\leq \delta \int_{\T^2} \lb \u, \nabla \u\rb^2 dx + c(\delta) (1+ |\nabla_{A} \u|_{\b L^4}^4) |1-|\u|^2|_{\b L^2}^2. 
	\end{align*}
	Hence, for a sufficiently small $ \delta >0 $, by Gronwall's inequality, 
	\begin{align*}
		\phi(\u(t)) 
		\le \phi(\u_0) e^{c \int_{\tau_0}^\tau \la 1+|\nabla_{A} \u(s)|_{\b L^4}^4 \ra ds}, \quad 
		t \in [\tau_0,\tau),
	\end{align*}
	where
	$ |\nabla_{A} \u|_{\b L^4} \lesssim |A|_{\b H^1} + |\nabla \u|_{\b L^2}^\frac{1}{2} |\u|_{\b H^2}^{\frac{1}{2}} < \infty $ on $ [\tau_0,\tau) $ by Sobolev embedding, and $ \phi(\u_0) = 0 $ yields the result.  
\end{proof}

Next, we verify measurability of the maximal time and the solution process. 
With an abuse of notation, we set $ \u = 0 $ and $ \tau = T $ on $ \P $-null sets. 
To prepare for shifted processes, we define 
\begin{align*}
	\theta_{\tau_0}(s) := (\tau_0 + s) \wedge T, \quad s \geq 0 .
\end{align*}
\begin{remark}
	By the flow property of $ \Y $,  
	\begin{align*}
		\Y(\tau_0+s)
		&=[\Y(\tau_0+s)\Y^{-1}(\tau_0)] \Y(\tau_0)\\
		&=\mc X(s)\Y(\tau_0),
	\end{align*}
	where $\mc X$ is an independent copy of $\Y$, and 
	\[\begin{aligned}
		|A(\tau_0+s)|=|\nabla \Y(\tau_0+s)|
		&\le|\nabla\mc X(s)||\Y(\tau_0)|+|\mc X(s)||\nabla \Y(\tau_0)|\\
		&\le |\nabla\mc X(s)|+|\nabla \Y(\tau_0)|.
	\end{aligned}\]
	Therefore, given $ \sup_{s \in [0,\tau_0]} |A(s)|_{\b H^2} \leq \lambda $, $ \P $-a.s. and an independent copy $ \tau_1 $ of $ \tau_0 $, we have
	\begin{align*}
		\sup_{s \in [0,\tau_1]} |A(\theta_{\tau_0}(s))|_{\b H^2} - |A(\tau_0)|_{\b H^2} \leq \lambda, 
		\quad \P\text{-a.s.}
	\end{align*}
\end{remark}

\begin{corollary}\label{Coro: tau stopping} 
	Let $ \tau_0 $, $ \u_0 $ and $ (\u,\tau) $ be given as in Corollary \ref{Coro: u S2} with $ z=A $.  
	Then, 
	\begin{enumerate}[(i),leftmargin=*]
		\item 
		$ \tau $ is an $\b F$-stopping time, 
		
		\item 
		$ \u(\theta_{\tau_0}(s) \wedge \zeta) $ is $ {\mc F}_{\tau_0+s} $-measurable, for any $ s \geq 0 $ and $ \b F $-stopping time $ \zeta $ such that $ \zeta \in [\tau_0,\tau) $, $ \P $-a.s.
	\end{enumerate}
\end{corollary}
\begin{proof}
	For every $ n \geq 1 $, let $ \psi_n: [0,\infty] \to [0,1] $ be a smooth cut-off function, with $ \psi_n(y) = 1 $ for $ y \in [0,n] $ and $ \psi_n(y) = 0 $ for $ y \geq 2n $. 
	Then for a.e.-$ \omega \in \Omega $, consider the following approximating equation 
	\begin{equation}\label{eq: un}
		\begin{aligned}
			\partial_t \u_n(t,\omega) 
			&= \psi_n(|\u_n(t)|_{\b H^2}) H(t,A,\u_n(t))(\omega), 
			\quad t \in [\tau_0(\omega),T], \\
			\u_n(\tau_0,\omega) 
			&= \u_0(\omega), 
		\end{aligned}
	\end{equation}
	where $ H(t,A,\u_n(t)) $ is a function of $ A(t) $ and $ \u_n(t) $. 
	For every $ \v \in \b H^2 $, 
	the maps 
	\begin{align*}
		\v &\mapsto \psi(|\v|_{\b H^2})\la \alpha \Delta - \v \times \Delta \ra \in \mathcal{L}(\b H^2,\b L^2), \\
		\v &\mapsto f(t,\v) = \psi(|\v|_{\b H^2}) F(t,A(\omega),\v) \in \b H^1,
	\end{align*}
	are locally Lipschitz, and 
	the operator $ \b H^2 \ni \varphi \mapsto \psi_n(|\v|_{\b H^2}) \la \alpha \Delta \varphi - \v \times \Delta \varphi \ra \in \b L^2 $ is sectorial. 
	By \cite[Theorem 6.3]{Amann}, there exists a unique maximal solution $ (\u_n,[\tau_0,\tau_n)) $ with values in $ C([\tau_0,t];\b H^2) $ for every $ t \in [\tau_0,\tau_n) $. 
	In fact, since $ \sup_{t \in [\tau_0,\tau_n)}|\u_n(t)|_{\b H^2} \leq 2n $ from construction, we have $ \tau_n = T $ and a global $ \b H^2 $-solution $ (\u_n,[\tau_0,T]) $, see \cite[Theorem 7.2]{Amann}. 
	
	For the solution $ \u_n = \u_n(A) $, varying the initial data and the parameter $ A $ (e.g. $ A_1 = A(\omega) $ and $ A_2 = A(\omega') $), it is not difficult to deduce that
	\begin{align*}
		&|\u_n(t;A_1,\u_{0,1}) - \u_n(t;A_2,\u_{0,2})|_{\b L^2}^2 \\
		&\leq c(t,n) \la |\u_{0,1} - \u_{0,2}|_{\b L^2}^2 + \int_{\tau_0}^t |A_1 - A_2|_{\b L^4}^2(s) \ ds \ra, 
	\end{align*}
	where the constant $ c(t,n) $ depends on $ |A_1|_{C([\tau_0,t];\b H^2)} $ and $ |A_2|_{C([\tau_0,t];\b H^2)} $,  
	for any $ t\in [\tau_0,T] $. 
	Thus, for a.e.-$ \omega \in \Omega $ and for any constant $ \lambda>0 $, the map
	\begin{align*}
		[0,T] \times {\mc E}_\lambda \times \b H^2 \ni (t,A,\u_0)(\omega) \mapsto \u_n(\theta_{\tau_0}(t); A,\u_0)(\omega) \in \b L^2
	\end{align*}
	is continuous. 	
	
	Let $ A^t(\cdot) := A(\cdot \wedge t) $ denote the stopped process at time $ t \in [0,T] $, which is still progressively measurable with $ A^t(s) \in \mc F_t \cap \mc F_s $ for all $ s\in[0,T] $.
	Let 
	\begin{align*}
		\tau_{A,\lambda} 
		&:= \inf\{ t \in [\tau_0,T]: |A(t)|_{\b H^2} > \lambda \}
	\end{align*}
	for any constant $ \lambda > 0 $. 
	Then $ \tau_{A,\lambda} $ is an $ \b F $-stopping time and 
	\begin{align*}
		\lim_{\lambda \to \infty} \tau_{A,\lambda} 
		&= T.
	\end{align*} 
	
	Since the evolution of $ \u_n(\omega) $ up to time $ t $ depends only on $ \u_0(\omega) $ and $ A^t(\omega) $, we have
	\begin{align*}
		\u_n(t;A,\u_0)(\omega) 
		&= \u_n(t; A^t,\u_0)(\omega), \quad t \in [\tau_0(\omega),T],\\	
		\u_n(\theta_{\tau_0}(s); A^{\tau_{A,\lambda}},\u_0)(\omega) 
		&= \u_n(\theta_{\tau_0}(s); A^{\theta_{\tau_0}(s) \wedge \tau_{A,\lambda}}, \u_0)(\omega), \quad s \geq 0. 
	\end{align*}	
	By the continuity of the solution map and the boundedness of $ A $ on $ [0,\tau_A^\lambda] $,
	\begin{align*}
		\{ \omega: \u_n(\theta_{\tau_0}(s); A^{\tau_{A,\lambda}}(\omega),\u_0(\omega)) \in B \} 
		\in {\mc F}_{\tau_0+s}, \quad \forall B \in \mathcal{B}(\b L^2).
	\end{align*}
	This implies that the $ \b H^2 $-valued process $ \{ \u_n(\theta_{\tau_0}(s); A^{\tau_{A,\lambda}},\u_0) : s \geq 0 \} $ is progressively measurable with respect to $ (\mc F_{\tau_0+s})_{s \geq 0} $. 
	As a result, for any constant $ \rho>0 $, 
	\begin{align*}
		\tau_{n,\rho}^\lambda
		&:= \la \inf\{ s \geq 0 : |\u_n(\theta_{\tau_0}(s); A^{\tau_{A,\lambda}},\u_0)|_{\b H^2} > \rho \} + \tau_0 \ra \wedge \tau_{A,\lambda} \\
		&= \inf\{ t \in [\tau_0,T] : |\u_n(t; A^{\tau_{A,\lambda}},\u_0)|_{\b H^2} > \rho \} \wedge \tau_{A,\lambda}
	\end{align*}
	is an $ {\b F} $-stopping time. 
	
	By the uniqueness of the $ \b H^2 $-maximal solution, $ \u(t; A,\u_0) = \u_n(t; A,\u_0) $ in $ \b H^2 $ for $ t \in [\tau_0, \tau_{n,n}^\lambda] $, $ \P $-a.s. 
	This implies that $ \u(\theta_{\tau_0}(s) \wedge \tau_{n,n}^\lambda; A,\u_0) \in \mc F_{\tau_0+s} $ for $ s \geq 0 $, and 
	\begin{align*}
		\tau_{n,n}^\lambda 
		&= \inf\{ t \in [\tau_0,T] : |\u(t;A^{\tau_{A,\lambda}},\u_0)|_{\b H^2} > n \} \wedge \tau_{A,\lambda} \\
		&= \inf\{ t \in [\tau_0, \tau_{A,\lambda}] : |\u(t;A,\u_0)|_{\b H^2} > n \} \\
		&=: \tau_n^\lambda \leq \tau.
	\end{align*}
	Then we have $ \lim_{n \to \infty} \lim_{\lambda \to \infty} \tau_n^\lambda = \tau $, $ \P $-a.s. and thus $ \tau $ is an $ \b F $-stopping time. 
\end{proof}

The rest of this section is devoted to the (pathwise) regularity of the solution $ u $. 

Heuristically, given a solution $ \u=\u(\omega) $ of \eqref{eq: u} in $ C([\tau_0,\tau);\b H^2) $, 
we have $ \u \in L^2(\tau_0,t; \b H^2) $, $ \partial_t u \in C([\tau_0,t]; \b L^2) $ and $ \nabla \u \in C([\tau_0,t]; \b L^p) $, leading to $ F(t,A,\u) \in L^p(\tau_0,t; \b L^p) $, for any $ t \in (\tau_0,\tau) $ and $ p \in [1,\infty) $.
We can improve the regularity of $ \u $ in $ [\tau_0,t] $ by noticing that \eqref{eq: u} is strongly parabolic and applying classic regularity results in \cite[Chapter VII]{Ladyzhenskaya_book} together with bootstrapping. 
This argument is well-studied (see \cite{GuoHong, Struwe_1985}), and follows similarly for our transformed equation provided that $ A = A(\omega) $ and $ \u_0 $ are sufficiently regular. 
In particular, if $ A \in C([\tau_0,T]; \b H^2) \cap L^2(\tau_0,T; \b H^{\sigma+1}) $ and $ \u(\rho) \in \b H^{\sigma+1} $ for some $ \rho \in [\tau_0,\tau) $ and $ \sigma \geq 1 $, then we have $ \u \in H^1(\rho,t; \b H^{\sigma}) \cap L^2(\rho,t; \b H^{\sigma+2}) $ for any $ t \in (\rho,\tau) $. 
Suppose that these high-order estimates are uniform in $ t $, then using the continuous embeddings such as \smash{$ H^1(\rho,t; \b H^{\sigma}) \hookrightarrow C^{\gamma}([\rho,t]; \b H^{\sigma}) $} for $ \gamma \in (0,1) $, we can extend the uniform continuity of $ \u $ to the interval $ [\rho, \tau) $ and thereby extending $ \u $ continuously to $ \tau $, violating the maximality of $ \tau $ when $ \sigma =2 $. 
This outlines the reason behind possible singularities in the solution that we construct in Theorem \ref{Theorem: Struwe u}, namely, smallness of local energy that ensures uniform (in time) higher-order bounds of $ \u $ cannot hold at $ \tau $. 

Next, we collect the aforementioned estimates of $ \u $: the case $ \sigma = 0 $ in Section \ref{Section: energy estimates} and the case $ \sigma =1,2 $ in Section \ref{Section: regularity small energy}.

\subsection{Energy estimates}\label{Section: energy estimates} 
For a.e.-$ \omega \in \Omega $, 
let $ (\u,[\tau_0,\tau))(\omega) $ be the unique $ \b H^2 $-maximal solution in Corollary \ref{Coro: u S2}. 
Its covariant derivative defines a gauged Dirichlet integral   
\[
E(\u, A)(\omega) = \frac{1}{2} \int_{\T^2} |\nabla_A \u|^2(\omega) \ dx,
\]
which is the governing energy depending on time through the process $A$. 
Roughness in time, however, limits its relevance for an exact energy law in the pathwise sense. 
On the other hand, $ m= Y u$ is (formally) a weak pathwise solution to the stochastic LLG equation and by \eqref{eq: defn gradA}, 
\begin{align*}
	E(\u, A)(\omega) 
	= \frac{1}{2} \int_{\T^2} |\nabla \m|^2(\omega) \ dx
\end{align*}  
has under appropriate conditions a mean energy bound from It\^o calculus. 

In the estimates below, we denote by $ B_r(x) $ the open ball in $ \R^2 $ of radius $ r $ centred at $ x $, and let $ B_r:= B_r(0) $. 

\begin{lemma}\label{Lemma: local energy} 
	Fix $ \omega \in \Omega $ such that $ A(\omega) \in C([0,T]; \b H^2) $. 
	For every $ t \in [\tau_0,\tau)(\omega) $ and $ r = r(\omega) \in (0,1) $, 
	there exists a constant $ c>0 $ such that
	\begin{enumerate}[leftmargin=*]
		\item[(i)]
		$ |\nabla \u(t)|_{\b L^2}^2 \leq \la |\nabla \u_0|_{\b L^2}^2 + c \int_0^t |\nabla \cdot A+ A^2|_{\b L^2}^2 \ ds \ra e^{c \int_0^t |A(s)|_{\b L^\infty}^2 ds}$,
		
		\item[(ii)]
		$ |\nabla \u(t)|_{L^2(B_r)}^2 \leq |\nabla \u_0|_{L^2(B_{2r})}^2 + S(\tau_0,t,r,A,\u) $, where 
		\begin{align*}
			&S(\tau_0, t,r,A,\u) \\ 
			&= - \alpha \int_{\tau_0}^t |\partial_s \u(s)|_{L^2(B_{r})}^2 ds + c \int_{\tau_0}^t |\nabla \cdot A(s) + A^2(s)|_{L^2(B_{2r})}^2 ds \\
			&\quad + c \sup_{s \in [\tau_0,t]} |\nabla \u(s)|_{L^2(B_{2r})}^2 \la r^{-2} t+\int_{\tau_0}^t |A(s)|_{L^\infty(B_{2r})}^2 ds \ra.
		\end{align*}
	\end{enumerate}	
\end{lemma}
\begin{proof}
	Applying It{\^o}'s lemma and using \eqref{eq: defn D2A}, 
	\begin{align*}
		\frac{1}{2} \partial_s |\nabla \u|_{\b L^2}^2
		&= \lb - \Delta \u , \u \times \Delta_A \u - \alpha \u \times (\u \times \Delta_A \u) \rb_{\b L^2} \\
		&\leq -\frac{\alpha}{2} |\u \times \Delta_A \u|_{\b L^2}^2 + c(\alpha)|A \cdot \nabla \u + (\nabla \cdot A) \u + A^2 \u|_{\b L^2}^2.
	\end{align*}
	Then part (i) follows from the fact $ |\u(s,x,\omega)|=1 $ a.e. on $ [\tau_0,t] \times \T^2 $, the estimate $ |A \cdot \nabla \u|_{\b L^2} \leq |A|_{\b L^\infty} |\nabla \u|_{\b L^2} $ and Gronwall's lemma. 
	Similarly, for $ \varphi \in C^\infty_0(B_{2r}) $ with $ \varphi \equiv 1 $ on $ B_r $ and $ |\nabla \varphi| \lesssim r^{-1} $ on $ B_{2r} $, 
	\begin{align*}
		\frac{1}{2} \partial_s |\varphi \nabla \u|_{\b L^2}^2
		&\leq -\frac{\alpha}{2} |\varphi \u \times \Delta_A \u|_{\b L^2}^2 + c(\alpha)|A \cdot \nabla \u + (\nabla \cdot A) \u + A^2 u|_{L^2(B_{2r})}^2 \\
		&\quad - 2 \lb \nabla \u, (\varphi \nabla \varphi) (\partial_t \u) \rb_{\b L^2},
	\end{align*}
	where $ |\varphi \partial_s \u|_{\b L^2} \lesssim |\varphi \u \times \Delta_A \u |_{\b L^2} $ on $ [\tau_0,t] $ thanks to the norm constraint. 
	Then part (ii) follows by noticing $ |f|_{L^2(B_r)} \leq |\varphi f|_{\b L^2} $. 	
\end{proof}

\begin{lemma}\label{Lemma: energy moments}
	For any $ \b F $-stopping time $ \zeta \in [\tau_0,\tau) $ and constant $ \gamma \in (0,\frac{1}{2}] $ and $ p \in [1,\infty) $, there exists a constant $ c>0 $ depending only on $ \gamma, p $ and $ T $ such that 
	\begin{align*}
		&\E \left[ \sup_{t \in [0,T]}|\nabla \u(\theta_{\tau_0} \wedge \zeta)|_{\b L^2}^{2p} 
		+ \la \int_{\tau_0}^{\zeta} |\u \times \Delta_A \u|_{\b L^2}^2 \ ds \ra^p \right] \\
		&+ \E \left[ |\u(\theta_{\tau_0} \wedge \zeta)|_{C^\gamma([0,T]; \b L^2)}^{2p} \right] \\
		&\leq c \la \E \left[|\nabla \u(\tau_0)|_{\b L^2}^{2p} + |A|_{C([0,T]; \b L^2)}^{2p} \right] + q^{4p}(2) + 1\ra.
	\end{align*} 
\end{lemma}
\begin{proof}
	By Corollary \ref{Coro: u S2} and \ref{Coro: tau stopping}, the process $ \{\u(\theta_{\tau_0}(s) \wedge \zeta): s \geq 0 \} $ is progressively measurable with respect to $ (\mc F_{\tau_0+s})_{s \geq 0} =: {\b F}_{\tau_0} $ and it takes values in $ \St $ for a.e.-$ (t,x,\omega) $. 
	Then we can apply It{\^o}'s lemma: for any $ t \in [\tau_0,\tau) $, 
	\begin{align*}
		\frac{1}{2}d|\nabla_A \u|_{\b L^2}^2(t) 
		&= -\alpha |\u \times \Delta_A\u|_{\b L^2}^2 \ dt \\
		&\quad + \frac{1}{2} \sum_{k=1}^\infty \lb \nabla_A \u, ((\nabla G_k) G_k-G_k \nabla G_k) \Y \u \rb_{L^2} \ dt \\
		&\quad + \frac{1}{2} \sum_{k=1}^\infty |(\nabla G_k) \Y \u|_{\b L^2}^2 \ dt \\
		&\quad + \sum_{k=1}^\infty \lb \nabla_A \u, (\nabla G_k) \Y \u \rb_{\b L^2} dW_k. 
	\end{align*}
	Estimating the drift part,
	\begin{align*}
		\lb \nabla_A \u, (\nabla G_k) G_k \Y\u \rb_{\b L^2}
		&\leq |\nabla_A \u|_{\b L^2} |\nabla g_k|_{\b L^4} |g_k|_{\b L^4} \\
		&\leq c(\delta) |g_k|_{\b L^4}^2 |\nabla g_k|_{\b L^4}^2 + \delta |\nabla_A \u|_{\b L^2}^2. 
	\end{align*}
	For the diffusion part, we shift the process before applying the Burkholder-Davis-Gundy inequality. 
	Note that $ \{ \tau_0+s \leq \zeta \} \in {\mc F}_{\tau_0+s} $ for every $ s \geq 0 $, and thus $ \zeta-\tau_0 $ is a stopping time with respect to $ (\mc F_{\tau_0+s})_{s \geq 0} =: {\b F}_{\tau_0} $. 
	Since for $ s \in [0,\zeta-\tau_0) $, 
	\begin{align*}
		&\int_{\tau_0}^{\tau_0+s} \lb \nabla_A\u(r), (\nabla G_k) (\Y\u)(r) \rb_{\b L^2} dW_k(r) \\
		&= \int_0^s \lb \nabla_A\u(\tau_0+r), (\nabla G_k) (\Y\u)(\tau_0+r) \rb_{\b L^2} dW_k(\tau_0+r)
	\end{align*}
	and $ \{ W_k(\tau_0+r): r \geq 0\} $ is an $ {\b F}_{\tau_0} $-martingale, 
	we have for $ p \in [1,\infty) $, 
	\begin{align*}
		&\E \left[ \sup_{s \in [0,\zeta-\tau_0]} \left| \sum_{k=1}^\infty \int_0^s \lb \nabla_A\u(\tau_0+r), (\nabla G_k) (\Y\u)(\tau_0+r) \rb_{\b L^2} \  dW_k(\tau_0+r) \right|^p \right] \\
		&\leq \sum_{k=1}^\infty \E \left[ \la \int_0^{\zeta-\tau_0} \sum_{k=1}^\infty \lb \nabla_A\u(\tau_0+r), (\nabla G_k) (\Y\u)(\tau_0+r) \rb_{\b L^2}^{2}\ dr \ra^\frac{p}{2} \right] \\
		&\leq c \la \sum_{k=1}^\infty |\nabla g_k|_{\b L^2}^2 \ra^p \E \left[ |\zeta-\tau_0|^p \right] + \delta \E \left[ \la \int_{\tau_0}^\zeta |\nabla_A \u|_{\b L^2}^2(r) \ dr \ra^p \right],
	\end{align*}
	where $ |\zeta-\tau_0| \leq 2T $, $ \P $-a.s. by construction. 
	Therefore, 
	\begin{align*}
		&\E \left[ \sup_{s \geq 0}|\nabla_A \u(\theta_{\tau_0}(s) \wedge \zeta)|_{\b L^2}^{2p} + \la \int_{\tau_0}^{\zeta} |\u \times \Delta_A \u|_{\b L^2}^2 \ ds \ra^p \right] \\
		&\leq \E \left[ |\nabla_A \u(\tau_0)|_{\b L^2}^{2p} \right] + c q^{4p}(2),
	\end{align*}
	where $ q(\sigma) $ is given in \eqref{eq_q3}. 
	
	As a result, together with the fact that $ \u $ takes values in $ \St $ a.e. (see Corollary \ref{Coro: u S2}), we obtain the moment estimate in $ C^\gamma([0,T]; \b L^2) $ and in $ C([0,T];\b H^1) $. 
	In more details, for any $ p \in [1,\infty) $, we have for the H{\"o}lder seminorm 
	\begin{align*}
		&\E \left[ \sup_{0\leq s_1<s_2 \leq T}\frac{\left|\u(\theta_{\tau_0}(s_2) \wedge \zeta)-\u(\theta_{\tau_0}(s_1) \wedge \zeta) \right|_{\b L^2}^{2p}}{|s_2-s_1|^p} \right] \\
		&\leq c \E \left[ \sup_{0\leq s_1<s_2 \leq T}\frac{\left|\int_{\theta_{\tau_0}(s_1) \wedge \zeta}^{\theta_{\tau_0}(s_2) \wedge \zeta} \u \times \Delta_A \u \ ds \right|^{2p}_{\b L^2}}{|s_2-s_1|^p} \right] \\
		&\leq c \E \left[ \sup_{0\leq s_1<s_2 \leq T} \la \int_{\theta_{\tau_0}(s_1) \wedge \zeta}^{\theta_{\tau_0}(s_2) \wedge \zeta} |\u \times \Delta_A \u|_{\b L^2}^2 \ ds \ra^p \right] \\
		&\leq c \E \left[ \la \int_{\tau_0}^{\zeta} |\u \times \Delta_A \u|_{\b L^2}^2 \ ds \ra^p \right] \\
		&\leq c\la \E \left[ |\nabla_A \u(\tau_0)|_{\b L^2}^{2p} \right] + q^{4p}(2) \ra,
	\end{align*}
	and for the $ \b L^2 $-norm of the gradient, 
	\begin{align*}
		&\E \left[ \sup_{s \geq 0} \left|\nabla \u(\theta_{\tau_0}(s) \wedge \zeta) \right|_{\b L^2}^{2p} \right] \\
		&\leq c(p) \E \left[ \sup_{s \geq 0}|\nabla_A \u(\theta_{\tau_0}(s) \wedge \zeta)|_{\b L^2}^{2p} + |A|_{C([0,T];\b L^2)}^{2p} \right] \\
		&\leq c(p) \la \E \left[ |\nabla_A \u(\tau_0)|_{\b L^2}^{2p} + |A|_{C([0,T];\b L^2)}^{2p} \right] + q^{4p}(2) \ra, 
	\end{align*}
	where $ |\nabla_A \u(\tau_0)|_{\b L^2}^2 \leq \frac{1}{2}|\nabla \u(\tau_0)|_{\b L^2}^2 + \frac{1}{2}|A(\tau_0)|_{\b L^2} $, 
	and $ A $ is uniformly bounded in $ L^{2p}(\Omega; C([0,T]; \b L^2)) $ by Lemma \ref{Lemma: Y flow}. 
\end{proof}

\subsection{Regularity under small local energy}\label{Section: regularity small energy}

\begin{hypothesis}\label{Hypo: small |Du|+|A|}
	Let $ \v \in L^\infty(0,t; \b H^1) $, $ \P $-a.s. 
	For a.e.-$ \omega \in \Omega $, there exist a radius $ r = r(\omega) \in (0,1) $, a non-empty interval $ I=I(\omega) \subseteq [0,t(\omega)] $, and a constant $ \eps \in (0,1) $, such that for some $ x_0 \in \T^2 $, 
	\begin{equation*}
		\sup_{s \in I} \la |\nabla \v(s,\omega)|_{L^2(B_{2r}(x_0))}^2 + |A(s,\omega)|_{L^2(B_{2r}(x_0))}^2 \ra \le \eps^2.
	\end{equation*}
\end{hypothesis}

\subsubsection{$ L^4 $-estimates}\label{Section: L4 estimates}	
The key ingredient is the local form of Ladyzhenskaya's interpolation inequality which follows from the Sobolev inequality $|h|_{L^2(\R^2)} \lesssim |\nabla h|_{L^1(\R^2)}$ applied to functions of the form $h=\varphi |f|^2$ for $\varphi \in C^\infty_0(\R^2; \R)$ and $f \in H^1(\R^2)$, yielding 
\begin{equation}\label{eq: lady}
	\int_{B_{r}} |f|^4 \varphi^2 \ dx  \le c_0 |f|^2_{L^2(B_{r})}  \int_{B_{r}} \la |\nabla f|^2 \varphi^2 + |f|^2 |\nabla \varphi|^2 \ra \ dx\,,
\end{equation}
for $\varphi$ supported on $ B_{r} $, where $ r>0 $ and the constant $ c_0 $ is independent of $ r $ and $ f $. 

In the following, we collect $ L^4 $ and subsequently $ H^1 $ estimates of $ A $, $ \nabla A $ and $ \nabla \u $ using \eqref{eq: lady} under Hypothesis \ref{Hypo: small |Du|+|A|}. 

For estimates of $ A $, we first observe that using \eqref{eq:div_curl}, integration by parts and Young's inequality, 
\begin{equation}\label{eq: |varphi DA| L2}
	\begin{aligned}
		|\varphi \nabla A|^2_{\b L^2} 
		&\lesssim \int_{B_r} \varphi^2 \la |\nabla \cdot A|^2 + |A|^4 \ra + \varphi |\nabla \varphi| |A| |\nabla A| \ dx \\
		&\leq c_1 \la |\varphi (\nabla \cdot A)|_{\b L^2}^2 + |\varphi^\frac{1}{2}A|_{\b L^4}^4 + r^{-2} |A|_{L^2(B_{2r})}^2 \ra, 
	\end{aligned}
\end{equation}
and 
\begin{equation}\label{eq: |varphi D2A| L2}
	\begin{aligned}
		|\varphi \nabla^2 A|_{\b L^2}^2 
		&\lesssim \int_{B_r} \varphi^2 \la |\nabla(\nabla \cdot A)|^2 + |\nabla[A_1, A_2]|^2 \ra + \varphi |\nabla \varphi| |\nabla A| |\nabla^2 A| \ dx \\
		&\leq c_1 \la |\varphi \nabla (\nabla \cdot A)|_{\b L^2}^2 + |\varphi^\frac{1}{2}A|_{\b L^4}^2 |\varphi^\frac{1}{2} \nabla A|_{\b L^4}^2 + r^{-2} |\nabla A|_{L^2(B_{2r})}^2 \ra, 
	\end{aligned}
\end{equation}
for every non-negative $ \varphi \in C_0^\infty(B_{2r}; \R) $ with $ |\nabla \varphi| \lesssim r^{-1} $ on $ B_{2r} $, and some constant $ c_1 $ independent of $ r $ and $ A $.  

\begin{lemma}\label{Lemma: A L4 + H1}
	Assume that $ A $ satisfies Hypothesis \ref{Hypo: small |Du|+|A|} for $ \eps $ sufficiently small such that $ c_0 c_1 \eps^2 \leq c_2 < 1 $ for some constant $ c_2 $ independent of $ \eps $. 
	Then for any $ s \in I $ and $ \varphi \in C_0^\infty(B_{2r}(x_0); \R) $ with $ 0 \leq \varphi \leq 1 $ and $ |\nabla \varphi| \lesssim r^{-1} $ on $ B_{2r}(x_0) $, 
	$ A = A(s,\omega) $ satisfies 
	\begin{enumerate}[leftmargin=*,parsep=0pt]
		\item[(i)]
		$ |\varphi^\frac{1}{2} A|_{\b L^4}^4 \lesssim \eps^2 (|\varphi \nabla \cdot A|_{\b L^2}^2 + r^{-2} \eps^2) $, 
		
		\item[(ii)]
		$ |\varphi \nabla A|_{\b L^2}^2 \lesssim |\varphi \nabla \cdot A|_{\b L^2}^2 + r^{-2} \eps^2 $,
		
		\item[(iii)]
		$ |\varphi^\frac{1}{2} \nabla A|_{\b L^4}^4 \lesssim \Psi(r, A, \eps) + r^{-6} \eps^4 $, 
		
		\item[(iv)]
		$ |\varphi \nabla^2 A|_{\b L^2}^2 \lesssim \Psi(r,A,\eps)+ r^{-2} \la 1 + r^{-4} \eps^4 \ra $,
	\end{enumerate}	 
	where
	\begin{align*}
		\Psi(r, A, \eps) 
		&:= r^2 |\nabla (\nabla \cdot A)|_{L^2(B_{2r}(x_0))}^4 \\
		&\quad + |\nabla \cdot A|_{L^2(B_{4r}(x_0))}^4 \la r^{-2} + \eps^2 |\nabla \cdot A|_{L^2(B_{2r}(x_0))}^2 \ra,
	\end{align*}
	for a.e.-$ \omega \in \Omega $. 
\end{lemma}
\begin{proof}
	For simplicity, let $ x_0 = 0 $. 
	For part (i), taking $ f = A $ in \eqref{eq: lady}, the term \smash{$ |\varphi \nabla A|_{\b L^2}^2 $} appearing on the right-hand side can be estimated by applying \eqref{eq: |varphi DA| L2}. As a result, the term \smash{$ |\varphi^\frac{1}{2}A|_{\b L^4}^4 $} appears on the right-hand side, which can be absorbed into the left provided that $ \eps $ is sufficiently small such that $ 1-c_0 c_1 \eps^2 > 0 $, and then the smallness condition in Hypothesis \ref{Hypo: small |Du|+|A|} yields the estimate. 
	Part (ii) follows from part (i) and \eqref{eq: |varphi DA| L2}. 
	
	For part (iii), taking $ f = \nabla A $ in \eqref{eq: lady} and applying the estimate in \eqref{eq: |varphi D2A| L2}, we have
	\begin{align*}
		|\varphi^\frac{1}{2} \nabla A|_{\b L^4}^4
		&\lesssim
		|\nabla A|_{L^2(B_{2r})}^2 \la |\varphi \nabla (\nabla \cdot A)|_{\b L^2}^2 + r^{-2} |\nabla A|_{L^2(B_{2r})}^2 \ra \\
		&\quad + |\nabla A|_{L^2(B_{2r})}^4 |\varphi^\frac{1}{2}A|_{\b L^4}^4,
	\end{align*}
	where \smash{$ |\nabla A|_{L^2(B_{2r})} \leq |\psi \nabla A|_{\b L^2} $} for another cut-off function $ \psi \in C_0^\infty(\R^2; \R) $ supported on $ B_{4r} $ with $ \psi = 1 $ on $ B_{2r} $ and \smash{$ |\nabla \psi| \lesssim r^{-1} $} on $ B_{4r} $. 
	Then using part (ii) with $ \psi $ in place of $ \varphi $, we obtain the estimate of \smash{$ |\varphi^\frac{1}{2} \nabla A|_{\b L^4} $}. 
	Part (iv) follows from part (iii), \eqref{eq: |varphi D2A| L2} and Young's inequality. 
\end{proof}

\begin{remark} 
	If \smash{$|\nabla \cdot A| \lesssim |A|^2$} (i.e. under the Coulomb gauge), 
	then under the smallness condition \smash{$ |A|_{L^2(B_{2r}(x_0))}^2<\eps^2 $}, one even has a reverse H\"older inequality
	\begin{align*}
		|\varphi^\frac{1}{2} A|_{\b L^4}^4 \lesssim \eps^4 \sup_{x \in B_{2r}(x_0)} |\nabla \varphi(x) |^2, 
	\end{align*}
	for any smooth cut-off function $ \varphi $ supported on $ B_{2r}(x_0) $. 
\end{remark}

Given $ (\u,[\tau_0,\tau))(\omega) $ the unique $ \b H^2 $-maximal solution of \eqref{eq: v(t,y)}, 
assume that $ \u $ satisfies Hypothesis \ref{Hypo: small |Du|+|A|}. 
Similar results hold for $ \nabla \u = \nabla \u(\omega) $ for a.e.-$ \omega \in \Omega $, where
\begin{equation}\label{eq: Du L4}
	|\varphi^\frac{1}{2} \nabla \u(s)|_{\b L^4}^4
	\leq c_0 \eps^2 \la |\varphi \nabla^2 \u(s)|_{\b L^2}^2 + r^{-2} \eps^2 \ra, \quad s \in I.
\end{equation} 
Then taking into account \eqref{eq: |F|} and Lemma \ref{Lemma: A L4 + H1}(i), we deduce below an estimate for $ \varphi \nabla^2 \u $ in $ L^2(0,t;\b L^2) $ and thereby refine the bound in \eqref{eq: Du L4} under time integral. 
\begin{proposition}\label{Prop: Du L4 + H1}
	Assume that $ \u $ and $ A $ satisfy Hypothesis \ref{Hypo: small |Du|+|A|} for $ \eps $ sufficiently small such that $ c_0 (1+c_1) \eps^2 \leq c_2 < 1 $ for some constant $ c_2 $ independent of $ \eps $.  
	Then for any $ \varphi \in C_0^\infty(B_{2r}(x_0); \R) $ with $ 0 \leq \varphi \leq 1 $ and $ |\nabla \varphi| \lesssim r^{-1} $ on $ B_{2r}(x_0) $, the processes $ \u = \u(\omega) $ and $ A = A(\omega) $ satisfy
	\begin{align*}
		&\int_{I} \la \eps^2 |\varphi \Delta \u|^2_{\b L^2} + |\varphi^\frac{1}{2} \nabla \u|_{\b L^4}^4 \ra ds \\
		&\lesssim (1+r^{-2} |I|) \eps^4 + \eps^2 \int_{I} |\nabla \cdot A|_{L^2(B_{2r}(x_0))}^2 \ ds,
	\end{align*}
	for a.e.-$ \omega \in \Omega $. 
\end{proposition}
\begin{proof} 
	For simplicity, let $ x_0 = 0 $ and $ I = [0,t] $. 
	Recall \eqref{eq: du = LL + F}. We have 
	\begin{align*}
		\frac{1}{2}\partial_s |\varphi \nabla \u(s)|_{\b L^2}^2 
		&= -\alpha |\varphi \Delta \u|_{\b L^2}^2 - \lb \varphi^2 \Delta \u, \alpha |\nabla \u|^2 u + F(s,A,\u) \rb_{\b L^2} \\
		&\quad -2 \lb \nabla \u, \varphi (\nabla \varphi) (\partial_s \u) \rb_{\b L^2},
	\end{align*}
	for $ s \in (0,t) $. 
	Here, $ |(\nabla \varphi) (\nabla \u) |_{L^2(B_{2r})} \lesssim r^{-1} \eps $ by Hypothesis \ref{Hypo: small |Du|+|A|}, and under the norm constraint, 
	\begin{align*}
		|\varphi \partial_s \u|_{\b L^2}
		&\leq | \varphi \Delta \u|_{\b L^2}^2 + |\varphi^\frac{1}{2}\nabla \u|^4_{\b L^4} + |\varphi F(s,A,\u)|_{\b L^2}^2,
	\end{align*}
	where by \eqref{eq: |F|} and Lemma \ref{Lemma: A L4 + H1}(i),
	\begin{align*}
		|\varphi F(s,A,\u)|_{\b L^2}^2
		&\leq 
		c(\delta) \la |\varphi (\nabla \cdot A)|_{\b L^2}^2 + r^{-2}\eps^4 \ra
		+ \delta |\varphi^\frac{1}{2} \nabla \u|_{\b L^4}^4
	\end{align*}
	Thus, for any small $ \delta \in (0,1) $, 
	\begin{align*}
		\partial_s |\varphi \nabla \u(t)|_{\b L^2}^2 
		&\leq 
		-(\alpha-\delta) |\varphi \Delta \u|_{\b L^2}^2 
		+ (\alpha+\delta) |\varphi^\frac{1}{2} \nabla \u|_{\b L^4}^4 \\
		&\quad + c(\delta) \la |\varphi (\nabla \cdot A)|_{\b L^2}^2 + r^{-2} \eps^2 \ra,
	\end{align*}
	where $ |\varphi^\frac{1}{2} \nabla \u|_{\b L^4}^4 $ is estimated in \eqref{eq: Du L4}. 
	Thus, for a sufficiently small $ \delta $ (depending only on $ c_2 $ and $ \alpha $), we have 
	\begin{align*}
		\int_0^s |\varphi \Delta \u|_{\b L^2}^2 \ ds
		\lesssim 
		|\varphi \nabla \u_0|_{\b L^2}^2 
		+ r^{-2} \eps^2 s + \int_0^s |\nabla \cdot A|_{L^2(B_{2r})}^2 \ ds,
	\end{align*}
	for $ s \in (0,t) $. 
	The estimates follow by noting $ |\nabla \u_0|_{L^2(B_{2r})} < \eps $, and then applying \eqref{eq: Du L4}. 
\end{proof}

The estimate in Proposition \ref{Prop: Du L4 + H1} becomes particularly strong in the parabolic scaling $t \sim r^2$.
For $r>0$ and a time-space point $z=(t,x)$, let
$ P_r(z):=[t-r^2,t] \times \bar{B_r(x)} $
denote the closed parabolic cylinder and accordingly $P_r:=P_r(0)$. 
Fix $ \omega \in \Omega $. Given a point $ x_0 $, the radius $ r = r(\omega) $ and the paths $ \u = \u(\omega) $ and $ A = A(\omega) $ in Proposition \ref{Prop: Du L4 + H1}, for $ z_0 = (t_0,x_0) \in [r^2,t] \times \T^2 $, we have
\begin{equation} \label{eq: D2u L2 Pr} 
	|\nabla \u|^4_{L^4(P_{r}(z_0))} + \eps^2 |\Delta \u|^2_{L^2(P_{r}(z_0))} 
	\lesssim
	\eps^2 \la \eps^2 + |\nabla \cdot A|^2_{L^2(P_{2r}(z_0))} \ra,	
\end{equation}
where the right-hand side is almost surely finite if $A \in C([0,T];\b H^1)$, $ \P $-a.s., for example, under \eqref{eq_q3} with $ \sigma =2 $.

\subsubsection{Higher-order estimates}\label{Section: higher-order estimates}

\begin{proposition}\label{Prop: u H2 est Pr}
	Assume that $ \u $ and $ A $ satisfy Hypothesis \ref{Hypo: small |Du|+|A|} with $ I = [t_0 - (2r)^2, t_0] $ for some $ t_0 \in ((2r)^2,t] $, and $ \eps $ as in Proposition \ref{Prop: Du L4 + H1}. 
	Let $ z_0 := (t_0, x_0) $ for a given $ x_0 \in \T^2 $, and define
	\begin{align*}
		\Lambda_\gamma := L^\infty(t_0 - \gamma^2 r^2, t_0; L^2(B_{\gamma r}(x_0))) \cap L^2(t_0 - \gamma^2 r^2, t_0; H^1(B_{\gamma r}(x_0))),
	\end{align*}
	for $ \gamma \in (0,2) $. 
	Then for $ \u = \u(\omega) $ and $ A = A(\omega) $, 
	\begin{enumerate}[leftmargin=*,parsep=0pt]
		\item[(i)]
		$ \Delta \u \in \Lambda_{\frac{1}{2}} $, 
		with bound that only depends on $ r^{-1} \eps $, $ r^{-2} |\nabla \cdot A|^2_{L^2(P_{2r}(z_0))} $ and 
		\begin{align*}
			\eps^2 \int_{t_0-r^2}^{t_0} \Psi \la \frac{r}{2},\eps,A \ra \ ds,
		\end{align*}
		
		\item[(ii)]
		$ \nabla \Delta \u \in \Lambda_{\frac{1}{4}} $, 
		with bound that only depends on 
		$ r^{-1}\eps $, 
		$ |\Delta(\nabla \cdot A)|_{L^2(P_\frac{r}{2}(z_0))}^2 $, 
		and 
		\begin{align*}
			\int_{t_0-\frac{r^2}{4}}^{t_0} \la |\nabla \cdot A|_{H^1(B_{r}(x_0))}^4 + \Psi \la \frac{r}{2},\eps,A \ra + \eps^2 \Psi^2 \la \frac{r}{4},\eps,A \ra \ra \ ds,
		\end{align*}
	\end{enumerate}
	for a.e.-$ \omega \in \Omega $. 
\end{proposition} 
\begin{proof}
	For simplicity, we shift time to the left and assume $(t_0, x_0)=0$. 
	For $ \rho>0 $, let $ \phi_\rho \in (0,1) $ be a cut-off function of the form 
	\begin{align*}
		\phi_\rho(t,x) = \eta_\rho(t) \varphi_\rho(x),
	\end{align*}
	where 
	$ \eta_\rho \in C^\infty(\R) $ is non-decreasing with $ \eta_\rho(t) =0 $ for $ t \leq -\rho^2 $ and $ \eta_\rho(t) =1 $ for $ t \geq -\frac{\rho^2}{4} $,  
	and
	$ \varphi_\rho \in C^\infty_0(B_\rho) $ with $ \varphi_\rho \equiv 1 $ on $ B_\frac{\rho}{2} $, and $ 0 \leq \varphi_\rho \leq 1 $ and $ |\nabla \varphi_\rho| \lesssim \rho^{-1} $ on $ B_\rho $. 
	
	{\it Part (i).}
	Take $ \rho = r $. 
	For brevity, we write $ \eta, \varphi $ instead of $ \eta_r, \varphi_r $ for part (i). 
	We have
	\begin{align*}
		\frac{1}{2}\partial_t |\phi \Delta \u |_{\b L^2}^2 
		&= -\lb \nabla(\phi^2 \Delta \u), \nabla \partial_t \u \rb_{\b L^2} + \lb \phi \Delta \u, (\partial_t \phi) \Delta \u \rb_{\b L^2} \\
		&=
		-\alpha |\phi \nabla \Delta \u|_{\b L^2}^2 
		- \alpha I_0 
		- (I_1+I_2) 
		+ I_3.
	\end{align*}
	where 
	\begin{align*}
		I_0 &= \lb (\nabla \phi^2) \Delta \u, \nabla \Delta \u \rb_{\b L^2}, \\
		I_1 &= \lb \nabla (\phi^2 \Delta \u), \alpha \nabla(|\nabla \u|^2 \u) - \nabla (\u \times \Delta \u) \rb_{\b L^2}, \\
		I_2 &= \lb \nabla (\phi^2 \Delta \u), \nabla F(t,A,\u)) \rb_{\b L^2}, \\
		I_3 &= \lb \phi \Delta \u, (\partial_t \phi) \Delta \u \rb_{\b L^2}.
	\end{align*}
	Since $ |\partial_t \phi| \leq |\eta'| \lesssim r^{-2} $, we have
	\begin{align*}
		\alpha |I_0| + |I_3|
		&\lesssim r^{-1} |\phi \nabla \Delta \u|_{\b L^2} |\Delta \u|_{L^2(B_r)} + |\Delta \u|_{L^2(B_r)}^2 |\eta'| \\
		&\leq \delta |\phi \nabla \Delta \u|_{\b L^2}^2 + c(\delta) r^{-2} |\Delta \u|_{L^2(B_r)}^2,
	\end{align*}
	and \eqref{eq: D2u L2 Pr} can be applied to bound $ |\Delta \u|_{L^2(P_r)}^2 $. 
	Then choosing a sufficiently small $ \delta > 0 $, 
	\begin{equation}\label{eq: |D2u(0)|-L2}
		\frac{1}{2} |\phi \Delta \u(t)|_{\b L^2}^2 
		+ \frac{\alpha}{2} |\phi \nabla \Delta \u|_{L^2(P_r)}^2 
		\leq 
		cr^{-2} |\Delta \u|_{L^2(P_{r})}^2 
		- \int_{-r^2}^t ( I_1 + I_2 ) \ ds,
	\end{equation}
	for any $ t \in (-r^2,0) $. 
	We estimate the nonlinear terms	$ I_1 $ and $ I_2 $ below. 
	
	\underline{$ I_1 $ estimate}
	
	By the interpolation inequality
	\begin{equation}\label{eq: phiD2-L4}
		\begin{aligned}
			|\phi \Delta \u|_{\b L^4} 
			&\lesssim 
			|\phi \nabla \Delta \u  + (\nabla \phi) \Delta \u |_{L^2(B_{r})}^\frac{1}{2} |\phi \Delta \u|_{L^2(B_{r})}^\frac{1}{2} + |\phi \Delta \u|_{L^2(B_{r})} \\
			&\lesssim 
			|\phi \nabla \Delta \u|_{\b L^2}^\frac{1}{2} |\phi \Delta \u|_{\b L^2}^\frac{1}{2} 
			+ r^{-\frac{1}{2}} |\Delta \u|_{L^2(B_{r})}^\frac{1}{2} |\phi \Delta \u|_{\b L^2}^\frac{1}{2},
		\end{aligned}
	\end{equation}
	and Young's inequality, we have
	\begin{equation}\label{eq: D321 + D221}
		\begin{aligned}
			&\la |\phi \nabla \Delta \u|_{\b L^2} + |\Delta \u|_{L^2(B_{2r})} \ra |\phi \Delta \u|_{\b L^4} |\nabla \u|_{L^4(B_{r})} \\
			&\leq
			\delta |\phi \nabla \Delta \u|_{\b L^2}^2 
			+ c\la |\nabla \u|_{L^4(B_{r})}^4 |\phi \Delta \u|_{\b L^2}^2 + r^{-2} |\Delta \u|_{L^2(B_{r})}^2 \ra,
		\end{aligned}
	\end{equation}
	for some small $ \delta>0 $. 
	This inequality will help us to bound the quasilinear precession (cross product) term and the quadratic gradient term in $ I_1 $. 
	
	Since $ |u(t,x)| = 1 $ in $ P_{2r} $ as a result of its $ H^2 $-regularity, 
	we have
	\begin{align*}
		&-\lb \nabla (\phi^2 \Delta \u), \nabla (\u \times \Delta \u) \rb_{\b L^2} \\
		&= -\lb \phi^2 \nabla \Delta \u, \nabla \u \times \Delta \u \rb_{\b L^2} 
		-\lb \phi (\nabla \phi) \Delta \u, \u \times \nabla \Delta \u \rb_{\b L^2} \\
		&\lesssim 
		|\phi \nabla \Delta \u|_{\b L^2} |\phi \Delta \u|_{\b L^4} |\nabla \u|_{L^4(B_{r})} 
		+ r^{-1} |\phi \nabla \Delta \u|_{\b L^2} |\Delta \u|_{L^2(B_{r})},
	\end{align*}
	and similarly, with $ |\nabla \u(t,x)|^2 = -(\u \cdot \Delta \u)(t,x) $, 
	\begin{align*}
		&\lb \nabla (\phi^2 \Delta \u), \nabla(|\nabla \u|^2 u) \rb_{\b L^2} \\
		&= \lb \phi^2 \nabla \Delta \u + 2 (\phi \nabla \phi) \Delta \u, (-\nabla \u \cdot \Delta \u - \u \cdot \nabla \Delta \u) \u + |\nabla \u|^2 \nabla \u\rb_{\b L^2} \\
		&\lesssim 
		|\phi \nabla \Delta \u|_{\b L^2} |\phi \Delta \u|_{\b L^4} |\nabla \u|_{L^4(B_{r})} 
		+ r^{-1} |\phi \Delta \u|_{\b L^4} |\Delta \u|_{L^2(B_{r})} |\nabla \u|_{L^4(B_{r})} \\
		&\quad + r^{-1} |\phi \nabla \Delta \u|_{\b L^2} |\Delta \u|_{L^2(B_{r})}, 
	\end{align*}
	where the right-hand sides can be addressed by \eqref{eq: D321 + D221}. 
	Hence, 
	\begin{align*}
		|I_1|
		&\leq 
		\delta |\phi \nabla \Delta \u|_{\b L^2}^2 
		+ c |\nabla \u|_{L^4(B_{r})}^4 |\phi \Delta \u|_{\b L^2}^2
		+ c r^{-2} |\Delta \u|_{L^2(B_{r})}^2.
	\end{align*}

	\underline{$ I_2 $ estimate}	
	\begin{align*}
		\lb \nabla(\phi^2 \Delta \u), \nabla F \rb_{\b L^2}
		&\lesssim \la |\phi \nabla \Delta \u |_{\b L^2} + r^{-1} |\Delta \u|_{L^2(B_{r})} \ra
		|\phi \nabla F|_{\b L^2}.
	\end{align*}
	In view of \eqref{eq: |DF|},
	\begin{align*}
		|\phi \nabla F|_{\b L^2} 
		&\lesssim  
		|\phi^\frac{1}{2} A|_{\b L^4} |\phi^\frac{1}{2} \nabla A|_{\b L^4} 
		+ |\phi \nabla (\nabla \cdot A)|_{\b L^2} \\
		&\quad + \la |\phi^\frac{1}{2} \nabla A|_{\b L^4} + |\phi^\frac{1}{4} A|_{\b L^8}^2 \ra |\phi^\frac{1}{2}\nabla \u|_{\b L^4} \\
		&\quad + |A|_{L^4(B_{r})} |\phi \Delta \u|_{\b L^4}.
	\end{align*}
	For the $ \b L^8 $-term, we apply \eqref{eq: lady} to $ f=|A|^2 $ and obtain
	\begin{align*}
		|\phi^\frac{1}{4} A|_{\b L^8}^4
		&\lesssim |A|_{L^4(B_{r})}^2 \la |\phi^\frac{1}{2}A|_{\b L^4} |\phi^\frac{1}{2} \nabla A|_{\b L^4} + r^{-1} |A|_{L^4(B_{r})}^2 \ra,
	\end{align*}
	which implies 
	\begin{align*}
		|\phi^\frac{1}{4} A|_{\b L^8}^4 |\phi^\frac{1}{2}\nabla \u|_{\b L^4}^2
		&\lesssim 
		|A|_{L^4(B_{r})}^4 \la r^{-2} + |\phi^\frac{1}{2}\nabla \u|_{\b L^4}^4 \ra
		+ |\phi^\frac{1}{2}A|^2_{\b L^4} |\phi^\frac{1}{2} \nabla A|^2_{\b L^4}. 
	\end{align*}
	Thus, 
	\begin{align*}
		|I_2| 
		&\leq 
		\delta |\phi \nabla \Delta \u|_{\b L^2}^2 
		+ c\la 1+ |A|_{L^4(B_{r})}^4 \ra |\phi \Delta \u|_{\b L^2}^2 \\
		&\quad + c r^{-2} |\Delta \u|_{L^2(B_{r})}^2 + c |\phi \nabla (\nabla \cdot A)|_{\b L^2}^2 \\ 
		&\quad + c \la |\phi^\frac{1}{2} A|_{\b L^4}^2 + |\phi^\frac{1}{2} \nabla \u|_{\b L^4}^2 \ra |\phi^\frac{1}{2} \nabla A|^2_{\b L^4} \\
		&\quad + c |A|_{L^4(B_{r})}^4 \la r^{-2} + |\phi^\frac{1}{2} \nabla \u|_{\b L^4}^4 \ra.
	\end{align*}

	\underline{Under Hypothesis \ref{Hypo: small |Du|+|A|}}
	
	Under the smallness condition, we can estimate the time integral of $ |I_1|+|I_2| $ using Lemma \ref{Lemma: A L4 + H1}, \eqref{eq: Du L4} and \eqref{eq: D2u L2 Pr}.  
	In particular, for the last two terms in the estimate of $ |I_2| $, we have
	\begin{align*}
		&\int_{-r^2}^0 \la |\phi^\frac{1}{2} A|^2_{\b L^4} + |\phi^\frac{1}{2} \nabla \u|_{\b L^4}^2 \ra |\phi^\frac{1}{2} \nabla A|^2_{\b L^4} \ ds \\
		&\lesssim 
		\int_{-r^2}^0 |\phi \Delta \u|_{\b L^2}^2 \ dt + |\nabla \cdot A|_{L^2(P_{r})}^2 + \eps^2 \int_{-r^2}^0 \Psi\la \frac{r}{2},\eps,A \ra ds 
		+ \eps^2 + r^{-4} \eps^6, 
	\end{align*}
	and
	\begin{align*}
		&\int_{-r^2}^0 |A|_{L^4(B_{r})}^4 \la r^{-2}+ |\phi^\frac{1}{2} \nabla \u|_{\b L^4}^4 \ra ds \\
		&\lesssim 	
		\eps^2 \int_{-r^2}^0 |A|_{L^4(B_{r})}^4 |\phi \Delta u|_{\b L^2}^2 \ ds 
		+ r^{-2} \eps^2 |\nabla \cdot A|_{L^2(P_{2r})}^2 + r^{-2} \eps^4.
	\end{align*}
	Then it holds by \eqref{eq: |D2u(0)|-L2} that for a sufficiently small $ \delta>0 $ and $ t \in (-r^2,0) $,  
	\begin{align*} 
		|\phi \Delta u(t)|_{\b L^2}^2 
		+ \frac{\alpha}{2} \int_{-r^2}^0 |\phi \nabla \Delta \u|_{\b L^2}^2 \ ds 
		&\lesssim 
		v_0(r) 
		+ \int_{-r^2}^0 v_1(r) |\phi \Delta \u|_{\b L^2}^2 \ ds,
	\end{align*}
	where 
	\begin{align*}
		v_0(r) 
		&\lesssim r^{-2} |\nabla \cdot A|_{L^2(P_{2r})}^2 
		+ \eps^2 \int_{-r^2}^0 \Psi\la \frac{r}{2},\eps,A \ra ds \\
		&\quad + \eps^2 \la 1 + r^{-4} \eps^4 + r^{-2} \eps^2 \ra, \\
		v_1(r)
		&\lesssim \int_{-r^2}^0 \la 1+|\nabla \u|_{L^4(B_r)}^4 + |A|_{L^4(B_r)}^4 \ra ds \\
		&\lesssim 
		\eps^2 \la |\nabla \cdot A|_{L^2(P_{2r})}^2 + \eps^2 \ra.
	\end{align*}
	The estimate follows immediately from Gronwall's inequality.

	{\it Part (ii).}
	Take $ \rho = \bar{r} := \frac{r}{2} $. 
	Similarly, we omit the subscript $ \bar{r} $ for $ \phi $, $ \eta $ and $ \varphi $ below. 
	\begin{align*}
		\frac{1}{2}\partial_t |\phi \nabla \Delta \u(t)|_{\b L^2}^2 
		&= -\lb \nabla \la \phi^2 \nabla \Delta \u \ra, \Delta \partial_t \u \rb_{\b L^2} 
		+ \lb \phi \nabla \Delta \u, (\partial_t \phi) \nabla \Delta \u \rb_{\b L^2} \\ 
		&= 
		-\alpha|\phi \Delta^2 \u|_{\b L^2}^2 - \alpha J_0 - (J_1 + J_2) + J_3,
	\end{align*}
	where 
	\begin{align*}
		J_0 &= \lb (\nabla \phi^2) \nabla \Delta \u, \Delta^2 \u \rb_{\b L^2}, \\
		J_1 &= \lb \nabla \la \phi^2 \nabla \Delta \u \ra, \alpha \Delta (|\nabla \u|^2 \u) + \Delta(\u \times \Delta \u) \rb_{\b L^2}, \\
		J_2 &= \lb \nabla \la \phi^2 \nabla \Delta \u \ra, \Delta F(t,A,\u) \rb_{\b L^2}, \\
		J_3 &= \lb \phi \nabla \Delta \u, (\partial_t \phi) \nabla \Delta \u \rb_{\b L^2}. 
	\end{align*}
	Similarly, we have
	\begin{align*}
		\alpha |J_0| + |J_3|
		&\leq \delta |\phi \Delta^2 \u|_{\b L^2}^2 + c(\delta) r^{-2} |\nabla \Delta \u|_{L^2(B_{\bar{r}})}^2,
	\end{align*}
	where $ |\nabla \Delta \u|_{L^2(P_{\bar{r}})}^2 $ is bounded by part (i). 
	Then choosing a sufficiently small $ \delta > 0 $, 
	\begin{equation}\label{eq: |D3u(0)|-L2}
		\frac{1}{2} |\phi \nabla\Delta \u(t)|_{\b L^2}^2 
		+ \frac{\alpha}{2} |\phi \Delta^2 \u|_{\b L^2}^2 
		\leq 
		cr^{-2} |\nabla \Delta \u|_{L^2(P_{\bar{r}})}^2 
		- \int_{-\bar{r}^2}^t ( J_1 + J_2 ) \ ds,
	\end{equation}
	for any $ t \in (-\bar{r}^2,0) $. 
	We estimate the nonlinear terms $ J_1 $ and $ J_2 $ below. 
	
	\underline{$ J_1 $ estimate}
	
	As in $ I_1 $, we first observe that
	\begin{align*}
		|\phi \nabla \Delta \u|_{\b L^4}
		&\lesssim 
		|\phi \Delta^2 \u|_{\b L^2}^\frac{1}{2} |\phi \nabla \Delta \u|_{\b L^2}^\frac{1}{2} + r^{-\frac{1}{2}} |\nabla \Delta \u|_{L^2(B_{\bar{r}})}^\frac{1}{2} |\phi \nabla \Delta \u|_{\b L^2}^\frac{1}{2}. 	
	\end{align*}
	Then with $ |\u(t,x)| = 1 $ in $ P_{2r} $, for the two terms of $ J_1 $ we have
	\begin{align*}
		&\lb \nabla (\phi^2 \nabla \Delta \u), \Delta(\u \times \Delta \u) \rb_{\b L^2} \\
		&= 2\lb \phi^2 \Delta^2 \u, \nabla \u \times \nabla \Delta \u \rb_{\b L^2} 
		+ 2\lb \phi(\nabla \phi) \nabla \Delta \u, \u \times \Delta^2 \u \rb_{\b L^2} \\
		&\lesssim 
		|\phi \Delta^2 \u|_{\b L^2} |\phi \nabla \Delta \u|_{\b L^4} |\nabla \u|_{L^4(B_{\bar{r}})} 
		+ r^{-1} |\phi \Delta^2 \u|_{\b L^2} |\nabla \Delta \u|_{L^2(B_{\bar{r}})}, 
	\end{align*}
	and 
	\begin{align*}
		&\lb \nabla(\phi^2 \nabla \Delta \u), \Delta (|\nabla \u|^2 \u) \rb_{\b L^2} \\
		&= \lb \nabla(\phi^2 \nabla \Delta \u), 2 (\nabla \u \cdot \nabla \Delta \u + |\nabla^2 \u|^2) \u + |\nabla \u|^2 \Delta \u + 4 (\nabla \u \cdot \nabla^2 \u) \nabla \u \rb_{\b L^2} \\
		&\lesssim 
		|\phi \Delta^2 \u|_{\b L^2} 
		\la |\phi \nabla \Delta \u|_{\b L^4} |\nabla \u|_{L^4(B_{\bar{r}})}
		+ |\phi^\frac{1}{2} \Delta \u|_{\b L^4}^2 \ra \\
		&\quad + r^{-1} \la |\phi \nabla \Delta \u|_{\b L^4} |\nabla \Delta \u|_{L^2(B_{\bar{r}})} |\nabla \u|_{L^4(B_{\bar{r}})} 
		+ |\phi \nabla \Delta \u|_{\b L^2} |\Delta \u|_{L^4(B_{\bar{r}})}^2 \ra.
	\end{align*}
	Thus, using the $ L^4 $-estimate of $ \phi \nabla \Delta \u $, we obtain
	\begin{align*}
		|J_1| 
		&\leq
		\delta |\phi \Delta^2 \u|_{\b L^2}^2 
		+ c(\delta)|\phi \nabla \Delta \u|_{\b L^2}^2 \la 1 + |\nabla \u|_{L^4(B_{\bar{r}})}^4 \ra \\
		&\quad + c(\delta) \la r^{-2} |\nabla \Delta \u|_{L^2(B_{\bar{r}})}^2 + |\Delta \u|_{L^4(B_{\bar{r}})}^4 \ra,
	\end{align*}
	where
	\begin{align*}
		|\Delta \u|_{L^4(P_{\bar{r}})}^4
		&\lesssim \sup_{t \in [-\bar{r}^2,0]}|\Delta \u(t)|_{L^2(B_{\bar{r}})}^2 \la |\Delta \u|_{L^2(P_{\bar{r}})}^2 + |\nabla \Delta \u|_{L^2(P_{\bar{r}})}^2 \ra.
	\end{align*}

	\underline{$ J_2 $ estimate}
	
	For the gauge term, 
	\begin{align*}
		|\phi \Delta F|_{\b L^2} 
		&\lesssim 
		|\Delta(\nabla \cdot A)|_{L^2(B_{\bar{r}})} + |\nabla A|_{H^1(B_{\bar{r}})}^2 \\
		&\quad + \la |\phi \Delta A|_{\b L^2} + |\phi \nabla(\nabla \cdot A)|_{\b L^2} + |\phi^\frac{1}{2} A|_{\b L^4} |\phi^\frac{1}{2}\nabla A|_{\b L^4} \ra |\nabla \u|_{L^\infty(B_{\bar{r}})} \\
		&\quad + \la |\nabla A|_{L^4(B_{\bar{r}})} + |A|_{L^8(B_{\bar{r}})}^2 + |\nabla \u|_{L^\infty(B_{\bar{r}})} \ra |\phi \Delta \u|_{\b L^4} \\
		&\quad + |A|_{L^4(B_{\bar{r}})} |\phi \nabla \Delta \u|_{\b L^4} 
	\end{align*}
	Applying \eqref{eq: phiD2-L4}, the previous $ \b L^4 $-estimate of $ \phi \nabla \Delta \u $, the $ \b L^8 $-estimate of $ A $ in the proof of part (i), and the inequality
	\begin{align*}
		|\nabla \u|_{L^\infty(B_{\bar{r}})}^2 &\lesssim |\nabla \Delta \u|_{L^2(B_{\bar{r}})} |\nabla \u|_{L^2(B_{\bar{r}})} + |\nabla \u|_{L^2(B_{\bar{r}})}^2,
	\end{align*}
	we obtain
	\begin{align*}
		|J_2|-\delta |\phi \Delta^2 \u|_{\b L^2}^2
		&\lesssim 
		\la 1 + |A|_{L^4(B_{\bar{r}})}^4 + |\Delta \u|_{L^2(B_{\bar{r}})}^2 \ra |\phi \nabla \Delta \u|_{\b L^2}^2 \\ 
		&\quad + \la r^{-2} + |\nabla \u|_{L^2(B_{\bar{r}})}^2 \ra |\nabla \Delta \u|_{L^2(B_{\bar{r}})}^2 \\
		&\quad + \la 1 + |\nabla \u|_{L^2(B_{\bar{r}})}^2 \ra |\Delta \u|_{L^2(B_{\bar{r}})}^4 \\
		&\quad + |\Delta(\nabla \cdot A)|_{L^2(B_{\bar{r}})}^2 + |\nabla A|_{H^1(B_{\bar{r}})}^2 + |\nabla A|_{L^4(B_{\bar{r}})}^4 \\
		&\quad + |\nabla u|_{L^2(B_{\bar{r}})}^2 |\phi \Delta A|_{\b L^2}^4 + |\phi \nabla(\nabla \cdot A)|_{\b L^2}^4 \\
		&\quad + (1+ |A|_{L^4(B_{\bar{r}})}^4) \la 1 + |\nabla A|_{L^4(B_{\bar{r}})}^4 \ra
		+ r^{-2} |A|_{L^4(B_{\bar{r}})}^8. 
	\end{align*}
	
	\underline{Under Hypothesis \ref{Hypo: small |Du|+|A|}}
	
	By the smallness assumption and Lemma \ref{Lemma: A L4 + H1}, 
	\begin{align*}
		|J_2|-\delta |\phi \Delta^2 \u|_{\b L^2}^2
		&\lesssim 
		\la 1 + |A|_{L^4(B_{\bar{r}})}^4 + |\Delta \u|_{L^2(B_{\bar{r}})}^2 \ra |\phi \nabla \Delta \u|_{\b L^2}^2 \\ 
		&\quad + r^{-2} |\nabla \Delta \u|_{L^2(B_{\bar{r}})}^2 + |\Delta \u|_{L^2(B_{\bar{r}})}^4 \\
		&\quad + |\Delta(\nabla \cdot A)|_{L^2(B_{\bar{r}})}^2 
		+ |\nabla \cdot A|_{H^1(B_{2\bar{r}})}^4  \\
		&\quad + \Psi(\bar{r}) + \eps^2 \Psi^2\la \frac{\bar{r}}{2} \ra + r^{-4} \eps^2 \la 1 + r^{-4} \eps^4 \ra + r^{-6} \eps^4 + 1. 
	\end{align*}
	Similarly, we choose a sufficiently small $ \delta > 0 $ and substitute the estimates into \eqref{eq: |D3u(0)|-L2}. Then the result follows from estimates in part (i).  
\end{proof}

\begin{proposition}\label{Prop: u H2 est global}
	Assume Hypothesis \ref{Hypo: small |Du|+|A|} holds with $ I = [0,t] $ and $ \eps $ as in Proposition \ref{Prop: Du L4 + H1}. 
	Then for $ \u = \u(\omega) $ and $ A = A(\omega) $, 
	\begin{enumerate}[leftmargin=*,parsep=2pt]
		\item[(i)]
		$ \Delta \u \in L^2(I; \b L^2) $, 
		with bound that only depends on $ r $, $ \eps $, $ |A|_{L^\infty(I; \b H^1)} $ and $ |\nabla \u_0|_{\b L^2} $,
		
		\item[(ii)]
		$ u \in L^\infty(\rho, t; \b H^3) \cap L^2(\rho, t; \b H^4) $ for arbitrary $ \rho \in (0,t) $, 
		with bound that only depends on $ r $, $ \eps $ and $ |A|_{L^\infty(I; \b H^2) \cap L^2(I; \b H^3)} $, 
		
		\item[(iii)]
		$ u \in L^\infty(I; \b H^3) \cap L^2(I; \b H^4) $ provided that $ \u_0 \in \b H^3 $, 
		with bound that only depends on $ r $, $ \eps $, $ |A|_{L^\infty(I; \b H^2) \cap L^2(I; \b H^3)} $ and $ |u_0|_{\b H^3} $.
	\end{enumerate}
\end{proposition}
\begin{proof}
	Using a covering argument for $ \T^2 $ by balls, we easily obtain from \eqref{eq: Du L4} that 
	\begin{align*}
		|\nabla \u|_{\b L^4}^4 \lesssim r^{-2} \eps^2 \la |\nabla^2 \u|_{\b L^2}^2 + r^{-2} |\nabla \u|_{\b L^2}^2 \ra. 
	\end{align*}
	Taking $ \varphi \equiv 1 $ with $ \nabla \varphi \equiv 0 $ in the proof of Proposition \ref{Prop: Du L4 + H1} and using the $ L^4 $-estimate above, we obtain a sharpened bound of $ \nabla u $ in $ L^4(I; \b L^4) \cap L^2(I; \b H^1) $ depending only on $ r $, $ \eps $, $ |A|_{\b H^1} $ and $ |\nabla \u_0|_{\b L^2} $, proving part (i). 
	This can then be used to deduce the high-order estimates in part (ii) for positive times using Proposition \ref{Prop: u H2 est Pr}, and part (iii) follows similarly by taking the cut-off functions $ \varphi_\rho = \eta_\rho \equiv 1 $ in the proof of Proposition \ref{Prop: u H2 est Pr}. 
\end{proof}

\section{Well-posedness in energy space}\label{Section: well-posedness}

In this section, we state the existence and uniqueness result for a weak solution of the transformed equation \eqref{eq: u} taking values in $ C([0,T]; \b L^2) \cap L^\infty(0,T; \b H^1) $ a.s.  
For clarity, we outline the main steps of the proofs below before stating the results.  

For the proof of existence, first, for a.e. $ \omega \in \Omega $, local well-posedness of \eqref{eq: u} in $\b H^1$ can be achieved by means of an approximation procedure. 
Recall that the approximability of $\b H^1$ initial data by smooth fields is guaranteed by the density result of Schoen and Uhlenbeck \cite{SU_density}. 
Then given smooth initial data, the existence and uniqueness of maximal solutions in $ \b H^2 $ with improved regularity (time integrability in $ \b H^4 $) under small energy assumptions follow from Section \ref{Section: higher-order estimates}. 
A compactness argument on a uniform interval is required for the convergence of these approximating solutions. 
If we proceed directly by using the bounds obtained in Proposition \ref{Prop: u H2 est global}, then we obtain a different convergent subsequence for each $ \omega $ (in a subset of $ \Omega $ of full measure), leaving measurability of the limiting function and the underlying time interval in doubt. 
In particular, the dependence of the limit on the stochastic process $ A $ is unclear. 
To circumvent this issue, we first apply the Skorohod theorem using the moment estimates in Lemma \ref{Lemma: energy moments} to establish the existence of an a.s.-convergent subsequence of the approximations and hence progressiveness of the limit. 
Then for each $ \omega $, a compactness argument can be applied to obtain regularity of the limit in $ \b H^1 $ starting from time $ 0 $, and beyond $ \b H^1 $ on a positive time interval defined by an energy-related stopping time. 
This allows us to extend the limit to a progressively measurable maximal $ \b H^2 $-solution defined at positive times, and analyse weak $ \b H^1 $ convergence and the blow-up scenario at the maximal time. 
Finally, iteration of the above steps leads us to a global $ \b H^1 $ solution.

\begin{theorem}\label{Theorem: Struwe u}   
	Assume that \eqref{eq_q3} holds for $ \sigma = 4 $. 
	For every $\u_0 \in H^1(\T^2; \St)$ there exist a stopping time $\tau^*>0$
	and a progressively measurable process $ \u $ that takes values in the space 
	\[
	H^1(0,\tau^*; \b L^2) \cap L^\infty\la 0,\tau^*; H^1(\T^2; {\b S}^2)\ra \cap L^2_{\rm loc}\la [0,\tau^*); \b H^2 \ra 
	\quad \text{$ \P $-a.s.} 
	\]
	and satisfies \eqref{eq: u} with $ \u(0) = \u_0 $, $ \P $-a.s. 
	Moreover, there exist random variables 
	$N < \infty$ ($ \in \N $), $x_1, \dots, x_{N} \in \T^2$ and $\u^* $ in the space $H^1(\T^2; \St)$, $ \P $-a.s. such that 
	\[
	\lim_{i \to \infty} \u(t_i) = \u^* \quad \text{weakly in} \quad H^1(\T^2; \R^3),
	\]
	and strongly in $ H^1(\T^2 \setminus \{x_1,\dots, x_{N}\}) $, for any sequence of stopping times $ (t_i)_i $ with $ t_i < \tau^* $ for all $ i \in \b N $ and $ t_i \nearrow \tau^* $, $ \P $-a.s. 
	Moreover, for every $r>0$ and a.e.-$ \omega \in \Omega $,
	\[
	\limsup_{t \nearrow \tau^*(\omega)} \frac{1}{2} \int_{B_r(x_n(\omega))} |\nabla \u(t)|^2(\omega) \, dx \geq 4 \pi,
	\quad \text{$n=1, \dots, N(\omega) $}.
	\]
\end{theorem}

Taking $ \u^* $ as new random initial data, the construction in Theorem \ref{Theorem: Struwe u} can be iterated. In contrast to the unperturbed case, the noise may bring energy influx into the system. 
We rely on the energy bound 
\begin{align*}
	\sup_{\ell \in \b N} \sup_{t \in [\tau_{\ell-1},\tau_{\ell}]} |\nabla \u(t)|_{\b L^2}^2 + 8 \pi \sum_{j=1}^\ell N_j < c(|\u_0|_{\b H^1}, |A|_{L^2(0,T; \b H^2)}), \quad \P\text{-a.s.}
\end{align*}
to show that the blow-up times (when truncated at $ T $) yield $ \lim_{\ell \to \infty} \tau_{\ell} = T $. 

\begin{corollary}\label{Corollary: Struwe u iteration}
	For every $ T>0 $, there exists a progressively measurable process 
	\begin{align*}
		\u \in H^1(0,T; \b L^2) \cap L^\infty(0,T; H^1(\T^2; \St)), \quad \P\text{-a.s.}
	\end{align*}
	and an increasing sequence of stopping times 
	$ (\tau_\ell)_{\ell \in \N} $
	such that 
	$ T \in (\tau_M, \tau_{M+1}) $ for some $ M \geq 0 $ with $ \tau_0 = 0 $, and that for all $ \ell \leq M $,
	\begin{enumerate}[leftmargin=*,parsep=2pt]
		\item[(a)]  
		$ \u|_{[\tau_\ell, \tau_{\ell+1} \wedge T)} $ is a local strong solution of \eqref{eq: u} as in Theorem \ref{Theorem: Struwe u} with $ \u^* = \u(\tau_{\ell+1}) $ and the number of singular points $ N_{\ell+1} $ at $ \tau_{\ell+1} $, 
		
		\item[(b)] 
		$ \sup_{\ell \in \b N} \sum_{j=1}^\ell N_j < \infty $, $ \P $-a.s. 
	\end{enumerate}
\end{corollary}

\begin{remark}
	The effect of $A$ becomes negligible at small scales. 
	For $ z_0=(t_0,x_0) $ and $ \lambda \in \R $, solutions $\u$ and connection forms $A$ in $ P_{\lambda r}(z_0) $ (if exist) obey the following scaling laws on $ P_r $, 
	\[
	\u_\lambda(t,x) = \u(t_0+ \lambda^2 t, x_0+\lambda x), \quad A_\lambda(t,x)= \lambda A(t_0+ \lambda^2 t, x_0+\lambda x).
	\]
	Note that $\nabla \u$ and $A$ share the same scaling behaviour - invariant in the $L^2(\R^2)$ norm. 
	Moreover, the corresponding perturbation has a higher order scaling behaviour
	$ 
	F_\lambda(t,x)=\lambda^2 F(t_0+ \lambda^2 t, x_0+\lambda x)
	$ 
	and becomes negligible on small scales $\lambda \ll 1$. 
	Thus, the existence and partial regularity result for the free harmonic Landau-Lifshitz equation is essentially reproduced apart from the fact that singularities are no longer controlled by the initial energy.
\end{remark}

For the uniqueness of solution, it is well-known that weak solutions $ \u:(0,\infty) \times \T^2 \to \St$ to the LLG in the energy space 
\begin{align*}
	\c E_T := \{ \v: [0,T] \times \T^2 \to \b S^2 : \partial_t \u \in L^2(0,T; \b L^2), \nabla \u \in L^\infty(0,T; \b L^2)\}
\end{align*}
for arbitrary finite $ T>0 $, are generally not unique, a phenomenon that is attributed to bubbles that form backwards in time. 
Freire's arguments for the uniqueness of energy-decreasing solutions to the harmonic map heat flow \cite{Freire}, extended to LLG in \cite{ChenGuo, Harpes}, can be adapted to the magnetic case \eqref{eq: u}, provided that the perturbation $F$ enjoys sufficient integrability. 
Note that if $\sigma \geq 2$, 
then $ F \in L^\infty(0,T; \b L^p)$ for all $p<2$ which is more than sufficient. 
In this case, energy influx is possible, thus the uniqueness condition is to be modified. 
Namely, we show that the solution $\u$ is unique in the class of weak solutions with right-continuous Dirichlet energy. 
\begin{theorem}\label{Theorem: u unique}
	Assume that \eqref{eq_q3} holds for some $ \sigma \geq 2 $. Fix $ \omega \in \Omega $, uniqueness of solution $ \u(\omega) $ to \eqref{eq: u} holds in the class 
	$$ \{ \v \in \c E_T: t \mapsto |\nabla_{A(\omega)} \v(t)|^2_{\b L^2} \text{ is right-continuous} \}. $$ 
\end{theorem}
Note that $ \c E_T \subset C([0,T]; \b L^2) $ so that the right-continuity condition can equivalently be phrased in terms of $\nabla \v$ rather than $\nabla_A \v$. 
Due to weak lower semi-continuity, the right continuity above is equivalent to
\[
\limsup_{s\searrow t} |\nabla \v(s)|^2_{\b L^2} \le |\nabla \v(t)|^2_{\b L^2}, \quad \forall t<T. 
\]
Here, energy evaluation is meant in the sense of traces, see \cite{Topping}.

Applying the inverse transformation, 
we obtain the existence of weak martingale solutions to \eqref{eq: sLLG} from the equivalence result (Lemma \ref{Lemma: u vs m weak}). 
These solutions are locally strong due to the interpolation embedding $ H^1(0,t; \b L^2) \cap L^2(0,t; \b H^2) \hookrightarrow C([0,t]; \b H^1) $. Moreover, they are pathwise unique thanks to Theorem \ref{Theorem: u unique}, yielding our main result Theorem \ref{Theorem: Struwe m} and \ref{Theorem: unique m}.

\section{Proof of existence}\label{Section: Proof existence}

The proof of Theorem \ref{Theorem: Struwe u} is divided into several steps: Section \ref{Section: step 1} -- \ref{Section: step 3} for the construction of solution $ \u $, and Section \ref{Section: step 6} for an estimate of $ \eps $. 
The proof of Corollary \ref{Corollary: Struwe u iteration} is in Section \ref{Section: step 5}.

\subsection{Approximations in $ C([0,T];\b H^2) $}\label{Section: step 1}
By the density result in \cite{SU_density}, there exists a sequence \smash{$ \{\u_0^{(k)}\}_{k \in \b N} $} in $ C^\infty(\T^2; \St)$ such that \smash{$ \u_0^{(k)} \to \u_0 $} strongly in $ H^1(\T^2; \St) $ as $ k \to \infty $. 
Given $ \eps >0 $ as in Proposition \ref{Prop: u H2 est global}, it is possible to select $r_0 \in (0,1)$ such that 
\begin{equation}\label{eq: u0k eps}
	\sup_{k \in \N}  \sup_{x \in \T^2} \left|\nabla \u^{(k)}_0 \right|^2_{L^2(B_{4r_0}(x))} \leq \frac{\eps^2}{4}.
\end{equation} 
For every $ k \in \N $ and a.e.-$ \omega \in \Omega $, given initial data \smash{$ \u_0^{(k)} $}, by Corollary \ref{Coro: u S2} there exists a unique $ \b H^2 $-maximal solution \smash{$ (\u^{(k)},[0,\tau^{(k)}))(\omega)$} of equation \eqref{eq: u}. 

Next, we show that the solutions $ \u^{(k)} $ and the parameter $ A $ satisfy Hypothesis \ref{Hypo: small |Du|+|A|} on a uniform interval.  
Using Lemma \ref{Lemma: local energy} and \eqref{eq: u0k eps}, there exists a constant $ c>0 $ such that for a.e.-$ \omega \in \Omega $, 
\begin{equation}\label{eq: uk epsh}
	\sup_{x \in \T^2} |\nabla \u^{(k)}(t,\omega)|_{L^2(B_{2r_0}(x))}^2 
	\leq 
	\frac{\eps^2}{4} + c h(t,\omega) t, \quad t \in [0,\tau^{(k)}(\omega)),
\end{equation} 
where $ h:[0,T] \times \Omega \to (0,\infty) $ is given by 
\begin{equation}\label{eq: h(t,u0,r0)}
	\begin{aligned}
		h(t)
		= h(t; A,\u_0) 
		&:= \la 1+\sup_{k \in \b N} |\nabla \u_0^{(k)}|_{\b L^2}^2 + \int_0^t |A|^4_{\b H^1} ds \ra \\
		&\ \ \quad \times e^{\int_0^t |A|_{\b H^2}^2 ds} \la r_0^{-2} + |A|^2_{C([0,t]; \b H^2)} \ra \\
		&\ \ \quad + 1+|A|_{C([0,t];\b H^1)}^4.
	\end{aligned}
\end{equation}
The process $ h $ is \smash{$ {\b F} $}-adapted and independent of $ k $. 
For a.e.-$ \omega \in \Omega $, $ h(\cdot,\omega) $ is increasing and continuous, since $ A $ is progressively measurable with values in $ C([0,T];\b H^2) $, $ \P $-a.s. 
Then we can define an $ \b F $-stopping time 
\begin{equation}\label{eq: defn tauh}
	\tau_h := \inf \left\{ t \in [0,T] : h(t)t = \frac{\eps^2}{4c} \right\},
\end{equation}
where $ \tau_h(\omega) = T $ if \smash{$ \sup_{t \in [0,T]} h(t,\omega)t < \frac{\eps^2}{4c} $}. 
By definition, we have either 
$ \tau_h(\omega) = T $ or 
$ h(T) \tau_h > h(\tau_h) \tau_h = \smash{\frac{\eps^2}{4c}} $ in the case $ \tau_h(\omega) < T $, implying that 
$ \tau_h > 0 $, $ \P $-a.s.  	
Coming back to \eqref{eq: uk epsh}, let $ I^{(k)} = [0,\tau^{(k)}) \cap [0,\tau_h] $ and $ r \in (0,r_0] $, then
\begin{equation}\label{eq: Duk eps}
	\sup_{k \in \b N} \sup_{t \in I^{(k)}}\sup_{x \in \T^2} |\nabla \u^{(k)}(t)|_{L^2(B_{2r}(x))}^2 
	\leq \frac{\eps^2}{4} + c h(\tau_h) \tau_h 
	\leq \frac{\eps^2}{2}, \quad \P\text{-a.s.}
\end{equation}	
Let $ \c N := \{ \omega \in \Omega: \tau_h(\omega) \geq \tau^{(k)}(\omega) \} $. 
If $ \P(\c N) >0 $, then $ \u^{(k)}(\omega) \in L^\infty(0,\tau^{(k)}; \b H^3) \cap L^2(0,\tau^{(k)};\b H^4) $ for any $ \omega \in \c N $ by Proposition \ref{Prop: u H2 est global}(iii) and \eqref{eq: Duk+A eps}. 
This implies that the $ \b H^2 $-solution \smash{$ (\u^{(k)},\tau^{(k)})(\omega) $} can be continuously extended up to (including) $ \tau^{(k)} $, which contradicts with the maximality of $ \tau^{(k)} $. 
Thus, $ \P(\c N) = 0 $, meaning $ \tau_h < \tau^{(k)} $, $ \P $-a.s. for all $ k \in \b N $.
By the regularity of $ A $, there exists $ r_{A} > 0 $, $ \P $-a.s. such that for any $ r \in (0, r_{A}] $, 
\begin{equation}\label{eq: A eps}
	\sup_{t \in [0,T]} \sup_{x \in \T^2} |A(t)|^2_{L^2(B_{2r}(x))} 
	\lesssim r^2 |A|^2_{C([0,T]; \b H^2)}
	< \frac{\eps^2}{2}, \quad \P\text{-a.s.}
\end{equation}
where the first inequality holds by the Sobolev embedding $ \b H^2 \hookrightarrow \b L^\infty $ and the fact $ |B_{2r}| \lesssim r^2 $. 
For example, we can take $ r_A^2 = \frac{\eps^2}{4} (1+|A|^{2}_{C([0,T]; \b H^2)})^{-1} $.

Define $ \u^{(k)}_h(t) := \u^{(k)}(t \wedge \tau_h) $ for $ t \in [0,T] $. 
Then by \eqref{eq: Duk eps} -- \eqref{eq: A eps}, 
\begin{equation}\label{eq: Duk+A eps}
	\sup_{t \in [0,T]}\sup_{x \in \T^2} \la |\nabla \u_h^{(k)}(t)|_{L^2(B_{2r_0}(x))}^2 + |A(t)|^2_{L^2(B_{2r_A}(x))} \ra
	< \eps^2, \quad \P\text{-a.s.} 
\end{equation} 
We will use $ \b H^1 $-moment estimates of \smash{$ \u^{(k)}_h $} and Skorohod theorem in Section \ref{Section: step 2} to secure adaptedness of the limiting process, and leave higher-order regularity proof to Section \ref{Section: step 3} which require (a similar condition to) \eqref{eq: Duk+A eps}.

\subsection{Convergence to an adapted process $ \tilde{\u} $}\label{Section: step 2}
For $ \beta \in \R $, $ \sigma \geq 2 $ and $ \gamma \in (0,\frac{1}{2}) $, we define the spaces 
\begin{align*}
	U_0^\gamma &:= C^\gamma([0,T];\b L^2) \cap C([0,T]; \b H^1), \\
	V_0^{\beta,\sigma} &:= C([0,T]; \b H^{\beta}) \times C([0,T];\b H^\sigma) \times C([0,T]; \ell^2),
\end{align*}
where $ U_0^\gamma $ is compactly embedded in $ C([0,T]; \b H^{-1}) $ by Arzela-Ascoli's theorem. 
Recall that $ \tau_h \in (0,\tau^{(k)}) $, $ \P $-a.s. 
By Lemma \ref{Lemma: energy moments}, \smash{$ \u^{(k)}_h $} is uniformly bounded in $ L^2(\Omega;U_0^\gamma) $, 
Thus, the laws \smash{$ \{\mathcal{L}(\u^{(k)}_h)\}_{k \in \b N} $} are tight on $ C([0,T];\b H^{-1}) $. 

Applying Prokhorov theorem and a version of Skorohod representation theorem \cite[Theorem 3.2]{BertiPratelliRigo},
there exist \smash{$ V_0^{-1,\sigma} $}-valued random variables, 
\smash{$ (\tilde{\u}, \tilde{\Y}, \tilde{W}) $} and \smash{$ \{(\tilde{\u}^{(k)}, \tilde{\Y}^{(k)},  \tilde{W}^{(k)})\}_{k} $}, 
defined on the same probability space $ (\Omega, \mc F, \P) $,   
such that
\begin{enumerate}[(i)]
	\ie
	\item 
	$ \mathcal{L}((\tilde{\u}^{(k)}, \tilde{\Y}^{(k)}, \tilde{W}^{(k)})) = \mathcal{L}((\u^{(k)}_h, \Y, W)) $ on $ V_0^{-1,\sigma} $, for all $ k \in \b N $,
	\item 
	$ (\tilde{\u}^{(k)}, \tilde{\Y}^{(k)}, \tilde{W}^{(k)}) \to (\tilde{\u}, \tilde{\Y}, \tilde{W}) $ in $ V_0^{-1,\sigma} $, $ \P $-a.s.
\end{enumerate}
By (i), \smash{$ (\tilde{\u}^{(k)}, \tilde{\Y}^{(k)}, \tilde{W}^{(k)}) $} and \smash{$ (\u^{(k)}_h, \Y, W) $} have same laws on \smash{$ V_0^{-1,\sigma} $} and any separable metric subspace of \smash{$ V_0^{-1,\sigma} $} by Kuratowski's theorem. 

Let $ \tilde{\mc F}_t $ be the augmentation of \smash{$ \sigma(\tilde{\u}^{(k)}(s), \tilde{W}^{(k)}(s), \tilde{\u}(s), \tilde{W}(s) ; k \in \b N, s\leq t) $} and let $ \tilde{\b F} := \smash{(\tilde{\mc F}_t)_{t \in [0,T]}} $. 
Then for every $ k \in \b N $, $ \tilde{W}^{(k)} $ defines a cylindrical Wiener process on \smash{$ (\Omega, \mc F, \tilde{\b F}, \P) $}, and $ \tilde{\Y}^{(k)} $ is the unique solution of \eqref{eq: mc Y} driven by $ \tilde{W}^{(k)} $ satisfying Lemma \ref{Lemma: Y flow}. 	
Let
\begin{align*}
	\tilde{A}^{(k)} := (\tilde{\Y}^{(k)})^\star \la \nabla \tilde{\Y}^{(k)} \ra,
\end{align*}
and $ r_{\tilde{A}^{(k)}} $ be chosen as in \eqref{eq: A eps} . 
Then $ \tilde{A}^{(k)} \in C([0,T]; \b H^2) $, $ \P $-a.s. and the laws $ \mathcal{L}((\tilde{\u}^{(k)}, \tilde{A}^{(k)}, \tilde{W}^{(k)})) = \mathcal{L}((\u_h^{(k)}, A, W)) $ on \smash{$ V_0^{2,3} $} since $ C([0,T];\b H^2) $ is a separable subspace of \smash{$ C([0,T];\b H^{-1}) $}. 

Repeating the argument in Section \ref{Section: step 1}, 
we consider $ h $ in \eqref{eq: h(t,u0,r0)} as a function of $ t $, $ A $ and $ \u_0 $, and define
\begin{align*}
	\tilde{h}^{(k)}(t) &:= h(t; \tilde{A}^{(k)}, \u_0), \quad
	\tau_{\tilde{h}^{(k)}} := \inf \left\{ t \in [0,T] : \tilde{h}^{(k)}(t)t = \frac{\eps^2}{4c} \right\}. 
\end{align*}
Here, $ \tau_{\tilde{h}^{(k)}} $ is a $ \tilde{\b F} $-stopping time since $ \tilde{h}^{(k)} $ is a continuous and $ \tilde{\b F} $-adapted process by construction (see the definitions of $ V_0^{2,3} $ and $ \tilde{\b F} $). 
The law of $ \tau_{\tilde{h}^{(k)}} $ depends only on the law of $ \tilde{h}^{(k)} $. Then for a.e.-$ \omega \in \Omega $, the maps 
$ C([0,T]; \R) \ni \tilde{h}^{(k)}(\omega) \mapsto \tau_{\tilde{h}^{(k)}}(\omega) \in [0,T] $ 
and 
$ C([0,T]; \b H^2) \ni \tilde{A}^{(k)}(\omega) \mapsto h(\cdot; \tilde{A}^{(k)}(\omega), \u_0) = \tilde{h}^{(k)}(\omega) \in C([0,T]; \R) $ 
are Borel measurable. 
Then we have $ \mathcal{L}(\tilde{h}^{(k)}) = \mathcal{L}(h) $ on $ C([0,T]; \R) $ and $ \smash{\P(\tau_{\tilde{h}^{(k)}} > 0)} = \P(\tau_h > 0) = 1 $. 
Alternatively, the fact $ \tau_{\tilde{h}^{(k)}} > 0 $, $ \P $-a.s. is clear from the definition of $ \tilde{h}^{(k)} $. 
As a result of the same (joint) laws, for any $ t \in [0,T] $, we have on $ \b L^2 $,
\begin{align*}
	\tilde{\u}^{(k)}(t)- \u^{(k)}_0- \int_0^{t \wedge \tau_{\tilde{h}^{(k)}}} H(s,\tilde{A}^{(k)}, \tilde{\u}^{(k)}(s)) \ ds = 0, \quad \P\text{-a.s.}
\end{align*}
To find a uniform interval on which $ \tilde{\u}^{(k)} $, $ k \in \b N $, are defined, let
\begin{align*}
	\tilde{h}(t) &:= \sup_{k \in \b N}\tilde{h}^{(k)}(t), \quad
	\tau_{\tilde{h}} := \inf \left\{ t \in [0,T] : \tilde{h}(t)t > \frac{\eps^2}{4c} \right\}. 
\end{align*}
The process $ \tilde{h} $ is still progressively measurable with respect to $ \tilde{\b F} $. 
Recall from \eqref{eq: h(t,u0,r0)} that for every $ k \in \b N $, the path $ t \mapsto \tilde{h}^{(k)}(t,\omega)t $ is increasing and $ \tilde{h}^{(k)}(t,\omega) \geq 1 $, for a.e.-$ \omega \in \Omega $. 
Thus, $ \tilde{h}(\cdot, \omega) $ is increasing for a.e.-$ \omega \in \Omega $, and the hitting time $ \tau_{\tilde{h}} $ (of a measurable subset of $ \R $) is again a positive $ \tilde{\b F} $-stopping time, and $ \tau_{\tilde{h}} \leq \tau_{\tilde{h}^{(k)}} $ for all $ k \in \b N $.  
In other words, $ \tilde{\u}^{(k)} $ on $ [0,\tau_{\tilde{h}}] $ is a local solution of the equation \eqref{eq: u} with initial condition \smash{$ \u_0^{(k)} $} and parameter $ \tilde{A}^{(k)} $, $ \P $-a.s. 
Moreover, 
\begin{equation}\label{eq: utilde eps}
	\sup_{t \in [0,T]}\sup_{x \in \T^2} \la |\nabla \tilde{\u}^{(k)}(t \wedge \tau_{\tilde{h}})|_{L^2(B_{2r_0}(x))}^2 + |\tilde{A}^{(k)}(t)|^2_{L^2(B_{2r_{\tilde{A}^{(k)}}}(x))} \ra
	< \eps^2, \quad 
	\P\text{-a.s.}
\end{equation} 
We shall use \eqref{eq: utilde eps} for convergences in smaller spaces under fixed $ \omega $ in Section \ref{Section: step 3}. 
Before doing so, we collect some properties of the limit $ (\tilde{\u}, \tilde{\Y}, \tilde{W}) $ which will be used in Section \ref{Section: step 4}. 

As a result of the $ \P $-a.s. convergence, $ \tilde{W} $ is a cylindrical Wiener process and $ \tilde{\Y} $ is the unique solution of \eqref{eq: mc Y} driven by $ \tilde{W} $. As usual, we define $ \tilde{A} = \tilde{\Y}^\star \nabla \tilde{\Y} $, where $ \tilde{A}^{(k)} \to \tilde{A} $ in $ C([0,T]; \b H^{\sigma-1}) $, $ \P $-a.s. 
By construction, the processes $ \tilde{\u} $ and $ \tilde{\u}(\cdot \wedge \tau_{\tilde{h}}) $ are $ \tilde{\b F} $-progressively measurable. 
The following map is lower semi-continuous
\begin{align*}
	C([0,T]; \b H^{-1}) \ni f \mapsto \sup_{t \in [0,T]} |f(t)|_{\b H^1}.
\end{align*}
Thus, by \eqref{eq: utilde eps} and the $ \P $-a.s. convergence of $ \tilde{\u}^{(k)} $ in $ C([0,T]; \b H^{-1}) $ and $ \tilde{A}^{(k)} $ in $ C([0,T]; \b H^{\sigma -1}) $ (see (ii)), we have
\begin{equation}\label{eq: Du+A eps}
	\sup_{t \in [0,\tau_{\tilde{h}}]} \sup_{x \in \T^2} \la |\nabla \tilde{\u}(t)|_{L^2(B_{2r}(x))}^2 + |\tilde{A}(t)|^2_{L^2(B_{2r_{\tilde{A}}}(x))} \ra
	< \eps^2, 
	\quad \forall r \in (0,r_0], 
	\quad \P\text{-a.s.}
\end{equation}
where $ r_{\tilde{A}^{(k)}} \to r_{\tilde{A}} $, $ \P $-a.s. 
Similarly, using the uniform integrability of $ \tilde{\u}^{(k)} $ in $ C([0,T]; \b H^1) $ under (i) and Fatou's lemma, we have for any $ p \in [1,\infty) $,  
\begin{equation}\label{eq: tilde u H1 moment part 1}
	\begin{aligned}
		\E \left[ \sup_{t \in [0,T]} \left|\tilde{\u}(t \wedge \tau_{\tilde{h}}) \right|_{\b H^1}^p \right] 
		&\leq \liminf_{k \to \infty} \E \left[ \sup_{t \in [0,T]} \left| \tilde{\u}^{(k)}(t \wedge \tau_{\tilde{h}}) \right|_{\b H^1}^p \right] \\
		&= \liminf_{k \to \infty} \E \left[ \sup_{t \in [0,T]} \left|\u_h^{(k)}(t \wedge \tau_{\tilde{h}}) \right|_{\b H^1}^p \right] 
		< \infty. 
	\end{aligned}
\end{equation}


\subsection{$ \b H^2 $-maximal solution}\label{Section: step 3}
For $ \rho,\gamma \in (0,1) $, define the intervals
\begin{align*}
	I_{\tilde{h}} := [0,\tau_{\tilde{h}}], \quad 
	I_{\tilde{h},\rho}:= [\rho \tau_{\tilde{h}},\tau_{\tilde{h}}], 
\end{align*}
and let
\begin{align*}
	U_1 &:= L^\infty(I_{\tilde{h}}; \b H^1) \cap L^2(I_{\tilde{h}}; \b H^2) \cap H^1(I_{\tilde{h}}; \b L^2), \\
	V_1 &:= C(I_{\tilde{h}}; \b H^\gamma) \cap L^2(I_{\tilde{h}}; \b H^1), \\
	U_2 &:= L^\infty(I_{\tilde{h},\rho};\b H^3) \cap L^2(I_{\tilde{h},\rho};\b H^4) \cap H^1(I_{\tilde{h},\rho}; \b H^2), \\
	V_2 &:= C(I_{\tilde{h},\rho}; \b H^{2+\gamma}) \cap L^2(I_{\tilde{h},\rho}; \b H^3).
\end{align*}  
For fixed $ \omega \in \Omega $, we have the usual vector spaces and compact embeddings  
$ U_1(\omega) \Subset V_1(\omega) $ and $ U_2(\omega) \Subset V_2(\omega) $. 
By Proposition \ref{Prop: u H2 est global}, for a.e.-$ \omega \in \Omega $, 
\begin{align*}
	&\sup_{k \in \b N} \left|\tilde{\u}^{(k)}(\cdot \wedge \tau_{\tilde{h}})(\omega) \right|_{U_1(\omega) \cap U_2(\omega)} \\
	&< c\la r_0, \eps, \sup_{k \in \b N} |\u_0^{(k)}|_{\b H^1}, \sup_{k \in \b N} |\tilde{A}^{(k)}(\omega)|_{C([0,T]; \b H^3)}\ra 
	< \infty,
\end{align*}
where the fact $ \sup_{k \in \b N} |\tilde{A}^{(k)}|_{C([0,T]; \b H^3)} < \infty $, $ \P $-a.s. follows from the a.s.-convergence to $ \tilde{A} $ and the finite moments of $ \tilde{A} $. 
Then there exists a subsequence of $ \{\tilde{\u}^{(k)}(\cdot \wedge \tau_{\tilde{h}})(\omega)\} $, denoted as $ \{v^{(k_j)}\}_j $, such that 
\begin{equation}\label{eq: vkj to v}
	v^{(k_j)} \stackrel{j \to \infty}{\to} v 
	\quad \begin{cases} 
		\text{weak-* in $ L^\infty(I_{\tilde{h}}(\omega); \b H^1) \cap L^\infty(I_{\tilde{h},\rho}(\omega); \b H^3) $} \\
		\text{weakly in $ L^2(I_{\tilde{h}}(\omega); \b H^2) \cap L^2(I_{\tilde{h},\rho}(\omega);\b H^4) $} \\
		\text{strongly in $ V_1(\omega) \cap V_2(\omega) $.}
	\end{cases}
\end{equation}
Recall (ii) in Section \ref{Section: step 2}, we have 
\begin{align*}
	\tilde{\u}^{(k)}(\cdot \vee \rho \tau_{\tilde{h}} \wedge \tau_{\tilde{h}}) \to
	\tilde{\u}(\cdot \vee \rho \tau_{\tilde{h}} \wedge \tau_{\tilde{h}}) \quad 
	\text{ in } C([0,T]; \b H^{-1}), \quad \P\text{-a.s.}
\end{align*}  
(including the case $ \rho=0 $.)
Then by the uniqueness of limit in $ C(I_{\tilde{h}}(\omega); {\b H}^{-1}) $,   
\begin{align*}
	\tilde{\u}(t,\omega) = \v(t), \quad \forall t \in I_{\tilde{h}}(\omega),
\end{align*}
Using the convergences \eqref{eq: vkj to v} above,  $ \tilde{\u} $ is an $ \b H^1 $-solution of \eqref{eq: u} on $ I_{\tilde{h}} $ and an $ \b H^2 $-solution on $ I_{\rho,\tilde{h}} $, $ \P $-a.s., with regularity 
\begin{align}
	&\tilde{\u}(\cdot \wedge \tau_{\tilde{h}})
	\in C(0,T; \b H^\gamma) \cap L^\infty(0,T; \b H^1) \cap L^2(0,T; \b H^2), \quad \P\text{-a.s.} \label{eq: utilde reg on Ih} \\
	&\tilde{\u}(\cdot \vee \rho\tau_{\tilde{h}} \wedge \tau_{\tilde{h}})
	\in C([0,T]; \b H^{2+\gamma}) \cap L^2(0,T; \b H^4), \quad \P\text{-a.s.} \label{eq: utilde reg on Ihrho}
\end{align} 

For a.e.-$ \omega \in \Omega $, thanks to the higher-order (particularly $ \b H^{2+\gamma} $) regularity of $ \tilde{\u}(\omega) $ given in \eqref{eq: utilde reg on Ihrho}, this local solution can be extended to a unique $ \b H^2 $-maximal solution, still denoted as $ \tilde{\u}(\omega)$, of \eqref{eq: u} in $ \b H^2 $ on $ [\rho \tau_{\tilde{h}}, \tilde{\tau}^*_1) $ where $ \tau_{\tilde{h}}(\omega) < \tilde{\tau}^*_1(\omega) \leq T $ and 
\begin{equation}\label{eq: tilde u in sphere}
	|\tilde{\u}(t,x,\omega)| = 1 , \quad \text{a.e-}(t,x) \in [0,\tilde{\tau}^*_1(\omega)) \times \T^2.
\end{equation}
Note that $ \tau_{\tilde{h}} $ is an $ \tilde{\b F} $-stopping time. 
Then by Corollary \ref{Coro: tau stopping} and Lemma \ref{Lemma: energy moments}, the maximal time $ \tilde{\tau}^*_1 $ is an $ \tilde{\b F} $-stopping time, and for any $ \tilde{\b F} $-stopping time $ \tilde{\zeta} \in [\tau_{\tilde{h}}, \tilde{\tau}_1^*) $, the (stopped) extended solution process $ \tilde{\u}(\cdot \wedge \tilde{\zeta}): [0,T] \times \Omega \to \b H^1 $ is $ \tilde{\b F} $-progressive, satisfying
\begin{equation}\label{eq: tilde u H1 moment part 2}
\begin{aligned}
	&\E \left[ \sup_{s \in [0,T]}\left|\nabla \tilde{\u}(\theta_{\tau_{\tilde{h}}}(s) \wedge \tilde{\zeta}) \right|^2_{\b L^2} + \int_{\tau_{\tilde{h}}}^{\tilde{\zeta}} |\tilde{\u} \times \Delta_{\tilde{A}} \tilde{\u}|_{\b L^2}^2(s) \ ds\right] \\
	&\leq \tilde{c}\la \E \left[ \left|\nabla \tilde{\u}(\tau_{\tilde{h}}) \right|_{\b L^2}^2 + |\tilde{A}|_{C([0,T]; \b L^2)}^2 \right] + q^4(2) \ra < \infty, 
\end{aligned}	
\end{equation}
for some constant $ \tilde{c}>0 $.

\subsection{Singularities at $ \tilde{\tau}_1^* $}\label{Section: step 4} 
If the smallness condition \eqref{eq: Du+A eps} is satisfied at $ \tilde{\tau}^*(\omega) $ and a positive $ r(\omega) \leq r_0 $ for some $ \omega \in \Omega $, then $ \tilde{\u}(\omega) $ has improved regularity with uniform (in time) bound on $ [\tau_{\tilde{h}}(\omega), \tilde{\tau}^*_1(\omega)) $ as in Proposition \ref{Prop: u H2 est global}(ii), which allows it be extended smoothly up to (including) $ \tilde{\tau}^*_1(\omega) $, contradicting the maximality of $ \tilde{\tau}^*_1(\omega) $.  
Thus, \eqref{eq: Du+A eps} is not satisfied at $ \tilde{\tau}^*(\omega) $, for a.e.-$ \omega \in \Omega $. 
In other words, there must exist a singular set $ S^*_1(\omega) $ given by
\begin{align*}
	S^*_1(\omega)
	&= \Big\{ x_n \in \T^2  
	\mid 
	\limsup_{t \nearrow \tilde{\tau}^*_1(\omega)} \la |\nabla \tilde{\u}(t)|_{L^2(B_{2r}(x_n))}^2 + |\tilde{A}(t)|_{L^2(B_{2r_{\tilde{A}}}(x_n))}^2 \ra(\omega) \\
	&\qquad \qquad \qquad> \eps^2,  \quad
	\forall r \in (0,r_0] \Big\}.
\end{align*} 
By \eqref{eq: A eps} and the way that $ r_{\tilde{A}} $ was chosen in Section \ref{Section: step 2},  
\begin{equation}\label{eq: limsup_t |Du|_L2Br}
	\limsup_{t \nearrow \tilde{\tau}^*_1(\omega)} \left|\nabla \tilde{\u}(t,\omega) \right|_{L^2(B_{2r}(x_n))}^2 > \frac{\eps^2}{2},
\end{equation}
for any $ r(\omega) \in (0, r_0] $ and thus for any $ r > 0 $.  

Now we identify the limit of $ \tilde{\u}(t) $ as $ t \to \tilde{\tau}^*_1 $. 
Let $ \{t_i\}_{i \in \mathbb{N}} $ be a sequence of increasing positive $ \tilde{\b F} $-stopping times such that $ 0 \leq t_i < \tilde{\tau}_1^* $ for all $ i $ and $ t_i \to \tilde{\tau}^*_1 $, $ \P $-a.s. as $ i \to \infty $. 
For instance, we can take 
\begin{align*}
	t_i = \inf\{ t \in [\tilde{\tau}_h, T] : |\tilde{\u}(t)|_{\b H^2} > i \}. 
\end{align*}
By the on-sphere constraint \eqref{eq: tilde u in sphere} and the moment estimates \eqref{eq: tilde u H1 moment part 1} and \eqref{eq: tilde u H1 moment part 2}, we have 
\begin{equation}\label{eq: E[|u(ti)| H1 norm] finite}
	\sup_{i \in \b N} \E \left[ \left|\tilde{\u}(t_i) \right|_{\b H^1}^2 \right]
	\leq \E \left[ \sup_{i \in \b N} \left|\tilde{\u}(t_i) \right|_{\b H^1}^2 \right] < \infty. 
\end{equation}
This implies that there exists a subsequence of $ (\tilde{\u}(t_i))_i $ that converges weakly in $ L^2(\Omega;\b H^1) $ to some $ \tilde{\u}^* \in L^2(\Omega; \b H^1) $. 
Moreover, fix $ \omega $, there exists a subsequence of $ (\tilde{\u}(t_i))_i(\omega) $ that converges weakly in $ \b H^1 $ to some $ \tilde{\v}^* \in \b H^1 $.  
In fact, these weak convergences for the entire sequences $ (\tilde{\u}(t_i))_i $ and $ (\tilde{\u}(t_i))_i(\omega) $, where the limit $ \tilde{\u}^* $ is $ \tilde{\mc F}_{\tilde{\tau}_1^*} $-measurable and $ \tilde{\u}^*(\omega)=\v^* $. 
We start from verifying that $ (\tilde{\u}(t_i))_{i} $ is a Cauchy sequence in $ L^2(\Omega; \b L^2) $:
\begin{align*}
	\E \left[ \left|\tilde{\u}(t_{i}) - \tilde{\u}(t_{j}) \right|_{\b L^2}^2 \right]
	&\leq \tilde{c} \E \left[ \left| \int_{t_{i}}^{t_{j}} \tilde{\u} \times \Delta_{\tilde{A}} \tilde{\u} \ ds \right|_{\b L^2}^2 \right] \\
	&\leq \tilde{c} \E \left[ \la \int_{t_{i}}^{t_{j}} |\tilde{\u} \times \Delta_{\tilde{A}} \tilde{\u} |_{\b L^2}^2 \ ds \ra^2 \right]^\frac{1}{2} \la \E \left[ |t_{j}-t_{i}| \right]\ra^\frac{1}{2},
\end{align*}
for any $ i, j \in \b N $, where 
\begin{align*}
	\lim_{i \to \infty} \E \left[ \int_{\tau_{\tilde{h}}}^{t_i} |\tilde{\u} \times \Delta_{\tilde{A}} \tilde{\u} |_{\b L^2}^2 \ ds \right] 
	&=\E \left[ \int_{\tau_{\tilde{h}}}^{\tilde{\tau}_1^*} |\tilde{\u} \times \Delta_{\tilde{A}} \tilde{\u} |_{\b L^2}^2 \ ds \right] < \infty .
\end{align*}
The inequalities above hold as in the proof of Lemma \ref{Lemma: energy moments} using the equation \eqref{eq: u}.  
Therefore, $ \tilde{\u}(t_i) $ converges strongly in $ L^2(\Omega; \b L^2) $ to some $ \tilde{\u}' \in L^2(\Omega; \b L^2) $. 
Then we have the following consequences. 
\begin{enumerate}
	\item[(i)] 
	$ \tilde{\u}' = \tilde{\u}^* $ by the uniqueness of weak limit in $ L^2(\Omega; \b L^2) $. 
	This applies to every subsequence that is weakly convergent in $ L^2(\Omega; \b H^1) $, leading to weak convergence in $ L^2(\Omega, \b H^1) $ of the entire sequence $ (\u(t_i))_i $. 
	
	\item[(ii)] 
	There exists a further subsequence of $ (\u(t_i))_i $ that converges in $ \b L^2 $ to the unique limit $ \tilde{\u}^* $, $ \P $-a.s. 
	The random variables $ \tilde{\u}(t_i) $, $ i \in \b N $, are $ \tilde{\mc F}_{\tilde{\tau}_1^*} $-measurable, and so does their $ \P $-a.s. limit $ \tilde{\u}^* $.
\end{enumerate}
For fixed $ \omega $, $ (\tilde{\u}(t_i)(\omega))_i $ is similarly a Cauchy sequence in $ \b L^2 $, since $ \tilde{\u}(t_i)(\omega) $ solves a deterministic equation, $ \int_{t_i}^{t_j} |\tilde{\u} \times \Delta_{\tilde{A}} \tilde{\u}|_{\b L^2}^2 \ ds < \infty $, $ \P $-a.s. and $ t_i \nearrow \tilde{\tau}_1^* $, $ \P $-a.s. 
Similar to consequence (i), we obtain that the entire sequence $ (\tilde{\u}(t_i))_i(\omega) $ converges weakly in $ \b H^1 $ to some $ \tilde{\v}^* \in \b H^1 $. 
Then by consequence (ii), $ (\tilde{\u}(t_i))_i(\omega) $ converges to $ \tilde{\u}^*(\omega) = \v^* $ weakly in $ \b L^2 $. 
Since this argument applies for a.e.-$ \omega $, we obtain the a.s. weak convergence 
\begin{align*}
	\tilde{\u}(t_i) \rightharpoonup \tilde{\u}^* \text{ in } \b H^1, \quad \P\text{-a.s.}
\end{align*}
for any sequence of stopping times $ t_i \nearrow \tilde{\tau}_1^* $, $ \P $-a.s.

By the weak lower semi-continuity of the $ \b H^1 $-norm and the inequality \eqref{eq: limsup_t |Du|_L2Br},
\begin{align}
	\E \left[ |\nabla \tilde{\u}^*_1|^2_{\b L^2} \right]
	&= \E \left[ \lim_{r \to 0} |\nabla \tilde{\u}^*_1|^2_{L^2(\T^2 \setminus \cup_{n} B_{2r}(x_n))} \right] \nonumber\\
	&\leq \liminf_{i\to \infty} \ \E \left[ \lim_{r \to 0} |\nabla \tilde{\u}(t_i)|^2_{L^2(\T^2 \setminus \cup_{n} B_{2r}(x_n))} \right] \nonumber \\
	&= \liminf_{i \to \infty} \ \E \left[ |\nabla \tilde{\u}(t_i)|^2_{\b L^2} - \lim_{r \to 0} \sum_{n \leq |S^*_1|} |\nabla \tilde{\u}(t_i)|^2_{B_{2r}(x_n)} \right] \nonumber\\
	&\leq \liminf_{i \to \infty} \ \E\left[ |\nabla \tilde{\u}(t_i)|^2_{\b L^2} \right] - \frac{\eps^2}{2}\E \left[ |S^*_1| \right]. \label{eq: E|Du*|<=supt|Du|-eps} 
\end{align}
Then by the fact $ \nabla \tilde{\u}^*_1 \in L^2(\Omega; \b H^1) $ and \eqref{eq: E[|u(ti)| H1 norm] finite}, 
we have $ \E[|S^*_1|] < \infty $. 

Alternatively, by Lemma \ref{Lemma: local energy}(i), we can also show that the $ \P $-a.s. $ \b L^2 $-convergent sequence has a $ \P $-a.s. weak $ \b H^1 $-convergent subsequence, and similarly,
\begin{equation}\label{eq: |Du*|<=supt|Du|-eps}
	|\nabla \tilde{\u}^*_1|^2_{\b L^2}
	\leq \liminf_{j \to \infty} |\nabla \tilde{\u}(t_{i_j})|^2_{\b L^2} - \frac{\eps^2}{2} |S^*_1| 
	\leq \sup_{t \in [0,\tilde{\tau}_1^*)} |\nabla \tilde{\u}(t)|^2_{\b L^2} - \frac{\eps^2}{2} |S^*_1|. 
\end{equation}

To improve the a.s. $ \b H^1 $-weak convergence, we note that the condition \eqref{eq: Du+A eps} is satisfied if the supremum is taken over $ x \in \T^2 \setminus S_1^* $ instead of $ x \in \T^2 $. 
This implies that away from the random singular locations, the restriction $ \tilde{\u} \vert_{\T^2 \setminus S_1^*} $ can be extended continuously up to $ \tilde{\tau}_1^* $, yielding $ \tilde{\u}(\omega) \in C([\tau_{\tilde{h}},\tilde{\tau}^*_1](\omega); H^2(\T^2 \setminus S_1^*(\omega))) $ for a.e.-$ \omega \in \Omega $. 
From the continuity in time, we have the following strong convergence
\begin{align*}
	\tilde{\u}(t_i)(\omega) \to \tilde{\u}^*(\omega) \text{ in } H^1(\T^2 \setminus S_1^*(\omega)),  
\end{align*}  
for a.e.-$ \omega \in \Omega $.


\subsection{Iterate to construct solution in $ C([0,T]; \b L^2) \cap L^\infty(0,T; \b H^1) $}\label{Section: step 5} 
We repeat Section \ref{Section: step 1} -- \ref{Section: step 4} using the initial point $ \tilde{\u}^*_1 $ to obtain a similar local $ \b H^1 $-solution.  
For clarity, we provide details of the second iteration, particularly for 
the time shift, 
and the random radius associated with $ \tilde{\u}^*_1 $. 

As in Section \ref{Section: step 1}, for a.e.-$ \omega \in \Omega $, there is a sequence of smooth approximations \smash{$ \{\tilde{\u}_1^{(k)}(\omega)\} $} converging to $ \tilde{\u}_1^*(\omega) $ strongly in $ \b H^1 $. 
The approximating sequence is constructed from a (deterministic) mollification of \smash{$ \tilde{\u}^*_1 $}, thus 
for every $ k \in \b N $, $ \smash{\tilde{\u}^{(k)}_1}: \Omega \to C^\infty(\T^2; \St) $ is similarly \smash{$ \tilde{\mc F}_{\tilde{\tau}_1^*} $}-measurable and the convergence is $ \P $-a.s. 
Similar to \eqref{eq: u0k eps}, now we select a $ \tilde{\mc F}_{\tilde{\tau}_1^*} $-measurable random radius $ \tilde{r}_1>0 $ such that 
\begin{align*}
	\sup_{k \in \b N} \sup_{x \in \T^2} \left| \nabla \tilde{\u}_1^{(k)} \right|_{L^2(B_{4\tilde{r}_1(x)})}^2 \leq \frac{\eps^2}{4}, \quad \P\text{-a.s.} 
\end{align*}
Let $ (\tilde{\u}^{(k)}_{*}, [\tilde{\tau}_1^*, \tilde{\tau}^{(k)}_{*})) $ denote the unique $ \b H^2 $-maximal solution of \eqref{eq: v(t,y)} with the parameter $ \tilde{A} $ and the initial data $ \tilde{\u}_1^* $ at $ \tilde{\tau}_1^* $. 
Recall $ h $ in \eqref{eq: h(t,u0,r0)}. Let 
\begin{align*}
	\tilde{h}_*(s)	
	&:= h\la s; \tilde{A}(\theta_{\tilde{\tau}_1^*}),\tilde{\u}_1^* \ra, 
	\quad s \geq 0,
\end{align*}
which depends only on \smash{$ ( \tilde{\u}_1^{(k)} )_{k \in \b N} $}, $ \tilde{r}_1 $ and the restriction of $ \tilde{A} $ to the interval \smash{$ {[\tilde{\tau}_1^*,\theta_{\tilde{\tau}_1^*}(s)]} $}. 
Similarly, the process $ \tilde{h}_* $ is independent of $ k $, increasing, continuous and $ \tilde{\mathbb{F}}_{\tilde{\tau}_1^*} $-adapted, where 
$ \tilde{\mathbb{F}}_{\tilde{\tau}_1^*} = (\tilde{\mc F}_{\tilde{\tau}_1^*+s})_{s \geq 0} $. 
Then we have an $ \tilde{\b F}_{\tilde{\tau}_1^*} $-stopping time
\begin{align*}
	\tilde{\tau}_{\tilde{h}_*} 
	:= \inf \left\{ s \geq 0 : \tilde{h}_*(s) s = \frac{\eps^2}{4c} \right\}, 
\end{align*}
where \smash{$ \tilde{\tau}_{\tilde{h}_*} > 0 $}, $ \tilde{\P} $-a.s. 
As a result, $ \tilde{\tau}_{\tilde{h}_*} + \tilde{\tau}_1^* $ is an $ \tilde{\b F} $-stopping time since it is equivalent to 
\begin{align*}
	\inf \left\{ t \in [\tilde{\tau}_1^*,T] : \tilde{h}_*(t-\tilde{\tau}_1^*) (t-\tilde{\tau}_1^*) = \frac{\eps^2}{4c} \right\}.
\end{align*}
Using again the maximality of $ \tilde{\tau}_*^{(k)} $ and Proposition \ref{Prop: u H2 est global}, we have 
$ \tilde{\tau}_{\tilde{h}_*} + \tilde{\tau}_1^* < \tilde{\tau}^{(k)}_* $, $ \tilde{\P} $-a.s. for every $ k \in \b N $. 	
Define
\begin{align*}
	\tilde{\u}^{(k)}_{\tilde{h}_*}(s) := \tilde{\u}^{(k)}_*(\theta_{\tilde{\tau}_1^*}(s \wedge \tilde{\tau}_{\tilde{h}_*})), \quad s \geq 0. 
\end{align*}
Then by Lemma \ref{Lemma: local energy}, a similar condition to \eqref{eq: Duk+A eps} holds for $ \la \smash{\tilde{\u}^{(k)}_{\tilde{h}_*}}, \tilde{r}_1, \tilde{A}, r_{\tilde{A}} \ra $ instead of $ (\smash{u^{(k)}_h}, r_0, A, r_A) $. 

As in Section \ref{Section: step 2}, \smash{$ \tilde{\u}^{(k)}_{\tilde{h}_*} $} is uniformly bounded in $ L^2(\tilde{\Omega}; U_0^\gamma) $ for any $ \gamma \in (0,\frac{1}{2}) $. 
Also, let
\begin{align*}
	\tilde{W}_*(s) := \tilde{W}(\tilde{\tau}_1^*+s) - \tilde{W}(\tilde{\tau}_1^*), \quad s \geq 0. 
\end{align*}
Then $ \tilde{W}_* $ is a new cylindrical Wiener process on $ (\Omega, \mc F, \P) $ independent of $ \tilde{\mc F}_{\tilde{\tau}_1^*} $, and $ \tilde{\Y}_*(\cdot) = \tilde{\Y}(\tilde{\tau}_1^*+\cdot) $ is the unique solution of \eqref{eq: mc Y} driven by $ \tilde{W}_* $ starting from $ \tilde{\Y}_*(0) = \tilde{\Y}(\tilde{\tau}_1^*) $. 
Applying again the version of Skorohod theorem in \cite{BertiPratelliRigo}, 
on the same probability space $ (\Omega, {\mc F}, \P) $, there exist \smash{$ V_0^{-1,\sigma} \times [0,1] $}-valued random variables
$ (\hat{\u},\hat{\Y},\hat{W},\hat{r}_1) $ and \smash{$ \{(\hat{\u}^{(k)}, \hat{\Y}^{(k)}, \hat{W}^{(k)}, \hat{r}_1^{(k)})\}_{k \in \b N} $} 
such that 
\begin{align*}
	\mathcal{L}((\hat{\u}^{(k)}, \hat{\Y}^{(k)}, \hat{W}^{(k)}, \hat{r}_1^{(k)})) = \mathcal{L}((\tilde{\u}^{(k)}_{\tilde{h}_*}, \tilde{\Y}_*, \tilde{W}_*, \tilde{r}_1)) 
	\text{ on $ V_0^{-1,\sigma} \times [0,1] $},
\end{align*}
for all $ k \in \b N $, and the a.s. convergence (ii) holds similarly.
Let $ \hat{\mc F}_t $ be the augmentation of $ \sigma(\smash{\hat{\u}^{(k)}}(s), \hat{W}^{(k)}(s), \hat{r}^{(k)}_1, \hat{\u}(s), \hat{W}(s), \hat{r}_1 ; k \in \b N, s\leq t) $ and $ \hat{\b F} := \smash{(\hat{\mc F}_t)_{t \in [0,T]}} $. 
In particular, $ \hat{r}_1^{(k)} $ is $ \hat{\mc F}_0 $-measurable, and $ \hat{\Y}^{(k)} $ is the unique solution of \eqref{eq: mc Y} driven by $ \hat{W}^{(k)} $ with initial distribution $ \mathcal{L}(\hat{\Y}^{(k)}(0)) = \mathcal{L}(\tilde{\Y}_*(0)) $. 

The rest of Section \ref{Section: step 2} -- \ref{Section: step 4} can be repeated, giving a solution $ \hat{\u} $ and an $ \hat{\b F} $-stopping time $ \hat{\tau}_2^*>0 $, $ \P $-a.s. with the following properties:
\begin{itemize}
	\item 
	for any $ \hat{\b F} $-stopping time $ \hat{\zeta} \in [0,\hat{\tau}_2^*) $, $ \hat{\u}(\cdot \wedge \hat{\zeta}) $ is 
	$ \hat{\b F} $-progressively measurable, 
	bounded in $ L^2(\Omega, L^\infty(0,T; \b H^1)) $ with $ |\hat{\u}(t,x,\omega)|=1 $ a.e. and 
	takes values in $ C((0,T]; \b H^2) $, $ \P $-a.s. 
	
	\item 
	for a.e.-$ \omega \in \Omega $, there exists a singular set $ S_2^*(\omega) $ such that $ \E[|S_2^*|] < \infty $, and 
	for every $ x \in S_2^*(\omega) $, 
	\begin{align*}
		\limsup_{t \nearrow \tilde{\tau}_1^*(\omega)} |\nabla \tilde{\u}(t,\omega)|^2_{L^2(B_{2r}(x))} > \frac{\eps^2}{2}, \quad \forall r > 0,
	\end{align*}
	
	\item 
	$ \hat{\u}(t) $ converges to an $ \hat{\mc F}_{\hat{\tau}_2^*} $-measurable random variable $ \hat{\u}^*_2 $ weakly in $ L^2(\Omega; \b H^1) $, as $ t \nearrow \hat{\tau}_2^* $. 
\end{itemize}

To concatenate the two solutions $ \tilde{\u} $ and $ \hat{\u} $, 
for every $ t \geq 0 $, let $ \bar{\mc F}_t $ be the completion of the $ \sigma $-algebra
\begin{align*}
	\sigma \la \tilde{\mc F}_t \cup \left\{ \{ \tilde{\tau_1}^* \in [s,t] \} \cap B: B \in \hat{\mc F}_{t-s}, \ s \in [0,t) \right\} \ra
\end{align*}
By definition, $ \tilde{\mc F}_t \subseteq \bar{\mc F}_t $, thus $ \bar{\tau}_1^*:= \tilde{\tau}_1^* $ is also an $ \bar{\b F} $-stopping time where $ \bar{\b F} := (\bar{\mc F}_t)_{t \geq 0} $. 
Similarly, $ \bar{\tau}_2^* := \tilde{\tau}_1^*+\hat{\tau}_2^* $ is an $ \bar{\b F} $-stopping time. 
Moreover, the increments of $ \hat{W} $ are all independent to $ \bar{\mc F}_{\tilde{\tau}_1^*} $. 
Then the process $ \bar{W} $ given by 
\begin{align*}
	\bar{W}(t) 
	:= \tilde{W}(t) \mathbbm{1}(\tilde{\tau}_1^* >t) 
	+ \la \tilde{W}(\tilde{\tau}_1^*) + \hat{W}(t) \ra \mathbbm{1}(\tilde{\tau}_1^* \leq t), \quad t \in [0,T],
\end{align*}
is a cylindrical Wiener process adapted with respect to $ \bar{\b F} $. 
The processes $ (\tilde{Y}, \hat{Y}) $ and $ (\tilde{A}, \hat{A}) $ can be similarly pasted together. 
As a result, we arrive at a $ \bar{\b F} $-progressively measurable solution on $ (\Omega, \bar{\mc F}, \P) $: $ (\bar{\u}, [0,\tilde{\tau}_1^*+ \hat{\tau}_2^*)) $ defined by
\begin{align*}
	\bar{\u}(t) 
	= \begin{cases}
		\tilde{\u}(t), &\quad t \in [0,\tilde{\tau}_1^*) = [0,\bar{\tau}_1^*), \\
		\hat{\u}(t-\tau_1^*)-\hat{\u}(0) + \tilde{\u}_1^*, &\quad t \in [\tilde{\tau}_1^*, \tilde{\tau}_1^*+\hat{\tau}_2^*) = [\bar{\tau}_1^*, \bar{\tau}_2^*).
	\end{cases}
\end{align*}
For simplicity, let us omit all accents. 

Continuing this procedure, we produce a weak solution $ \u(t) $ in $ \b H^1 $ on $ [0, \tau_\ell^*) $ by the end of the $ \ell $-th iteration. 
Then we only need to show that this solution $ u $ is global. 
Let $ N_\ell := |S_\ell^*| $ denote the (random) number of singular points as $ t \to \tau^*_{\ell} $. 
Recall \eqref{eq: E[|u(ti)| H1 norm] finite}, \eqref{eq: |Du*|<=supt|Du|-eps} and Lemma \ref{Lemma: local energy} -- \ref{Lemma: energy moments}. 
Inductively,
\begin{align*}
	|\nabla \u^*_{\ell}|_{\b L^2}^2
	&\leq
	\la |\nabla \u_0|_{\b L^2}^2 + c \int_0^{\tau^*_{\ell}} |\nabla \cdot A + A^2|_{\b L^2}^2 \ ds \ra e^{c |A|^2_{L^2(0,\tau^*_{\ell}; \b L^\infty)}} \\
	&\quad - \frac{\eps^2}{2}\sum_{j=1}^{\ell} N_j, \quad \P\text{-a.s.} \\
	\E \left[ |\nabla \u^*_{\ell}|_{\b L^2}^2 \right]
	&\lesssim |\nabla \u_0|_{\b L^2}^2 + \E \left[ |A|^2_{C([0,T]; \b L^2)} \right] + q^4(2) - \frac{\eps^2}{2} \E \left[ \sum_{j=1}^\ell N_j \right].
\end{align*}
By construction, we have $ \tau^*_\ell \leq T $ and thus
\begin{align*}
	\sup_{\ell \in \b N} \la |\nabla \u^*_{\ell}|_{\b L^2}^2 + \sum_{j=1}^\ell N_j \ra < \infty, \quad \P\text{-a.s.}
\end{align*}
The sequence $ (\tau_\ell^*)_{\ell} $ is positive and increasing. 
Let $ \tau^*:= \lim_{\ell \to \infty} \tau^*_\ell $, and let $ \u^* $ denote the weak limit of $ \u_\ell^* $ in $ L^2(\Omega; \b H^1) $. 
If $ \tau^* < T $, then no more iteration can be performed with $ |\u_\ell^*|_{\b H^1} \to \infty $, leading to a contradiction. 
Thus, $ \tau^* = T $ and we arrive at a global weak solution $ \u $ in $ H^1(0,T;\b L^2) \cap L^\infty(0,T; \b H^1) $ of \eqref{eq: u} with a finite number of time-space singular points, $ \P $-a.s.

\subsection{Value of $ \eps $}\label{Section: step 6} 
We now provide a geometric characterization of singularities and deduce a precise numerical value for $\eps$. 
Again, for fixed $ \omega \in \Omega $, after translation and dilation we may $\u=\u(\omega) \in C([-4,0); H^2(\bar{B_2} \setminus \{0\}; {\b S}^2))$.
If $(0,0)$ is a singularity then by virtue of Proposition \ref{Prop: u H2 est Pr},
\[
|\nabla \u(t_k)|^2_{L^2(B_{r_k}(x_k))} = \sup_{x \in B_1} |\nabla \u(t_k)|^2_{L^2(B_{r_k}(x))} =\frac{\eps^2}{4}
\]
for suitable sequences $t_k \nearrow 0$, $x_k \to 0$ and $r_k \searrow 0$. 
Moreover, invoking Lemma \ref{Lemma: local energy} we find $0<\delta_0 \le 1/4$ such that
\[
\sup_{x \in B_{\frac{1}{2}}(x_k)}  \int_{B_{\frac{r_k}{2}}(x)}|\nabla \u(t)|^2 dx \le \frac{\eps^2}{2}, \quad \forall t \in [t_k -  r_k^2 \delta_0,t_k],
\]
for sufficiently large $k$. 
The blow-up solutions $\u_k$ given by
\[
\u_k(t,x)=\u(t_k +r_k^2 t, x_k+r_k x), \quad (t,x) \in [-\delta_0,0] \times \R^2
\]
admit for $x \in B_{1/2r_k}$ and $t \in [- \delta_0/2, 0]$ a uniform higher order bounds
according to Proposition \ref{Prop: u H2 est Pr}.
We consider $\u_k$ as a solution of the perturbed Landau-Lifshitz equation
\[
\partial_t \u_k = \u_k \times \left( \alpha \partial_t \u_k - \Delta \u_k \right) + F_k
\]
where $|F_k(t)|_{\b L^2} = O(r_k)$ uniformly for all admissible $t$. 
It follows from the energy inequality for $\u$ that 
$\int_{-\delta_0}^0 \int_{\R^2} |\partial_t \u_k|^2 \, dxdt \to 0$ as $k \to \infty$, hence
$\v_k=(\partial_t \u_k )(\tau_k)$ and $\w_k = F_k(\tau_k)$ converge to zero in $L^2(\R^2)$
for some sequence $\tau_k \nearrow 0$.
Note that $\tilde{\u}_k = \u_k(\tau_k)$ is an almost harmonic map in the sense that
\[ 
\tilde{\u}_k \times \Delta \tilde{\u}_k = \alpha \, \tilde{\u}_k \times \v_k - \v_k + \w_k
\]
and subconvergence strongly in $H^1_{\rm loc}(\R^2;\R^3)$ to a harmonic map $\tilde{\u}$ of finite energy in $\R^2$.  To show that $\tilde{\u}$ is non-constant we invoke Lemma \ref{Lemma: local energy} for $u_k$
\[
\int_{B_1} |\nabla \u_k(0)|^2 \, dx -  \int_{B_{2}} |\nabla \u_k(\tau_k )|^2 \, dx  \to 0 \quad \text{as $k \to \infty$,}
\]
so that by local strong convergence $\int_{B_{2}} |\nabla \tilde{\u}|^2 \, dx >0$. By virtue of the well-known theory about harmonic maps $\frac 1 2 \int_{\R^{2}} |\nabla \tilde{\u}|^2 \, dx$ is a positive multiple of $4 \pi$. 
Hence by letting $s_k = t_k+r^2_k \tau_k \to 0$ we have for arbitrary $r_0>0$
\[
\int_{B_2(0)} |\nabla \u(s_k)| \, dx \geq \int_{B_1(x_k)} |\nabla \u(s_k)| \, dx 
= \int_{B_{1/r_k}} |\nabla \tilde{\u}_k| \, dx \geq   \int_{B_{r_0}} |\nabla \tilde{\u}_k| \, dx 
\]
for $k > k_0$ depending on $r_0$
which implies by strong $L^2(B_{r_0})$ compactness of $\nabla \tilde{\u}_k$
\[
\liminf_{k \to \infty} \frac 1 2 \int_{B_2(0)} |\nabla \u(s_k)|^2 \geq \frac 1 2 \int_{B_{r_0}} |\nabla \tilde{\u}| \, dx = 4\pi |q| +o(1). 
\]
This concludes the proof of Theorem \ref{Theorem: Struwe u} and Corollary \ref{Corollary: Struwe u iteration}.

\begin{remark}
	If the noise is finite-dimensional with commuting $ G_j $ matrices (e.g. $ g_j(t) = f(x) h $ for a fixed constant vector $ h \in \R^3 $, $ \forall j \leq N $),  
	then $ A $ can be expressed in terms of $ \sin(W_j) $ and $ \cos(W_j) $ as in Example \ref{Example: 1D}, so that for a sufficiently smooth $ f $. 
	As a result, we can obtain a non-random bound for $ |A(\omega)|_{C([0,T],\b H^\sigma)} $ for a.e.-$ \omega \in \Omega $. 
	This could simplify the derivations above, since an explicit expression of $ \tau_h $ in terms of $ r_0 $ would be available. 
\end{remark}

\section{Proof of uniqueness}\label{Section: Proof uniqueness}

We sketch the proof of Theorem \ref{Theorem: u unique} in this section. 
For simplicity, we fix $ \omega \in \Omega $ such that $ A(\omega) $ is sufficiently regular, and write $ \u = \u(\omega) $ (resp., $ t = t(\omega) $). 

Right-continuity of Dirichlet energy of $ \u $ is the key to control the most singular nonlinear term by decomposing it into a smooth and small part. Thereby it implies improved space-time regularity 
\[
\nabla \u \in L^4(0,t; \b L^4) \quad \text{}
\]
on small forward in time intervals $ (0,t] $, 
which is sufficient to perform a local Gronwall argument, cf. \cite{Chen_LLBVP}. 
In fact, by $L^2$-continuity and translation properties of $A$, it is enough to establish local uniqueness near $t=0$. 

We present the crucial step towards regularity using right-continuity and the decomposition argument, which is not properly described in the corresponding literature on LLG as pointed out in \cite{Harpes}.

\begin{lemma}\label{Lemma: D2u in L4/3}
	Under the assumption of Theorem \ref{Theorem: u unique}, there exists (random) $ t>0 $ such that 
	$ \Delta \u \in L^2(0,t;\b L^{\frac{4}{3}}) $.
\end{lemma}

\begin{proof}[Sketch of proof of Lemma \ref{Lemma: D2u in L4/3}] 
	The main idea is to interpret the magnetic LLG (in Gilbert form) as an
	inhomogeneous magnetic heat equation 
	\[
	\partial_t \u - \alpha \left( \Delta \u - \nabla^\perp  B_\eps : \nabla \u \right) = \boldsymbol f.
	\]
	where the magnetic term in $B_\eps$ is arising from the geometric nonlinearity - not the transformation. 
	The key observation (cf. \cite{Freire}) is that by virtue of Wente's inequality, the magnetic Laplacian
	above is uniformly elliptic in the sense that 
	\[
	- \langle \Delta \u + \nabla^\perp  B_\eps : \nabla \u, \u \rangle_{\b L^2} \geq \frac{1}{2} \left(1- c \, |\nabla B_\eps|_{\b L^2} \right) |\nabla \u |_{\b L^2}^2,  
	\]
	where the operator $ -(\Delta \u + \nabla^\perp B_\eps: \nabla \u)(s) $ is coercive for any $ s \in [0,t] $, provided $ |\nabla B_\eps|_{L^\infty(0,t;\b L^2)}< \eps$ and parabolic theory applies.

	In the case of the our magnetic LLG equation, the source term $\boldsymbol{f}$ contains 
	\begin{enumerate}[leftmargin=*,parsep=1pt]
		\item[(i)]
		the Hamiltonian terms that in Gilbert form are given by $ \u \times \partial_t \u \in L^2(0,T; \b L^2)$, 
		
		\item[(ii)]
		the pertubations arising from $A$ that are bounded as follows
		\[
		|\nabla A|+ |A||\nabla \u| + |A|^2 \in L^\infty(0,T; \b L^p), \quad \forall p<2,
		\]
		
		\item[(iii)]
		and the terms arising from the $\mathfrak{so}(3)$ valued Hodge decomposition
		\[
		\Psi = \nabla^\perp B + \Phi, \quad \nabla^\perp \cdot \Phi =0,  
		\] 
		or in the H\'elein factorization of the geometric nonlinearity 
		\[
		|\nabla \u|^2 \u= \Psi :\nabla \u, \quad  
		\Psi = \u \otimes \nabla \u - \nabla \u \otimes \u.
		\]
	\end{enumerate}

	For terms in (iii), using the equation it follows that 
	\[
	|\nabla \Phi|_{\b L^2} = |\nabla \cdot \Psi|_{\b L^2} 
	\]
	which is uniformly bounded in time.
	Hence, by Sobolev embedding, $\Phi \in L^\infty(0,t; \b L^p)$ for all $p<\infty$ and will be attributed to $\boldsymbol{f}$. The final step is a decomposition of $B$. 
	According to our assumption $\Psi$ and 
	hence $B$ is right-continuous in time with values in $\b H^1$. This 
	implies by letting 
	\[
	B_\eps(t) = B(t) - \Pi_n B(0) 
	\]
	where $\Pi_n$ is a suitable spectral truncation so that $B-B_\eps \in C^\infty$ will 
	contribute to $\boldsymbol{f}$, while
	\begin{eqnarray*}
		|\nabla B_\eps(t) |_{\b L^2} &\le&  |\nabla( B(t)-B(0) )|_{\b L^2} + | \nabla(B(0)-\Pi_n B(0)) |_{\b L^2} \\
		&\le&  | \nabla(B(t)-B(0))|_{\b L^2} + \| 1-\Pi_n \| | \nabla \u(0) |_{\b L^2} < \eps
	\end{eqnarray*}
	for $n$ sufficiently large and $t$ sufficiently small. Since 
	\[
	\boldsymbol{f} \in L^\infty(0,t; \b L^p) \cap L^2(0,t; \b L^2) \cap L^2(0,t; \b L^{\frac{4}{3}})
	\]
	on finite time intervals, the claim follows from parabolic theory.
\end{proof}

In a second step based on the Landau-Lifshitz formulation, integrability in time can be improved to the power $4$, and the requisite gradient regularity $\nabla \u \in L^4(0,t; \b L^4)$ follows by Sobolev embedding. 
To conclude the proof of the theorem, it is routine to derive the  following differential inequality for the increment $\w=\u - \v$
of two weak solutions $\u$ and $\v$ in this class 
\[
\partial_t |\w|^2_{\b L^2} + 2 \alpha |\nabla_A \w|^2_{\b L^2} 
\lesssim \int_{\T^2}\la |\nabla_A U | |\nabla_A \w||\w| + |\nabla_A U|^2 |\w|^2 \ra dx
\]
with $U=(\u, \v)$. 
By Young's and Ladyzhenskaya's interpolation inequality, 
we obtain
\begin{align*}
	\partial_t |\w|^2_{\b L^2} +  \alpha |\nabla \w|^2_{\b L^2} 
	&\lesssim \left( |\nabla \u|_{\b L^4}^2 + |A|_{\b L^4}^2 \right) |\w|_{\b L^4}^2  \\
	&\lesssim \left( |\nabla \u|_{\b L^4}^2 + |\nabla A|_{L^2} |A|_{\b L^2}  \right) |\nabla \w|_{\b L^2}
	|\w|_{\b L^2},
\end{align*}
which can be recast into
\[
\partial_t |\w|^2_{\b L^2} +  \frac{\alpha}{2} |\nabla  \w|^2_{\b L^2} 
\lesssim 
\left( |\nabla \u|_{\b L^4}^4 + |\nabla A|_{\b L^2}^2 |A|_{\b L^2}^2  \right) | \w|_{\b L^2}^2.
\]
As the first factor on the right is integrable in time, Gronwall implies $\w=0$ on a small time interval, concluding the proof of this lemma.

\section{Concluding remarks}\label{Section: Concluding remarks}
With Theorem \ref{Theorem: Struwe u} and \ref{Theorem: u unique}, the inverse transformation yields our main result on solutions $ m $ of the stochastic LLG equation \eqref{eq: sLLG}. 

Let us briefly comment on the regularity of noise required in the proofs in Sections \ref{Section: Proof existence} and \ref{Section: Proof uniqueness}. 
For the existence proof, we set $ \sigma = 4 $ to ensure that $ A \in C([0,T]; \b H^2) \cap L^2(0,T; \b H^3) $, $ \P $-a.s. and thereby use Proposition \ref{Prop: u H2 est global}. 
For the uniqueness proof, $ \sigma = 2 $ is sufficient. 
As in Remark \ref{Remark: Coulomb gauge}, an improved control of $ \nabla \cdot A $ would lower the regularity requirements. 
Considered as time-space forms over $ \R \times \T^2 $, they satisfy the caloric gauge $ A_0=0 $, but roughness in time limits the class of admissible gauge transformations, which is a crucial aspect from the point of view of PDE analysis. 
Although we do not aim for an optimal regularity result here, it is worth mentioning some heuristic limitations from the basic energy law and scaling property.
Starting from a strong (in PDE sense) solution $ \m $ of \eqref{eq: sLLG} in the simplest case of one-dimensional Wiener process,  
\[
d\m = -\m \times \Delta \m \ dt + \alpha (\Delta \m +|\nabla \m|^2\m) \ dt
+ (\m \times g) \circ dW. 
\]
with $g = f h$ for some $ f: \T^2 \to \R $ and constant $ h \in \R^3 $, it follows from It\^o's lemma that 
\begin{align*}
	&E(\m(t)) - E(\m(0)) + \alpha \int_0^t \left| \Delta \m + |\nabla \m|^2 \m  \right|_{\b L^2}^2(s) \ ds \\
	&= \int_0^t \langle \nabla \m, \m \times \nabla g \rangle_{\b L^2}(s) \ dW(s) + \frac{1}{2} \int_0^t | \m \times \nabla g|_{\b L^2}^2(s) \ ds. 
\end{align*}
This suggests $\sigma=1$ as borderline regularity for weak solutions based on energy concepts. This situation has been considered in \cite{Hocquet1} in the case of HMHF. 
On the other hand, 
for every $\lambda>0$, the process $ W_\lambda $ given by  $W_\lambda(t)=\lambda^{-1}W(t\lambda^2)$ is again a cylindrical Wiener process on $ \ell^2 $. 
Starting from a strong (in PDE sense) solution of \eqref{eq: sLLG} on $[t_0,t_0+r^2] \times B_r(x_0)$, 
let
\begin{align*}
	\m_{\lambda}(t,x)&=\m(t_0 + \lambda^2 t, x_0+ \lambda x), \\
	G_\lambda(\varphi) l &= \sum_{j=1}^\infty \lb l, h_j \rb \varphi \times g_{j,\lambda},
\end{align*}
where $ g_{j,\lambda}(x)=g_j(x_0+\lambda x) $, $ j \in \b N $, for any $\lambda >0$. 
Then
\[
d\m_\lambda 
= - \m_\lambda \times \Delta \m_\lambda \ dt 
+ \alpha (\Delta \m_\lambda + |\nabla \m_\lambda|^2 \m_\lambda) \ dt 
+ \lambda G_\lambda(\m_\lambda) \circ dW_\lambda. 
\]
Notably, in dimension 2, dilation is norm preserving in $L^\infty \cap H^1$. In the blow-up regime $\lambda<1$, the contributions from noise are uniformly controlled. 
This suggests a borderline regularity 
$\sigma=0$ where the noise is equally active on every scale.

\end{document}